\documentclass[final]{mcom-l}

\usepackage{url}
\usepackage{color}
\usepackage[colorlinks=true,linkcolor=blue,citecolor=blue,urlcolor=blue]{hyperref}
\usepackage{cleveref}
\usepackage{enumitem}
\usepackage{booktabs}
\usepackage{multirow}
\usepackage{graphicx}
\usepackage{subfigure}
\usepackage{bbding}

\DeclareMathOperator{\R}{\mathbb{R}}
\DeclareMathOperator{\N}{\mathbb{N}}

\DeclareMathOperator{\setC}{\mathcal{C}}
\DeclareMathOperator{\setA}{\mathcal{A}}
\DeclareMathOperator{\setQC}{\mathcal{QC}}

\DeclareMathOperator{\argmin}{argmin}

\DeclareMathOperator{\PD}{PD}

\DeclareMathOperator{\SC}{SC}

\newcommand{\EiCP}[1]{\text{EiCP}(#1)}
\newcommand{\AEiCP}[1]{\text{AEiCP}(#1)}
\newcommand{\SEiCP}[1]{\text{SEiCP}(#1)}

\newcommand{\IPOPT}{\texttt{IPOPT}}
\newcommand{\KNITRO}{\texttt{KNITRO}}
\newcommand{\FILTERSD}{\texttt{FILTERSD}}

\newcommand{\cmark}{\Checkmark}
\newcommand{\xmark}{\XSolidBrush}
\newcommand{\RAND}[1]{\texttt{RAND(#1)}}
\newcommand{\MaxIT}{\texttt{MaxIT}}

\usepackage{algorithm,algorithmic}

\newenvironment{proc}[2]
{\begin{algorithm}[#1]\floatname{algorithm}{#2:}}
	{\end{algorithm}\addtocounter{algorithm}{-1}}

\newtheorem{theorem}{Theorem}[section]
\newtheorem{lemma}[theorem]{Lemma}

\newtheorem{proposition}[theorem]{Proposition}

\theoremstyle{definition}

\newtheorem{hypothesis}[theorem]{Hypothesis}

\theoremstyle{remark}
\newtheorem{remark}[theorem]{Remark}

\numberwithin{equation}{section}

\graphicspath{{figs/}}

\begin{document}

\title[Accelerated DCA for AEiCP]{Accelerated DC Algorithms for the Asymmetric Eigenvalue Complementarity Problem}


\author{Yi-Shuai Niu}
\address{Department of Combinatorics and Optimization, Faculty of Mathematics, University of Waterloo, Canada }
\email{niuyishuai82@hotmail.com}
\thanks{The author was supported by the Natural Science Foundation of China (Grant No: 11601327).}

\subjclass[2020]{Primary 65F15, 90C33, 90C30; Secondary 90C26, 93B60, 47A75}

\date{}


\keywords{Accelerated DC algorithms, asymmetric eigenvalue complementarity problem, difference-of-convex-sums-of-square decomposition}

\begin{abstract}
	We are interested in solving the Asymmetric Eigenvalue Complementarity Problem (AEiCP) by accelerated Difference-of-Convex (DC) algorithms. Two novel hybrid accelerated DCA: the Hybrid DCA with Line search and Inertial force (HDCA-LI) and the Hybrid DCA with Nesterov's extrapolation and Inertial force (HDCA-NI), are established. We proposed three DC programming formulations of AEiCP based on Difference-of-Convex-Sums-of-Squares (DC-SOS) decomposition techniques, and applied the classical DCA and 6 accelerated variants (BDCA with exact and inexact line search, ADCA, InDCA, HDCA-LI and HDCA-NI) to the three DC formulations. Numerical simulations of 7 DCA-type methods against state-of-the-art optimization solvers IPOPT, KNITRO and FILTERSD, are reported.
\end{abstract}

\maketitle

\section{Introduction} \label{sec:Intro}
Asymmetric Eigenvalue Complementarity Problem (AEiCP) consists of finding complementary eigenvectors $x\in \R^n\setminus\{0\}$ and complementary eigenvalues $\lambda\in \R$ such that 
\begin{equation}\label{eq:aeicp}
	\begin{cases}
		w = \lambda B   x - A  x,\\
		x^{\top}  w = 0,\\
		0\neq x\geq 0, w\geq 0.
	\end{cases}\tag{AEiCP}
\end{equation}
where $x^{\top}$ is the transpose of $x$, $A\in \R^{n\times n}$ is an asymmetric real matrix, and $B$ is a \emph{positive definite} ($\PD$) matrix (unnecessarily to be symmetric). If $A$ and $B$ are both symmetric, then the problem is called \emph{symmetric EiCP} (SEiCP). Throughout the paper, we will denote the solution set of problem \eqref{eq:aeicp} (resp. SEiCP) by $\AEiCP{A,B}$ (resp. $\SEiCP{A,B}$). The Eigenvalue Complementarity Problem appeared in the study of static equilibrium states of mechanical systems with unilateral friction \cite{Costa01}, and found interest in the Spectral Theory of Graphs \cite{fernandes2017complementary,seeger2018complementarity}.

In \cite[Theorem 2.1]{Judice09}, it is known that \eqref{eq:aeicp} can be expressed as a finite variational inequality on the unit simplex $\Omega := \{x\in \R^n: e^{\top} x = 1, x\geq 0\}$, which involves finding $\bar{x}\in \Omega$ satisfying:
$$\left\langle \left( \frac{x^{\top} A x}{x^{\top} B x}B - A \right) x , x - \bar{x}\right\rangle\geq 0, \quad \forall x\in \Omega.$$
If $B$ is a \emph{strictly copositive} ($\SC$) matrix (i.e., $x^{\top} B x >0$ for all $x\geq 0$ and $x\neq 0$), then the variational inequality is guaranteed to have a solution \cite{facchinei2003finite}. 
Hence, \eqref{eq:aeicp} has a solution since $B\in \PD$ (a special case of $B\in \SC$). Moreover, \eqref{eq:aeicp} has at most $n2^{n-1}$ distinct $\lambda$-solutions \cite[Proposition 3]{Queiroz03}, and it has a positive complementary eigenvalue if $A$ is a copositive matrix (i.e., $x^{\top}  A   x \geq  0, \forall x \geq 0$) and $-A$ is a $R_0$ matrix (i.e., $M\in R_0$ if and only if $[x\geq 0, M  x\geq 0, x^{\top} Mx =0]\Rightarrow x=0$) \cite[Theorem 2.2]{Judice09}. In particular, when both $A,B\in \PD$, then all complementary eigenvalues of \eqref{eq:aeicp} are positive. If $A$ is not $\PD$, then we can find some $\mu>0$ such that $A+\mu B\in\PD$ and, as shown in \cite[Theorem 2.1]{niu2023accelerated}, we have: $$(x,\lambda)\in \EiCP{A,B} \Leftrightarrow (x, \lambda + \mu)\in \EiCP{A+\mu B,B}.$$ 
Without loss of generality, throughout the paper, we assume that:	
\begin{hypothesis}\label{hyp1}
	Both $A$ and $B$ are PD matrices for \eqref{eq:aeicp}.
\end{hypothesis}

In fact, solving the SEiCP is much easier than the AEiCP. For the SEiCP, there are three well-known equivalent formulations (see e.g., \cite{niu2023accelerated}), namely the Rayleigh quotient formulation
\begin{equation}\label{prob:seicp_frac}
	\max\left\{\frac{x^{\top}   A   x}{x^{\top}   B   x}: x\in \Omega \right\}, \tag{RP}
\end{equation}
the logarithmic formulation
\begin{equation}\label{prob:seicp_log}
	\max\left\{ \ln(x^{\top}   A   x) - \ln(x^{\top}   B   x): x\in \Omega \right\}, \tag{LnP}
\end{equation}
and the quadratic formulation
\begin{equation}\label{prob:seicp_qcqp}
	\max\{x^{\top}   A   x:  x^{\top}   B   x\leq 1, x\geq 0 \}. \tag{QP}
\end{equation}
These formulations are equivalent to SEiCP in the sense that if $\bar{x}$ is a stationary point of \eqref{prob:seicp_qcqp}, \eqref{prob:seicp_frac} or \eqref{prob:seicp_log}, then  $$(\bar{x},(\bar{x}^{\top}A\bar{x})/(\bar{x}^{\top}B\bar{x})) \in \SEiCP{A,B}.$$
However, they are no longer equivalent to \eqref{eq:aeicp}. A simple convincing counterexample 
is given below: Let $n=2$, $B=I_2$ and $$A = \begin{bmatrix}
	-1 & 1\\
	-2 & 2 
\end{bmatrix}.$$
Then, the vector $\bar{x}=[0,1]^{\top}$ is a stationary point of \eqref{prob:seicp_qcqp}, \eqref{prob:seicp_frac} and \eqref{prob:seicp_log}, but $$(\bar{x},(\bar{x}^{\top}A\bar{x})/(\bar{x}^{\top}B\bar{x})) = ([0,1]^{\top},2)\notin \AEiCP{A,B}$$ because 
$$\bar{w} = (\bar{x}^{\top}A\bar{x})/(\bar{x}^{\top}B\bar{x}) B \bar{x} - A \bar{x} = \begin{bmatrix}
	-1\\ 0
\end{bmatrix}\leq 0.$$
In fact $([0,1]^{\top},2)\in \SEiCP{A+A^{\top},B+B^{\top}}$, but $\SEiCP{A+A^{\top},B+B^{\top}}$ and $\AEiCP{A,B}$ are not equivalent for \eqref{eq:aeicp}. Hence, we can not get a solution of \eqref{eq:aeicp} by searching a stationary point of \eqref{prob:seicp_qcqp}, \eqref{prob:seicp_frac} or \eqref{prob:seicp_log}.

Concerning solution algorithms for \eqref{eq:aeicp}, a number of algorithms have been designed such as Enumerative method \cite{Judice09,Fernandes14}, Semismooth Newton method \cite{Adly11}, Descent algorithm for minimizing the regularized gap function \cite{Bras12}, Difference-of-Convex (DC) algorithm \cite{Niu12,Niu15,niu2019improved,niu2023accelerated}, Splitting algorithm \cite{iusem2019splitting}, Active-set method \cite{Bras17}, Projected-Gradient algorithm \cite{Judice08}, and Alternating Direction Method of Multipliers (ADMM) \cite{judice2021alternating}. Most of these methods can be applied to both AEiCP and SEiCP. However, despite the theoretical ability of the Enumerative method to always find a solution to \eqref{eq:aeicp}, its practical application often proves computationally demanding or even intractable for larger instances within reasonable CPU time. The other algorithms face the same drawback, that is they attempt to solve nonconvex optimization formulations using local optimization approaches, which may fail for many instances in \eqref{eq:aeicp}, particularly when solution accuracy is on demand and at least one of the matrices $A$ and $B$ is ill-conditioned. 

In this paper, we will focus on solving three nonlinear programming (NLP) formulations of \eqref{eq:aeicp} via several accelerated DCA. First, we present in \Cref{sec:accDCAs} the classical DCA and three accelerated variants of DCA (BDCA \cite{BDCA_S,Niu19higher}, ADCA \cite{phan2018accelerated} and InDCA \cite{de2019inertial,you2022refined}) for solving the convex constrained DC programming problem with polynomial objective function over a closed convex set. Then, the spotlight is on establishing two novel hybrid accelerated DCA, namely the \emph{Hybrid DCA with Line search and Inertial force} (HDCA-LI) and the \emph{Hybrid DCA with Nesterov's extrapolation and Inertial force} (HDCA-NI). A rigorous convergence analysis for both HDCAs is furnished. In \Cref{sec:DCformulations}, we delve into establishing three DC programming formulations for \eqref{eq:aeicp} by leveraging the DC-SOS decomposition techniques \cite{niu2018difference}, acclaimed for generating superior quality DC decompositions of polynomials. Following this, we discuss in \Cref{sec:BDCAforAEiCP,sec:ADCA&InDCAforAEiCP} the application of the aforementioned accelerated variants of DCA in solving the three DC formulations of \eqref{eq:aeicp}. Finally, in \Cref{sec:Simulations}, we report numerical simulations comparing the proposed $7$ DCA-type algorithms with cutting-edge NLP solvers, namely IPOPT, KNITRO, and FILTERSD, and demonstrate their numerical performance on some large-scale and ill-conditioned instances.

\section{Accelerated DC Algorithms} \label{sec:accDCAs}
Consider the convex constrained DC (Difference-of-Convex) program defined by 
\begin{equation}
	\label{prob:P}
	\min_{x\in \setC} \{f(x):=g(x) - h(x)\}, \tag{P}
\end{equation}
where $\setC$ is a nonempty closed and convex subset of a finite dimensional Euclidean space endowed with an inner product $\langle \cdot , \cdot \rangle$ and an induced norm $\|\cdot\|$. Both $g$ and $h$ are convex polynomials such that $g$ (resp. $h$) is $\rho_g$-convex with $\rho_g\geq 0$ (resp. $\rho_h$-convex with $\rho_g\geq 0$) over $\setC$, i.e., $g(\cdot)-\frac{\rho_g}{2}\|\cdot\|^2$ (resp. $h(\cdot)-\frac{\rho_h}{2}\|\cdot\|^2$) is convex over $\setC$. If $\rho_g>0$ (resp. $\rho_h>0$), then $g$ (resp. $h$) is called $\rho_g$ (resp. $\rho_h$) strongly convex over $\setC$.

In this section, we briefly introduce the standard DCA and three existing accelerated variants of DCA, namely BDCA, ADCA, and InDCA, for solving \eqref{prob:P}. Following that, we propose two novel hybrid accelerated DCAs: the \emph{Hybrid DCA with Line search and Inertial force} (HDCA-LI) and the \emph{Hybrid DCA with Nesterov's extrapolation and Inertial force} (HDCA-NI). We also provide their convergence analysis. Some characteristics of these DCA-type algorithms are summarized in Table \ref{tab:alldcalgos}, where the columns \verb|Linesearch|, \verb|Inertial|, and \verb|Nesterov| indicate whether the algorithm requires a line search, heavy-ball inertial force, and Nesterov's extrapolation, respectively. The \verb|monotone| column represents the monotonicity of the generated sequence $\{f(x^k)\}$.
\begin{table}[tbhp]
	\caption{Some characteristics of the DCA-type algorithms.}
	\label{tab:alldcalgos}
	\centering
	\begin{tabular}{l|cccc} \toprule
		Algorithm & Linesearch & Inertial & Nesterov & Monotone \\
		\midrule
		DCA \cite{Pham97} & \xmark & \xmark & \xmark & \cmark\\
		BDCA \cite{Niu19higher} & \cmark & \xmark & \xmark & \cmark\\
		ADCA \cite{phan2018accelerated} & \xmark & \xmark & \cmark & \xmark\\
		InDCA \cite{de2019inertial,you2022refined} & \xmark & \cmark & \xmark & \xmark\\
		HDCA-LI  & \cmark & \cmark & \xmark & \xmark\\
		HDCA-NI  & \xmark & \cmark & \cmark & \xmark\\
		\bottomrule
	\end{tabular}
\end{table}

\subsection{DCA}\label{subsec:DCA}
One of the most renowned algorithm for solving \eqref{prob:P} is called DCA, which is first introduced by Pham Dinh Tao in 1985 as an extension of the
subgradient method \cite{Pham1986algorithms}, and extensively developed by 
Le Thi Hoai An and Pham Dinh Tao since 1994 (see \cite{Pham97,Pham98,Pham05,DCA30} and the references therein).

DCA applied to \eqref{prob:P} consists of constructing a sequence $\{x^k\}$ by solving the convex subproblems:
\begin{equation}\label{alg:DCA}
	\boxed{x^{k+1}\in \argmin \{ g(x) - \langle x, \nabla h(x^k) \rangle : x\in \setC \},} \tag{DCA}
\end{equation}
that is, minimizing the convex majorant (cf. surrogate) of the DC function $f$ at iterate $x^k$ in form of $g(x) - [h(x^k) + \langle \nabla h(x^k),x-x^k \rangle]$ by linearizing $h$ at $x^k$.

DCA is a descent method enjoys the next convergence properties:
\begin{theorem}[Convergence theorem of DCA for \eqref{prob:P}, see e.g., \cite{Pham97,niu2022convergence,Niu19higher}]\label{thm:convDCA}
	Let $\{x^k\}$ be a well-defined sequence generated by DCA for problem \eqref{prob:P} from an initial point $x^0\in \R^n$. If $\inf\{f(x): x\in\setC\} > -\infty$ and the sequence $\{x^k\}$ is bounded, then
	\begin{itemize}[leftmargin=12pt]
		\item[$\bullet$] \textbf{(monotonicity and convergence of $\{f(x^k)\}$)} the sequence $\{f(x^k)\}$ is non-increasing and convergent.
		\item \textbf{(sufficiently descent property)}
		$$f(x^k) - f(x^{k+1}) \geq \frac{\rho_g+\rho_h}{2}\|x^k - x^{k+1}\|^2, \quad \forall k\geq 1.$$
		\item[$\bullet$] \textbf{(square summable property)} if $\rho_g+\rho_h>0$, then 
		$$\sum_{k\geq 0} \|x^{k+1}-x^k\|^2 < \infty.$$
		\item \textbf{(subsequential convergence of $\{x^k\}$)} any cluster point $x^*$ of the sequence $\{x^k\}$ is a critical point of \eqref{prob:P}, i.e., $\nabla h(x^*) \in \partial (g+\chi_{\setC})(x^*)$.
		\item \textbf{(convergence of $\{x^k\}$)} if $f$ is a KL function, $\rho_g+\rho_h>0$, $h$ has locally Lipschitz continuous gradient over $\setC$, and $\setC$ is a semi-algebraic set, then the sequence $\{x^k\}$ is convergent.
	\end{itemize}
\end{theorem}
\begin{remark}\label{rmk:DCA} 
	\begin{itemize}[leftmargin=12pt]
		\item DCA is a monotone algorithm in the sequence $\{f(x^k)\}$.
		\item The symbol $\partial$ denotes the subdifferential of the convex function \cite{rockafellar1970} and $\chi_{\setC}$ is the indicator function defined as $\chi_{\setC}(x) = 0$ if $x\in \setC$ and $\chi_{\setC}(x) = \infty$ otherwise.
		\item $\inf\{f(x): x\in\setC\} > -\infty$ does not necessarily imply that the sequence $\{x^k\}$ is bounded. For instance, consider $g(x)=e^x$, $h(x)=0$ and $\setC = \R$. Clearly, $\inf\{f(x): x\in\setC\} =0 > -\infty$. Then, DCA starting from any initial point $x^0\in \R$ will generate a sequence $\{x^k\}$ that tends to $\infty.$ Conversely, it is obvious that the existence of a bounded sequence $\{x^k\}$ does not imply that $\inf_{x\in \setC}f(x) > -\infty$.
		\item Without the assumption $\inf_{x\in \setC}f(x) > -\infty$, then the convergence of $\{f(x^k)\}$ may not be true; without the boundedness assumptions of $\{x^k\}$, then a cluster point of $\{x^k\}$ may not exist; without the assumption $\rho_g+\rho_h>0$, then $\|x^{k+1}-x^k\|\to 0$ may not be true, and thus a limit point of $\{x^k\}$ may not be a critical point of \eqref{prob:P} (see the corresponding counterexamples in \cite{niu2022convergence}). 
		\item The strong convexity assumption in $g$ or $h$ is not restrictive and can be easily guaranteed in practice, since we can introduce the regularity term $\frac{\rho}{2}\|x\|^2$ (for any $\rho>0$) into $g$ and $h$ for getting a satisfactory DC decomposition as $$\underbrace{(g(x) + \frac{\rho}{2}\|x\|^2)}_{=\tilde{g}(x)} - \underbrace{(h(x) + \frac{\rho}{2}\|x\|^2)}_{=\tilde{h}(x)},$$
		where both $\tilde{g}$ and $\tilde{h}$ are $\rho$-strongly convex. 
		\item The KL (Kurdyka-{\L}ajasiewicz) assumption plays an important role to guarantee the convergence of the whole sequence $\{x^k\}$, which is an immediate consequence of \cite[Theorem 5]{niu2022convergence} (see also \cite{le2018convergence}). The KL function is a function satisfying the well-known KL property (see e.g., \cite[Definition 3]{niu2022convergence}), which is ubiquitous in optimization and its applications, e.g., the semialgebraic, subanalytic, log and exp are KL functions (see \cite{kurdyka1998gradients,bolte2007lojasiewicz,attouch2013convergence} and the
		references therein). In particular, for a closed convex semi-algebraic set $\setC$ (i.e., $\setC$ is a set of polynomial equations and inequalities) and the objective function $f$ has a polynomial DC decomposition $g-h$ where $g$ and $h$ are $\rho_g$- and $\rho_h$-convex polynomials over $\setC$ with $\rho_g+\rho_h>0$, then the convergence of $\{x^k\}$ is guaranteed.
	\end{itemize}
\end{remark}

\subsection{BDCA}\label{subsec:BDCA}
BDCA is an accelerated DCA by introducing the (exact or inexact) line search into DCA for acceleration. It was first proposed by Artacho et al. in 2018 \cite{BDCA_S} for unconstrained smooth DC program with the Armijo line search, then extended by Niu et al. in 2019 \cite{Niu19higher} for general convex constrained smooth and nonsmooth DC programs with both exact and inexact line searches.  

The general idea of BDCA in \cite{Niu19higher} is to introduce a line search along the \emph{DC descent direction} (a feasible and descent direction generated by two consecutive iterates of DCA)  $d^k:= z^{k}-x^{k}$ to find a better candidate $x^{k+1}$, where $$z^k\in\argmin \{g(x) - \langle y^k,x \rangle\},\quad y^k\in \partial h(x^k).$$ The line search can be performed either exactly ($=$) or inexactly ($\approx$) to find an optimal solution or an approximate solution of the line search problem 
$$\alpha_k = \text{or} \approx \argmin \{f(z^k + \alpha d^k): z^k + \alpha d^k \in \setC, \alpha \geq 0\}$$
with $f(z^k + \alpha d^k)\leq f(z^k)$.
Then we update $x^{k+1}$ by
$$x^{k+1} = x^k + \alpha_k d^k.$$

BDCA for problem \eqref{prob:P} is summarized in \Cref{alg:BDCA}. 
\begin{algorithm}[ht!]
	\caption{BDCA for \eqref{prob:P}}
	\label{alg:BDCA}
	\begin{algorithmic}[1]
		\REQUIRE $x^0\in \R^n$, $\bar{\alpha}>0$;
		\FOR{$k=0,1,\ldots$}
		\STATE $z^{k}\in \argmin\{g(x)-\langle x, \nabla h(x^k) \rangle :x\in \setC\}$;\label{codeline:bdca_P_convexsubprob}
		\STATE initialize $x^{k+1}\leftarrow z^k$;
		\STATE $d^k\leftarrow z^{k}-x^{k}$; \label{codeline:bdca_P_dk}
		\IF{$\setA(z^{k})\subset \setA(x^{k})$ and $\langle \nabla f(z^k), d^k\rangle < 0$}\label{codeline:bdca_P_linesearchcond}
		\STATE $x^{k+1}\leftarrow \text{LineSearch}(z^k,d^k,\bar{\alpha})$; \label{codeline:bdca_P_linesearch}
		\ENDIF\label{codeline:bdca_P_endif}
		\ENDFOR
	\end{algorithmic}
\end{algorithm}

Some comments regarding \Cref{alg:BDCA} include:
\begin{itemize}[leftmargin=12pt]
	\item DCA for \eqref{prob:P} is just BDCA without the codes between line \ref{codeline:bdca_P_dk} to \ref{codeline:bdca_P_endif}.
	\item For proceeding the line search, we have to check the conditions $\setA(z^k)\subset \setA(x^k)$ and $\langle \nabla f(z^k), d^k\rangle < 0$ in line \ref{codeline:bdca_P_linesearchcond}, where $\setA(x^k)$ is the active set of $\setC$ at $x^k$. These are sufficient conditions for a DC descent direction, but not necessary when $\setC$ is not a polyhedral set. In fact, it was proven in \cite{Niu19higher} that $d^k$ is a `potentially' descent direction of $f$ at $z^{k}$ such that 
	\begin{itemize}[leftmargin=12pt]
		\item if $f'(z^{k};d^k)<0$, then $d^k$ is a descent direction of $f$ at $z^k$;
		\item if $\setC$ is polyhedral, then $\setA(z^{k})\subset \setA(x^{k})$ is a necessary and sufficient condition for $d^k$ to be a feasible direction at $z^k$;
		\item if $\setC$ is not polyhedral, then $\setA(z^{k})\subset \setA(x^{k})$ is just a necessary (but not always sufficient) condition for $d^k$ to be a feasible direction at $z^k$.
	\end{itemize}
	\item 
	The line search procedure \verb|LineSearch|($z^k$,$d^k$,$\bar{\alpha}$) in line \ref{codeline:bdca_P_linesearch} can be either exact or inexact, depending on the problem size and structure. Generally speaking, finding an exact solution for $\alpha_k$ is computationally expensive or intractable when the problem size is too large or the objective function $f$ is too complicated. In this case, we suggest finding an approximate solution for $\alpha_k$ using an inexact line search procedure such as Armijo-type, Goldstein-type, or Wolfe-type \cite[Chapter 3]{wright2006numerical}. The third parameter $\bar{\alpha}$ denotes a given upper bound for the line search stepsize such that $\alpha_k\in [0,\bar{\alpha}]$ verifying $f(x^{k+1})\leq f(z^k)$ with
	$$x^{k+1} = z^k + \alpha_k d^k, \quad k=0,1,2, \ldots.$$ Note that for exact line search with unbounded set $\setC$, if without $\bar{\alpha}$ (i.e., $\bar{\alpha}=\infty$), then the sequence $\{\alpha_k\}$ may be unbounded. Consequently, we impose $\bar{\alpha}$ to ensure the boundedness of the sequence $\{\alpha_k\}$, which is crucial for the well-definedness of the sequence $\{x^{k}\}$ and the convergence of the algorithm BDCA. See successful examples of BDCA with inexact Armijo-type line search in \cite{Niu19higher} and with exact line search in \cite{zhang2022boosted} for higher-order moment portfolio optimization problem, as well as BDCA with exact line search in \cite{niu2023accelerated} for symmetric eigenvalue complementarity problem.
\end{itemize}
BDCA enjoys the next convergence theorem:
\begin{theorem}[Convergence theorem of BDCA \Cref{alg:BDCA}, see \cite{niu2022convergence,Niu19higher}]\label{thm:convBDCA}
	Let $\{x^k\}$ be a well-defined and bounded sequence generated by BDCA \Cref{alg:BDCA} for problem \eqref{prob:P} from an initial point $x^0\in \R^n$. If $\rho_g+\rho_h>0$ and the solution set of \eqref{prob:P} is non-empty, then 
	\begin{itemize}[leftmargin=12pt]
		\item \textbf{(monotonicity and convergence of $\{f(x^k)\}$)} the sequence $\{f(x^k)\}$ is non-increasing and convergent.
		\item \textbf{(convergence of $\{\|x^k-z^{k}\|\}$ and  $\{\|x^k-x^{k-1}\|\}$)} $$\|x^k-z^{k}\|\xrightarrow{k\to \infty} 0 \quad \text{ and } \quad \|x^k-x^{k-1}\|\xrightarrow{k\to \infty} 0.$$
		\item \textbf{(subsequential convergence of $\{x^k\}$)} any cluster point of the sequence $\{x^k\}$ is a critical point of \eqref{prob:P}.
		\item \textbf{(convergence of $\{x^k\}$)} furthermore, if $f$ is a KL function, $\setC$ is a semi-algebraic set, and $h$ has locally Lipschitz continuous gradient over $\setC$, then the sequence $\{x^k\}$ is convergent.
	\end{itemize}
\end{theorem}
\begin{remark}
	BDCA, like the standard DCA, is a \emph{descent method}. \Cref{thm:convBDCA} is an immediate consequence of the general convergence theorem of BDCA for convex constrained nonsmooth DC programs established in \cite{Niu19higher}. \Cref{rmk:DCA} for standard DCA is adopted here as well.
\end{remark}

\subsection{ADCA}\label{subsec:ADCA}
ADCA is an accelerated DCA by introducing the Nesterov's acceleration into DCA, which is introduced by Phan et al. in 2018 \cite{phan2018accelerated}. 

The basic idea of ADCA is to compute $\nabla h(v^k)$ instead of $\nabla h(x^k)$ at iteration $k$, where $v^k$ is a more promising point than $x^k$ in the sense that $v^k$ is better than one of the last $q$ iterates $\{x^{k-q},\ldots,x^k\}$ in terms of objective function, i.e., 
\begin{equation}
	\label{eq:condadca}
	f(v^k)\leq \max_{t=\max\{0,k-q\},\ldots,k} f(x^t).
\end{equation}
The candidate $v^k$ is computed by Nesterov's extrapolation:
$$v^k = x^k + \frac{\theta_k -1}{\theta_{k+1}}\left(x^k - x^{k-1} \right),$$
where $$\theta_{0}=1, \quad \theta_{k+1} = \frac{1+\sqrt{1+4 \theta_k^2}}{2},\quad \forall k\geq 1.$$ 
If the condition \eqref{eq:condadca} is not satisfied, then we set $v^k=x^k$ as in DCA.

ADCA applied to problem \eqref{prob:P} is described below:
\begin{algorithm}[ht!]
	\caption{ADCA for \eqref{prob:P}}
	\label{alg:ADCA}
	\begin{algorithmic}[1]
		\REQUIRE $x^0\in \R^n$, $q\in \N$;
		\STATE \textbf{Initialization:} $x^{-1}=x^0$, $\theta_0=1$;
		\FOR{$k=0,1,\ldots$}
		\STATE $\theta_{k+1} \leftarrow \frac{1+\sqrt{1+4 \theta_k^2}}{2}$;
		\STATE $v^k \leftarrow x^k + \frac{\theta_k -1}{\theta_{k+1}}\left(x^k - x^{k-1} \right);$
		\IF{$f(v^k) > \max\{f(x^t) : t=\max\{0,k-q\},\ldots,k\}$}
		\STATE $v^k \leftarrow x^k;$
		\ENDIF
		\STATE $x^{k+1}\in \argmin\{g(x)-\langle x, \nabla h(v^k) \rangle :x\in \setC\}$;\label{codeline:adca_P_convexsubprob}
		\ENDFOR
	\end{algorithmic}
\end{algorithm}

\begin{remark}
	ADCA is a \emph{non-monotone} algorithm due to the introduction of Nesterov's extrapolation. It acts as a descent method when $q=0$. In this case, we choose a better candidate in terms of the objective value between $v^k$ and $x^k$ to compute $x^{k+1}$. If $q>0$, then ADCA can increase the objective function and consequently escape from a potential bad local minimum. A higher value of $q$ increases the chance of escaping bad local minima. 
\end{remark}
ADCA enjoys the next convergence theorem:
\begin{theorem}[Convergence theorem of ADCA \Cref{alg:ADCA}, see \cite{phan2018accelerated}]\label{thm:convADCA}
	Let $\{x^k\}$ be the sequence generated by ADCA \Cref{alg:ADCA} for problem \eqref{prob:P} from any initial point $x^0\in \R^n$ and with $q=0$. If $\rho_g+\rho_h>0$, $\inf\{f(x): x\in\setC\} > -\infty$ and the sequence $\{x^k\}$ is bounded, then any cluster point of $\{x^k\}$ is a critical point of \eqref{prob:P}.
\end{theorem}

\subsection{InDCA}\label{subsec:InDCA}
Inertial DCA (cf. InDCA) is an accelerated DCA by incorporating the momentum (heavy-ball type inertial force) into standard DCA. This method is first introduce by de Oliveira et al. in 2019 \cite{de2019inertial} and a refined version (RInDCA) is established by Niu et al. in \cite{you2022refined} with enlarged inertial stepsize for faster convergence. 

The basic idea of InDCA is to introduce an inertial force $\gamma (x^k - x^{k-1})$ (for some inertial stepsize $\gamma>0$) into $y^k\in \partial h(x^k)$ to get the next convex subproblem:
$$x^{k+1}\in \argmin \{g(x) - \langle y^k + \gamma (x^k - x^{k-1}),x \rangle : x\in \setC\}, \quad y^k\in \partial h(x^k).$$
InDCA applied to problem \eqref{prob:P} is described below:
\begin{algorithm}[ht!]
	\caption{InDCA for \eqref{prob:P}}
	\label{alg:InDCA}
	\begin{algorithmic}[1]
		\REQUIRE $x^0\in \R^n$;
		\STATE \textbf{Initialization:} $x^{-1}=x^0$, $\gamma\in [0,\frac{\rho_g + \rho_h}{2})$;
		\FOR{$k=0,1,\ldots$}
		\STATE $x^{k+1}\in \argmin\{g(x)-\langle \nabla h(x^k) + \gamma (x^k - x^{k-1}), x \rangle:x\in \setC\}$; \label{codeline:indca_P_convexsubprob}
		\ENDFOR
	\end{algorithmic}
\end{algorithm}

\begin{remark}
	Compared to the stepsize interval $[0,\rho_h/2)$ for $\gamma$ suggested in \cite{de2019inertial}, an enlarge stepsize interval $[0,(\rho_g + \rho_h)/2)$ was proven to be adequate for the convergence of InDCA \cite{you2022refined}. In practice, a larger stepsize leads to faster convergence.
\end{remark}
\begin{remark}
	Like ADCA, the sequence $\{f(x^k)\}$ generated by InDCA is not necessarily monotone. Numerical comparison between ADCA and InDCA has been performed on Image Denoising \cite{you2022refined} and Nonnegative Matrix Factorization \cite{phan2021dca}. Both experiments showed that ADCA outperformed InDCA with smaller stepsize in $[0,\rho_h/2)$. But InDCA with an enlarged stepsize can outperform ADCA (both speed and quality) on some instances of the Image Denoising dataset (see \cite{you2022refined}).
\end{remark}

The convergence of InDCA to problem \eqref{prob:P} is described below:

\begin{theorem}[Convergence theorem of InDCA \Cref{alg:InDCA}, see \cite{you2022refined,de2019inertial}]\label{thm:convInDCA}
	Let $\{x^k\}$ be the sequence generated by InDCA \Cref{alg:InDCA} for problem \eqref{prob:P} from any initial point $x^0\in \R^n$. If $\gamma \in [0,\frac{\rho_g + \rho_h}{2})$, $\rho_g+\rho_h>0$, $\inf\{f(x):x\in \setC\}>-\infty$ and the sequence $\{x^k\}$ is bounded, then 
	\begin{itemize}[leftmargin=12pt]
		\item \textbf{(sufficiently descent property)} $$\begin{aligned}
			f(x^{k+1}) + \frac{\rho_g + \rho_h - \gamma}{2} \|x^{k+1}-x^k\|^2 \leq& f(x^{k}) + \frac{\rho_g + \rho_h - \gamma}{2} \|x^{k}-x^{k-1}\|^2\\ 
			&- \frac{\rho_g + \rho_h - 2\gamma}{2} \|x^{k}-x^{k-1}\|^2.
		\end{aligned}$$
		\item \textbf{(convergence of $\{\|x^k-x^{k-1}\|\}$)} $ \|x^k-x^{k-1}\|\xrightarrow{k\to \infty} 0.$
		\item \textbf{(subsequential convergence of $\{x^k\}$)} any cluster point of the sequence $\{x^k\}$ is a critical point of \eqref{prob:P}.
		\item \textbf{(convergence of $\{x^k\}$)} furthermore, if $f$ is a KL function and $\setC$ is a semi-algebraic set, then the sequence $\{x^k\}$ is convergent.
	\end{itemize}
\end{theorem}

\subsection{HDCA}\label{subsec:HDCA}
We can combine the line search, the heavy-ball inertial force and the Nesterov's extrapolation with DCA to obtain some enhanced hybrid accelerated DCA algorithms. Here, we propose two hybrid methods: the \emph{Hybrid DCA with Line search and Inertial force} (HDCA-LI) and the \emph{Hybrid DCA with Nesterov's extrapolation and Inertial force} (HDCA-NI). 

The reason to do so is that: the inertial force intends to accelerate the gradient $\nabla h(x^k)$ by an inertial force $\gamma (x^k-x^{k-1})$; while both line search and Nesterov's extraplation play a similar rule by accelerating $x^k$ using a `potentially' better candidate in form of $x^k + \beta_k (x^k-x^{k-1})$ for some $\beta_k\geq 0$. Hence, by combining the inertial force with either line search or Nesterov's extrapolation, we may enhance the acceleration of DCA. 

\paragraph{\textbf{HDCA-LI}}
The Hybrid DCA with Line search and Inertial force accelerations (HDCA-LI) is described below:
\begin{algorithm}[ht!]
	\caption{HDCA-LI for \eqref{prob:P}}
	\label{alg:HDCA-LI}
	\begin{algorithmic}[1]
		\REQUIRE $x^0\in \R^n$, $\bar{\alpha}>0$;
		\STATE \textbf{Initialization:} $x^{-1}=x^0$, $\gamma\in [0,\frac{\rho_g + \rho_h}{1+(1+\bar{\alpha})^2})$; 
		\FOR{$k=0,1,\ldots$}
		\STATE $z^{k}\in \argmin\{g(x)-\langle \nabla h(x^k) + \gamma (x^k - x^{k-1}), x \rangle:x\in \setC\}$; \label{codeline:hdca-li_P_convexsubprob}
		\STATE $d^k\leftarrow z^{k}-x^{k}$;
		\IF{$\setA(z^{k})\subset \setA(x^{k})$ and $\langle \nabla f(z^k), d^k\rangle < 0$}
		\STATE $x^{k+1}\leftarrow \text{LineSearch}(z^k,d^k,\bar{\alpha})$; \ENDIF
		\ENDFOR
	\end{algorithmic}
\end{algorithm}

Some comments on HDCA-LI include:
\begin{itemize}[leftmargin=12pt]
	\item HDCA-LI is \emph{non-monotone} if $\gamma\neq 0$ and it reduces to BDCA if $\gamma=0$. 
	\item The upper bound of the inertial stepsize $\gamma$ is set as $(\rho_g + \rho_h)/(1+(1+\bar{\alpha})^2)$, differing from $\rho_h/2$ in InDCA \cite{de2019inertial} and $(\rho_g+\rho_h)/2$ in its refined version \cite{you2022refined}. The computation of this upper bound will be discussed in the convergence analysis of HDCA-LI (Theorem \ref{thm:convHDCA-LI}). This upper bound also reveals a trade-off between the inertial stepsize and the line search stepsize, that is a larger upper bound for the line search stepsize leads to a smaller upper bound for the inertial stepsize. 
\end{itemize}

Now, we establish the convergence theorem of HDCA-LI (similar to InDCA \cite{you2022refined} and BDCA \cite{Niu19higher}) based on the Lyapunov analysis. The reader is also referred to \cite{niu2022convergence} for a general framework to establish convergence analysis of DCA-type algorithm.
\begin{theorem}[Convergence theorem of HDCA-LI \Cref{alg:HDCA-LI}]\label{thm:convHDCA-LI}
	Let $\{x^k\}$ be the sequence generated by HDCA-LI \Cref{alg:HDCA-LI} for problem \eqref{prob:P} from any initial point $x^0\in \R^n$. Suppose that $\rho_g+\rho_h>0$, $\gamma \in [0,\frac{\rho_g + \rho_h}{1+(1+\bar{\alpha})^2})$, $\inf\{f(x):x\in \setC\}>-\infty$ and the sequence $\{x^k\}$ is bounded. Let $$E_k := f(x^{k}) + \frac{\rho_g + \rho_h - \gamma}{2(1+\bar{\alpha})^2}  \|x^{k}-x^{k-1}\|^2, \forall k=1,2,\ldots.$$ Then 
	\begin{itemize}[leftmargin=12pt]
		\item \textbf{(sufficiently descent property)} $$
		E_{k+1} \leq E_k - \frac{\rho_g + \rho_h - \gamma(1+(1+\bar{\alpha})^2)}{2(1+\bar{\alpha})^2} \|x^{k}-x^{k-1}\|^2, \forall k=1,2,\ldots.$$
		\item \textbf{(convergence of $\{\|x^k-x^{k-1}\|\}$ and $\{\|d^k\|\}$)} $$ \|x^k-x^{k-1}\|\to 0 \quad \text{and} \quad \|d^k\|\to 0 \text{ as }k\to \infty.$$
		\item \textbf{(subsequential convergence of $\{x^k\}$)} any cluster point $x^*$ of the sequence $\{x^k\}$ is a critical point of \eqref{prob:P} (i.e., $\nabla h(x^*) \in \nabla g(x^*) + N_{\setC}(x^*)$).
	\end{itemize}
\end{theorem}
\begin{proof}
	\textbf{(sufficiently descent property):} 
	By the first order optimality condition to the convex problem 
	$$z^{k}\in \argmin\{g(x)-\langle \nabla h(x^k) + \gamma (x^k - x^{k-1}), x \rangle:x\in \setC\},$$
	we get $$\nabla h(x^k)+\gamma(x^k-x^{k-1}) \in \partial (g + \chi_{\setC})(z^{k}) = \nabla g(z^{k}) + N_{\setC}(z^{k}),$$
	where $N_{\setC}(z^k)$ stands for the normal cone of $\setC$ at $z^k$.
	Hence,
	$$\nabla h(x^k)+\gamma(x^k-x^{k-1}) - \nabla g(z^{k}) \in N_{\setC}(z^{k})$$
	implies that
	\begin{equation}
		\label{eq:ineq01}
		\langle \nabla h(x^k)+\gamma(x^k-x^{k-1}) - \nabla g(z^{k}), x^k - z^k \rangle \leq 0.
	\end{equation}
	Then by the $\rho_g$-convexity of $g$, we get
	\begin{eqnarray*}
		g(x^k) &\geq& g(z^{k}) + \langle \nabla g(z^k),x^k-z^{k}\rangle+\frac{\rho_g}{2}\| z^k-x^k\|^2\\
		&\overset{\eqref{eq:ineq01}}{\geq}& g(z^{k}) + \langle \nabla h(x^k) + \gamma(x^k-x^{k-1}),x^k-z^{k}\rangle+\frac{\rho_g}{2}\| z^k-x^k\|^2,
	\end{eqnarray*}
	that is 
	\begin{equation}
		\label{eq:stronglycvxofg}
		g(x^k) \geq g(z^{k}) + \langle \nabla h(x^k) + \gamma(x^k-x^{k-1}),x^k-z^{k}\rangle+\frac{\rho_g}{2}\| z^k-x^k\|^2.
	\end{equation} 
	On the other hand, it follows from the $\rho_h$-convexity of $h$ that 
	\begin{equation}
		\label{eq:stronglycvxofh}
		h(z^k) \geq h(x^k) + \langle \nabla h(x^k),z^k-x^k \rangle +\frac{\rho_h}{2}\| z^k-x^k\|^2.
	\end{equation}
	Summing (\ref{eq:stronglycvxofg}) and (\ref{eq:stronglycvxofh}), we get  
	\begin{equation}
		\label{eq:ineq02}
		f(x^k) \geq f(z^k) + \gamma\langle x^k-x^{k-1},x^k-z^k\rangle + \frac{\rho_g + \rho_h}{2}\| z^k-x^k\|^2.
	\end{equation}
	By applying $\langle x^k-x^{k-1},x^k-z^k\rangle \geq -(\| x^k-x^{k-1}\|^2 + \| x^k-z^k\|^2)/{2}$ to \eqref{eq:ineq02},
	\begin{equation}
		\label{eq:ineq03}
		f(z^k) + \frac{\rho_g + \rho_h-\gamma}{2}\| z^k-x^k\|^2 
		\leq f(x^{k}) + \frac{\gamma}{2}\| x^k-x^{k-1}\|^2.
	\end{equation}
	The (exact or inexact) line search procedure ensures that $f(z^k)\geq f(x^{k+1})$ and   
	\begin{equation}
		\label{eq:updatexk1}
		x^{k+1} = z^k + \alpha_k (z^k-x^k), \forall k=1,2,\ldots,
	\end{equation}
	under the assumption that $\alpha_k\in[0,\bar{\alpha}]$. Then  
	$$\|x^{k+1}-x^k\| \overset{\eqref{eq:updatexk1}}{=} \|z^k + \alpha_k (z^k-x^k) - x^k\| = (1+\alpha_k)\| z^k-x^k\|\leq(1+\bar{\alpha})\| z^k-x^k\|.$$
	It follows from $\|x^{k+1}-x^k\| \leq (1+\bar{\alpha})\| z^k-x^k\|$, $f(z^k)\geq f(x^{k+1})$ and $\rho_g + \rho_h > \gamma$ that
	\begin{equation}
		\label{eq:ineq04}
		f(z^k) + \frac{\rho_g + \rho_h-\gamma}{2}\| z^k-x^k\|^2 \geq f(x^{k+1}) + \frac{\rho_g + \rho_h-\gamma}{2(1+\bar{\alpha})^2}\| x^{k+1}-x^k\|^2. 
	\end{equation}
	Therefore, we get from \eqref{eq:ineq03} and \eqref{eq:ineq04} that
	$$\begin{aligned}
		f(x^{k+1}) + \frac{\rho_g + \rho_h-\gamma}{2(1+\bar{\alpha})^2}\| x^{k+1}-x^k\|^2 \leq & f(x^{k}) + \frac{\gamma}{2}\| x^k-x^{k-1}\|^2\\
		=& f(x^{k}) + \frac{\rho_g + \rho_h-\gamma}{2(1+\bar{\alpha})^2}\| x^k-x^{k-1}\|^2 \\
		&- \frac{\rho_g + \rho_h - \gamma (1+(1+\bar{\alpha})^2)}{2(1+\bar{\alpha})^2}\| x^k-x^{k-1}\|^2.
	\end{aligned}$$
	Taking $E_k = f(x^{k}) + \frac{\rho_g + \rho_h - \gamma}{2(1+\bar{\alpha})^2}  \|x^{k}-x^{k-1}\|^2$. Then 
	\begin{equation}
		\label{eq:descentlemma}
		\boxed{E_{k+1} \leq E_k - \frac{\rho_g + \rho_h - \gamma(1+(1+\bar{\alpha})^2)}{2(1+\bar{\alpha})^2} \|x^{k}-x^{k-1}\|^2, \forall k=1,2,\ldots.}
	\end{equation}
	\textbf{(convergence of $\{\|x^k-x^{k+1}\|\}$):} For all $\gamma \in [0,\frac{\rho_g + \rho_h}{1+(1+\bar{\alpha})^2})$, we have 
	$$\frac{\rho_g + \rho_h - \gamma(1+(1+\bar{\alpha})^2)}{2(1+\bar{\alpha})^2} > 0.$$
	Then, it follows from \eqref{eq:descentlemma} that the sequence $\{E_k\}$ is non-increasing. The assumption $\inf\{f(x):x\in \setC\}>-\infty$ and $E_k = f(x^{k}) + \frac{\rho_g + \rho_h - \gamma}{2(1+\bar{\alpha})^2}  \|x^{k}-x^{k-1}\|^2 \geq f(x^{k})$ ensure that the sequence $\{E_k\}$ is lower bounded. Consequently, the non-increasing and lower bound of the sequence ${E_k}$ ensure its convergence. Let $E_k\to E^*$ as $k\to \infty$.    
	Summing \eqref{eq:descentlemma} for $k$ from $1$ to $\infty$, we get 
	$$\sum_{k=1}^{\infty}\|x^{k}-x^{k-1}\|^2 \leq \frac{2(1+\bar{\alpha})^2}{\rho_g + \rho_h - \gamma(1+(1+\bar{\alpha})^2)} (E^1 - E^*) < \infty.$$
	Therefore,  $$\boxed{\|x^{k}-x^{k-1}\| \xrightarrow{k\to \infty}0.}$$ 
	\textbf{(convergence of $\{\|d^k\|\}$):} It follows immediately from 
	$$0\leq \|d^k\| = \| z^k-x^k\|\leq (1+\alpha_k)\| z^k-x^k\| = \|x^{k+1}-x^k\| \xrightarrow{k\to \infty} 0$$
	that 
	$$\boxed{\|d^k\| \xrightarrow{k\to \infty} 0.}$$
	\textbf{(subsequential convergence of $\{x^k\}$):}
	The boundedness of the sequence $\{x^k\}$ implies that the set of its cluster points is non-empty, and for any cluster point $x^*$ there exists a convergent subsequence denoted by $\{x^{k_j}\}_{j\in \N}\subset \setC$ converging to $x^*$. The closedness of $\setC$ indicates that the limit point $x^*$ belongs to $\setC$. Then, we get from $x^{k_j}\to x^*$, $\|z^k-x^k\|\xrightarrow{k\to \infty} 0$ and $\|x^{k}-x^{k-1}\| \xrightarrow{k\to \infty}0$ that $$z^{k_j} \xrightarrow{j\to \infty} x^* \quad \text{ and }\quad x^{k_j-1} \xrightarrow{j\to \infty} x^*.$$ It follows from the first order optimality condition of the convex subproblem $$\nabla h(x^k)+\gamma(x^k-x^{k-1}) \in \nabla g(z^{k}) + N_{\setC}(z^{k}),$$
	the closedness of the graph of $\partial \chi_{\mathcal{C}} = N_{\setC}$, and the continuity of $\nabla g$ and $\nabla h$ that 
	$$\lim_{j\to\infty} [\nabla h(x^{k_j}) + \gamma (x^{k_j} - x^{k_j-1})] \in \nabla g(\lim_{j\to \infty}z^{k_j}) + N_{\setC}(\lim_{j\to \infty}z^{k_j}).$$
	That is,
	$$\boxed{\nabla h(x^*) \in \nabla g(x^*) + N_{\setC}(x^*).}$$
	Hence, any cluster point of the sequence $\{x^k\}$ is a critical point of \eqref{prob:P}.
\end{proof}
\begin{remark}
	The sequential convergence of the sequence $\{x^k\}$ can also be established under certain regularity conditions, notably the \emph{{\L}ojasiewicz subgradient inequality} or the \emph{Kurdyka-{\L}ojasiewicz (KL) property}. Gratifyingly, both of these conditions are naturally met for the DC program \eqref{prob:P} with polynomial convex components $g$ and $h$, as well as a semi-algebraic convex set $\setC$. For specific examples illustrating the establishment of the convergence of $\{x^k\}$ and the rate of convergence of both $\{f(x^k)\}$ and $\{x^k\}$, the reader is directed to \cite[Theorem 2, Lemma 3, Theorem 7, Theorem 8]{niu2022convergence}. Here, we admit theses results (i.e., the convergence of $\{x^k\}$ and the rate of convergence of $\{f(x^k)\}$ and $\{x^k\}$ under the KL property) and omit their proofs.
\end{remark}

\paragraph{\textbf{HDCA-NI}}
The \emph{Hybrid DCA with Nesterov's extrapolation and Inertial force accelerations} (HDCA-NI) is described in \Cref{alg:HDCA-NI}. 
\begin{algorithm}[ht!]
	\caption{HDCA-NI for \eqref{prob:P}}
	\label{alg:HDCA-NI}
	\begin{algorithmic}[1]
		\REQUIRE $x^0\in \R^n$, $q\in \N$, $\bar{\beta}\in (0,1)$;
		\STATE \textbf{Initialization:} $x^{-1}=x^0$, $\theta_0=1$, $\beta_{0}=0$, $\delta = (1-\bar{\beta}^2)(\rho_g + \rho_h)/4$; 
		\FOR{$k=0,1,\ldots$}
		\STATE $\theta_{k+1} \leftarrow (1+\sqrt{1+4 \theta_k^2})/2$;
		\STATE $\beta_{k} \leftarrow (\theta_k -1)/\theta_{k+1}$;
		\IF{$\beta_k>\bar{\beta}$}
		\STATE $\theta_k\leftarrow 1;$
		\STATE $\beta_{k}\leftarrow\bar{\beta};$
		\ENDIF
		\STATE $v^k \leftarrow x^k + \beta_k\left(x^k - x^{k-1} \right)$;
		\STATE $\gamma_k\in [0, ((\rho_g+\rho_h)(1-\beta_k^2) - 4 \delta)/(3-\beta_k^2)];$
		\IF{$v^k\notin \setC$ or $f(v^k) + \frac{\rho_g + \rho_h-\gamma_k}{4}\| x^k-x^{k-1}\|^2 > \max\{f(x^t) + \frac{\rho_g + \rho_h-\gamma_t}{4}\| x^{t} - x^{t-1}\|^2 : t=\max\{0,k-q\},\ldots,k\}$}
		\STATE $v^k \leftarrow x^k;$
		\ENDIF
		\STATE $x^{k+1}\in \argmin\{g(x)-\langle \nabla h(v^k) + \gamma_k (x^k - x^{k-1}), x \rangle:x\in \setC\}$; \label{codeline:hdca-ni_P_convexsubprob}		
		\ENDFOR
	\end{algorithmic}
\end{algorithm}
It should be noted that if $\beta_{k} > \bar{\beta}$, we reset $\theta_k = 1$ to ensure that the sequence $\{\beta_k\}_k\subset[0, \bar{\beta}]$ for any provided upper bound $\bar{\beta} < 1$ (typically chosen close to 1, such as $\bar{\beta} = 0.9$). The inertial stepsize $\gamma_{k}$ can be any value within $[0, \frac{(\rho_g+\rho_h)(1-\beta_k^2) - 4 \delta}{3-\beta_k^2}]$ and is suggested to be equal to $\frac{(\rho_g+\rho_h)(1-\beta_k^2) - 4 \delta}{3-\beta_k^2}$ for better inertial force acceleration. 

The convergence analysis of HDCA-NI is established in a manner analogous to that of HDCA-LI (as seen in \Cref{thm:convHDCA-LI}) and ADCA \cite{phan2018accelerated} as below:
\begin{theorem}[Convergence theorem of HDCA-NI \Cref{alg:HDCA-NI}]\label{thm:convHDCA-NI}
	Let $\{x^k\}$ be the sequence generated by HDCA-NI \Cref{alg:HDCA-NI} for problem \eqref{prob:P} from any initial point $x^0\in \R^n$. Let 
	$$\begin{cases}
		\phi(k) := \argmin\{\| x^{t}-x^{t-1}\|^2: t=k,\ldots,k+q\},\\
		c_k: = \frac{(\beta_k^2-3)\gamma_k + (1-\beta_k^2) (\rho_g + \rho_h)}{4},\\
		E_k := \max \{f(x^{t}) + \frac{\rho_g + \rho_h-\gamma_t}{4}\| x^{t} - x^{t-1}\|^2: t=\max\{0,k-q\},\ldots, k\},
	\end{cases}$$
	for all $k\in \N$. Suppose that $\rho_g+\rho_h>0$, $\bar{\beta}\in (0,1)$, $\delta = (1-\bar{\beta}^2)(\rho_g + \rho_h)/4>0$, $\inf\{f(x):x\in \setC\}>-\infty$ and the sequence $\{x^k\}$ is bounded. Then 
	\begin{itemize}[leftmargin=12pt]
		\item \textbf{(sufficiently descent property)} $$
		E_{k+1+q} \leq E_k - \delta \|x^{\phi(k)}-x^{\phi(k)-1}\|^2, \forall k=0,1,2,\ldots.$$
		\item \textbf{(convergence of $\{\|x^{\phi(k)}-x^{\phi(k)-1}\|\}$ and $\{\|v^{\phi(k)}-x^{\phi(k)-1}\|\}$)}
		$$ \|x^{\phi(k)}-x^{\phi(k)-1}\|\xrightarrow{k\to \infty} 0\quad \text{ and }\quad \|v^{\phi(k)}-x^{\phi(k)-1}\|\xrightarrow{k\to \infty} 0.$$
		\item \textbf{(subsequential convergence of $\{x^k\}$)} Let $q=0$. Then any cluster point $x^*$ of the sequence $\{x^k\}$ is a critical point of \eqref{prob:P} (i.e., $\nabla h(x^*) \in \nabla g(x^*) + N_{\setC}(x^*)$).
	\end{itemize}
\end{theorem}
\begin{proof}
	\textbf{(sufficiently descent property):} 
	By the first order optimality condition to the convex subproblem 
	$$x^{k+1}\in \argmin\{g(x)-\langle \nabla h(v^k) + \gamma_k (x^k - x^{k-1}), x \rangle:x\in \setC\},$$
	we get $$\nabla h(v^k)+\gamma_k(x^k-x^{k-1}) \in \nabla g(x^{k+1}) + N_{\setC}(x^{k+1}).$$
	Hence,
	$$\nabla h(v^k)+\gamma_k(x^k-x^{k-1}) - \nabla g(x^{k+1}) \in N_{\setC}(x^{k+1}),$$
	which implies that
	\begin{equation}
		\label{eq:ineq01-bis}
		\langle \nabla h(v^k)+\gamma_k(x^k-x^{k-1}) - \nabla g(x^{k+1}), v^k - x^{k+1} \rangle \leq 0,
	\end{equation}
	where $v^k,x^{k+1}\in \setC$ by their definitions.\\
	Then by the $\rho_g$-convexity of $g$ over $\setC$ and $v^k,x^{k+1}\in \setC$, we get
	\begin{eqnarray*}
		g(v^k) &\geq& g(x^{k+1}) + \langle \nabla g(x^{k+1}),v^k-x^{k+1}\rangle+\frac{\rho_g}{2}\| v^k-x^{k+1}\|^2\\
		&\overset{\eqref{eq:ineq01-bis}}{\geq}& g(x^{k+1}) + \langle \nabla h(v^k) + \gamma_k(x^k-x^{k-1}),v^k-x^{k+1}\rangle+\frac{\rho_g}{2}\| v^k-x^{k+1}\|^2,
	\end{eqnarray*}
	that is 
	\begin{equation}
		\label{eq:stronglycvxofg-bis}
		g(v^k) \geq g(x^{k+1}) + \langle \nabla h(v^k) + \gamma_k(x^k-x^{k-1}),v^k-x^{k+1}\rangle+\frac{\rho_g}{2}\| v^k-x^{k+1}\|^2.
	\end{equation} 
	On the other hand, it follows from the $\rho_h$-convexity of $h$ over $\setC$ that 
	\begin{equation}
		\label{eq:stronglycvxofh-bis}
		h(x^{k+1}) \geq h(v^k) + \langle \nabla h(v^k),x^{k+1}-v^k \rangle +\frac{\rho_h}{2}\| v^k-x^{k+1}\|^2.
	\end{equation}
	Summing (\ref{eq:stronglycvxofg-bis}) and (\ref{eq:stronglycvxofh-bis}), we get  
	\begin{equation}
		\label{eq:ineq02-bis}
		f(v^k) \geq f(x^{k+1}) + \gamma_k\langle x^k-x^{k-1},v^k-x^{k+1}\rangle + \frac{\rho_g + \rho_h}{2}\| v^k-x^{k+1}\|^2.
	\end{equation} 
	By applying $\langle x^k-x^{k-1},v^k-x^{k+1}\rangle \geq -(\| x^k-x^{k-1}\|^2 + \| v^k-x^{k+1}\|^2)/{2}$ to \eqref{eq:ineq02-bis}, we obtain
	\begin{equation}
		\label{eq:ineq03-bis}
		f(x^{k+1}) + \frac{\rho_g + \rho_h-\gamma_k}{2}\| v^k-x^{k+1}\|^2 
		\leq f(v^k) + \frac{\gamma_k}{2}\| x^k-x^{k-1}\|^2, \forall k=1,\ldots.
	\end{equation}
	Let $\beta_k = \frac{\theta_k-1}{\theta_{k+1}}$. As $v^k$ can take either $x^k$ or $x^k + \beta_k (x^k-x^{k-1})$, then we have 
	$$\|v^k - x^{k+1}\|^2 = \begin{cases}
		\|x^k - x^{k+1}\|^2, & \text{if } v^k = x^k;\\
		\|x^k - x^{k+1} + \beta_k (x^k-x^{k-1})\|^2, & \text{if } v^k = x^k + \beta_k (x^k-x^{k-1}).
	\end{cases}$$
	We get from the inequalities $$\|x^k - x^{k+1} + \beta_k (x^k-x^{k-1})\|^2 \geq \frac{1}{2}\|x^k - x^{k+1}\|^2 - \|\beta_k (x^k-x^{k-1})\|^2$$
	and
	$$\|x^k - x^{k+1}\|^2 \geq \frac{1}{2}\|x^k - x^{k+1}\|^2 - \|\beta_k (x^k-x^{k-1})\|^2$$
	that 
	\begin{equation}
		\label{eq:ineqvk-xkp1}
		\|v^k - x^{k+1}\|^2 \geq \frac{1}{2}\|x^k - x^{k+1}\|^2 - \beta_k^2\|x^k-x^{k-1}\|^2.
	\end{equation}
	It follows from \eqref{eq:ineq03-bis}, \eqref{eq:ineqvk-xkp1} and $\rho_g + \rho_h - \gamma_k > 0$ (see \Cref{lem:propEk&ck}-(i)) that 
	$$\begin{aligned}
		f(x^{k+1}) + \frac{\rho_g + \rho_h-\gamma_k}{4}\| x^{k+1} - x^k\|^2 \leq & f(v^k) + \left(\frac{\gamma_k}{2} - \frac{\beta_k^2(\rho_g + \rho_h-\gamma_k)}{4}\right) \| x^k-x^{k-1}\|^2\\
		=& f(v^{k}) + \frac{\rho_g + \rho_h-\gamma_k}{4}\| x^k-x^{k-1}\|^2 \\
		&- c_k\| x^k-x^{k-1}\|^2,
	\end{aligned}$$
	where 
	\begin{equation}
		\label{eq:ineqcoef}
		c_k: = \frac{(\beta_k^2-3)\gamma_k + (1-\beta_k^2) (\rho_g + \rho_h)}{4} \geq \delta > 0
	\end{equation}
	due to \Cref{lem:propEk&ck}-(ii).\\
	Observing that $E_k \geq f(v^k) + \frac{\rho_g + \rho_h-\gamma_k}{4}\| x^k-x^{k-1}\|^2$ for all $k=0,1,2,\ldots$ (see \Cref{lem:propEk&ck}-(v)) and $c_k\geq \delta$, we get for all $k=0,1,2,\ldots$,
	\begin{equation}
		\label{eq:ineq04-bis}
		f(x^{k+1}) + \frac{\rho_g + \rho_h-\gamma_k}{4}\| x^{k+1} - x^k\|^2 
		\leq E_k - \delta\| x^k-x^{k-1}\|^2.
	\end{equation}
	Now, we can prove by induction that for all $t=0, \ldots,q$, the next inequality holds 
	\begin{equation}
		\label{eq:ineq*}
		f(x^{k+1+t}) + \frac{\rho_g + \rho_h-\gamma_{k+t}}{4}\| x^{k+1+t}-x^{k+t}\|^2 
		\leq E_k - \delta\| x^{k+t}-x^{k+t-1}\|^2.
	\end{equation}
	First, it follows from \eqref{eq:ineq04-bis} that the claim holds for $t=0$. Suppose that it holds for $t=0,\ldots,p-1$ with $1\leq p\leq q$. Then, we get from $\delta>0$ that  
	\begin{equation}
		\label{eq:ineq05-bis}
		f(x^{k+1+t}) + \frac{\rho_g + \rho_h-\gamma_{k+t}}{4}\| x^{k+1+t}-x^{k+t}\|^2 
		\leq E_k, \forall t=0,\ldots,p-1.
	\end{equation}
	Replacing $k$ by $k+p$ in \eqref{eq:ineq04-bis}, then 
	$$\begin{aligned}
		f(x^{k+1+p}) + \frac{\rho_g + \rho_h-\gamma_{k+p}}{4}\| x^{k+1+p}-x^{k+p}\|^2 
		\leq& E_{k+p} - \delta\| x^{k+p}-x^{k+p-1}\|^2\\
		\leq&E_{k} - \delta\| x^{k+p}-x^{k+p-1}\|^2,
	\end{aligned}$$
	where the second inequality is due to $E_{k+p}\leq E_k$ with $1\leq p\leq q$ and with \eqref{eq:ineq05-bis} (see \Cref{lem:propEk&ck}-(vi)). Hence, we proved by induction that \eqref{eq:ineq*} holds for all $t=0,\ldots,q$. 
	Therefore,  
	\begin{eqnarray*}
		E_{k+1+q}& = & \max\{ f(x^{k+1+t}) + \frac{\rho_g + \rho_h-\gamma_{k+t}}{4}\| x^{k+1+t}-x^{k+t}\|^2 : t=0,\ldots,q \}\\ 
		&\overset{\eqref{eq:ineq*}}{\leq}& E_k - \delta\min\left\{\| x^{k+t}-x^{k+t-1}\|^2: t=0,\ldots,q\right\}\\
		&=&E_{k} - \delta\| x^{\phi(k)}-x^{\phi(k)-1}\|^2,
	\end{eqnarray*}
	That is the required inequality 
	\begin{equation}
		\label{eq:descentlemma-bis}
		\boxed{E_{k+1+q} \leq E_k - \delta \|x^{\phi(k)}-x^{\phi(k)-1}\|^2, \forall k=0,1,2,\ldots.}
	\end{equation}
	\textbf{(convergence of $\{\|x^{\phi(k)}-x^{\phi(k)-1}\|\}$):} 
	Summing \eqref{eq:descentlemma-bis} for $k$ from $0$ to $N$ (with $N\geq q+1$), we get 
	\begin{eqnarray*}
	\sum_{k=0}^{N}\|x^{\phi(k)}-x^{\phi(k)-1}\|^2 &\leq& \frac{1}{\delta} \sum_{t=0}^{q}(E_t - E_{N+t+1})\\
	&\leq& \frac{q+1}{\delta } (E_q - \inf\{ f(x): x\in \setC\}),
	\end{eqnarray*}
	where the second inequality is derived from $E_t\leq E_q, \forall t=0,\ldots, q$ (see \Cref{lem:propEk&ck}-(iv)) and $E_k\geq \inf\{ f(x): x\in \setC\} > -\infty, \forall k=1,2, \ldots$ (see \Cref{lem:propEk&ck}-(iii)). \\
	Passing $N$ to $\infty$, then 
	$$\sum_{k=1}^{\infty}\|x^{\phi(k)}-x^{\phi(k)-1}\|^2\leq \frac{q+1}{\delta} (E_q - \inf\{ f(x): x\in \setC\}) < \infty.$$
	Therefore,  
	\begin{equation}
		\label{eq:convsubseq}
		\boxed{\|x^{\phi(k)}-x^{\phi(k)-1}\| \xrightarrow{k\to \infty}0.}
	\end{equation} 
	\textbf{(convergence of $\{\|v^{\phi(k)}-x^{\phi(k)-1}\|\}$):}
	By the definition of $v^k$, we have 
	\begin{equation}
		\label{eq:vk-xk-1}
		\|v^k - x^{k-1}\|^2 = \begin{cases}
			\|x^k - x^{k-1}\|^2, & \text{if } v^k = x^k;\\
			\left(1+\beta_k\right)^2\|x^k - x^{k-1}\|^2, & \text{if } v^k = x^k + \beta_k (x^k-x^{k-1}).
		\end{cases}
	\end{equation}
	Since 
	\begin{equation}
		\label{eq:ineq06}
		\left(1+\beta_k\right)^2 \leq (1+\bar{\beta})^2 < 4, \forall k=0,1,2,\ldots.
	\end{equation} 
	It follows from \eqref{eq:convsubseq}, \eqref{eq:vk-xk-1} and \eqref{eq:ineq06} that 
	$$\boxed{\|v^{\phi(k)} - x^{\phi(k)-1}\|\xrightarrow{k\to \infty}0.}$$
	\textbf{(subsequential convergence of $\{x^k\}$):} Let $q=0$. Then we have 
	$$\phi(k) = k\quad \text{ and } \quad E_k = f(x^{k}) + \frac{\rho_g + \rho_h-\gamma_k}{4}\| x^{k} - x^{k-1}\|^2.$$
	Hence, the previously established sufficiently descent property turns to 
	$$E_{k+1} \leq E_k - \delta \|x^{k}-x^{k-1}\|^2, \quad k=0,1,2,\ldots,
	$$
	and we have 
	$$
			\|x^{k}-x^{k-1}\|\xrightarrow{k\to \infty} 0\quad \text{ and }\quad \|v^{k}-x^{k-1}\|\xrightarrow{k\to \infty} 0.
	$$
	The boundedness of the sequence $\{x^k\}_{k\geq 1}\subset \setC$ implies that its set of cluster points is nonempty. The closedness of $\setC$ indicates that all cluster points of $\{x^k\}$ belong to $\setC$. Then for any cluster point $x^*$ of the sequence $\{x^k\}$, there exists a convergent subsequence denoted by $\{x^{k_j}\}_{j\in \N}$ converging to $x^*$. We get from $x^{k_j}\xrightarrow{j\to \infty} x^*$, $\|x^{k}-x^{k-1}\| \xrightarrow{k\to \infty}0$ and $\|v^{k} - x^{k-1}\|\xrightarrow{k\to \infty}0$ that $$x^{k_j-1} \xrightarrow{j\to \infty} x^*, \quad x^{k_j+1} \xrightarrow{j\to \infty} x^* \text{ and } v^{k_j} \xrightarrow{j\to \infty} x^*.$$ It follows from the first order optimality condition of the convex subproblem $$\nabla h(v^k)+\gamma_k(x^k-x^{k-1}) \in \nabla g(x^{k+1}) + N_{\setC}(x^{k+1}),$$
	the closedness of the graph of $N_{\setC}$, the continuity of $\nabla g$ and $\nabla h$, and the boundedness of $\gamma_k$  ($0\leq \gamma_k < \rho_g+\rho_h$) that 
	$$\lim_{j\to\infty} \nabla h(v^{k_j}) + \gamma_{k_j} (x^{k_j} - x^{k_j-1}) \in \nabla g(\lim_{j\to \infty}x^{k_j+1}) + N_{\setC}(\lim_{j\to \infty}x^{k_j+1}).$$
	That is,
	$$\boxed{\nabla h(x^*) \in \nabla g(x^*) + N_{\setC}(x^*).}$$
	Hence, any cluster point $x^*$ of $\{x^k\}$ is a critical point of \eqref{prob:P}.
\end{proof}
\begin{remark}
	The convergence of $\{x^k\}$ and the rate of convergence for both $\{f(x^k)\}$ and $\{x^k\}$ under the Kurdyka-{\L}ojasiewicz property can be established in a similar way as in \cite{niu2022convergence}. Therefore, we admit their results and omit their discussions.
\end{remark}
\begin{remark}
	We only establish the subsequential convergence of $\{x^k\}$ for the case where $q=0$, however, we observe in practice that the sequence $\{x^k\}$ seems converges as well for $q>0$, and often benefits from a better acceleration and superior computed solutions compared to the case when $q=0$, despite the convergence analysis for $q>0$ remains an open challenge. Therefore, a pragmatic suggestion is to initially set $q>0$ to take advantage of better acceleration, then switch to $q=0$ after some iterations to ensure the convergence. For example, we can switch to $q=0$ when  $\beta_k>\bar{\beta}$ (with $\bar{\beta} = 0.9$) is met for the first time.
\end{remark}
The next lemma is required in the proof of \Cref{thm:convHDCA-NI}:
\begin{lemma}\label{lem:propEk&ck}
	Under the assumptions of \Cref{thm:convHDCA-NI}, we have 
	\begin{itemize}[leftmargin=12pt]
		\item[(i)] $\gamma_{k} < \rho_g+\rho_h$ for all $k=0,1,2,\ldots$.
		\item[(ii)] $c_k \geq \delta > 0$ for all $k=0,1,2,\ldots$.
		\item[(iii)] $E_k\geq \inf \{f(x): x\in \setC\}>-\infty,$ for all $k=1,2,\ldots.$
		\item[(iv)] $E_k\leq E_q,$ for all $k=0,\ldots, q$.
		\item[(v)] $E_k \geq f(v^k) + \frac{\rho_g + \rho_h-\gamma_k}{4}\| x^k-x^{k-1}\|^2,$ for all $k=0,1,2,\ldots$.
		\item[(vi)] Suppose that $E_k\geq f(x^{t}) + \frac{\rho_g + \rho_h-\gamma_t}{4}\| x^{t} - x^{t-1}\|^2$ for all $t=k+1,\ldots, k+p$ with $1\leq p\leq q$. Then $E_{k+p}\leq E_k$.
	\end{itemize}
\end{lemma}
\begin{proof}
	(i) For all $k=0,1,2,\ldots$, we get from the definition of $\delta=\frac{(1-\bar{\beta}^2)(\rho_g + \rho_h)}{4}$ and the inequalities $0\leq \beta_k\leq \bar{\beta}<1$ that $$\gamma_{k} \leq \frac{(\rho_g+\rho_h)(1-\beta_k^2) - 4 \delta}{3-\beta_k^2} = \frac{(\rho_g+\rho_h)(\bar{\beta}^2-\beta_k^2)}{3-\beta_k^2} <\frac{(\rho_g+\rho_h)(1-\beta_k^2)}{3-\beta_k^2} < \rho_g + \rho_h.$$
	(ii) For all $k=0,1,2,\ldots$, we get from the definition of $c_k$ that 
	$$c_k = \frac{(\beta_k^2-3)\gamma_k + (1-\beta_k^2) (\rho_g + \rho_h)}{4}.$$
	Then, taking any $\gamma_{k} \in [0, \frac{(\rho_g+\rho_h)(1-\beta_k^2) - 4 \delta}{3-\beta_k^2}]$ and bearing in mind that $0\leq \beta_k\leq\bar{\beta}<1$, we have 
	$$\frac{(\beta_k^2-3)\gamma_k + (1-\beta_k^2) (\rho_g + \rho_h)}{4} \geq  \frac{(\beta_k^2-3)\frac{(\rho_g+\rho_h)(1-\beta_k^2) - 4 \delta}{3-\beta_k^2} + (1-\beta_k^2) (\rho_g + \rho_h)}{4}=\delta>0.$$
	Hence,
	$$c_k\geq \delta > 0, \quad k=0,1,2,\ldots.$$ 
	(iii) For all $k=1,2,\ldots,$ we have 
	$$\begin{aligned}
		E_k =& \max \{f(x^{t}) + \frac{\rho_g + \rho_h-\gamma_t}{4}\| x^{t} - x^{t-1}\|^2: t=\max\{0,k-q\},\ldots, k\}\\
		\geq& \max \{f(x^{t}): t=\max\{0,k-q\},\ldots, k\}\\
		\geq& \inf \{f(x): x\in \setC\} > -\infty.
	\end{aligned}$$
	Note that $k=0$ is not considered because $x^0$ may not be a point in $\setC$. Hence the second-to-last inequality may not hold when $k=0$.\\ 
	(iv) For all $k=0, \ldots, q$, we have 
	$$\begin{aligned}
		E_q =& \max \{f(x^{t}) + \frac{\rho_g + \rho_h-\gamma_t}{4}\| x^{t} - x^{t-1}\|^2: t=0,\ldots, q\} \\
		\geq& \max \{f(x^{t}) + \frac{\rho_g + \rho_h-\gamma_t}{4}\| x^{t} - x^{t-1}\|^2: t=0,\ldots, k\}= E_k.
	\end{aligned}$$
	(v) There are two possible cases for $v^k$: \\
	$\bullet$ $v^k=x^k$, then $$f(v^k) + \frac{\rho_g + \rho_h-\gamma_k}{4}\| x^k-x^{k-1}\|^2 = f(x^k) + \frac{\rho_g + \rho_h-\gamma_k}{4}\| x^k-x^{k-1}\|^2 \leq E_k.$$
	$\bullet$ $v^k=x^k + \frac{\theta_k -1}{\theta_{k+1}}\left(x^k - x^{k-1} \right)$, this will only occur when $f(v^k) + \frac{\rho_g + \rho_h-\gamma_k}{4}\| x^k-x^{k-1}\|^2 \leq E_k$. \\
	(vi) For any $p=1,\ldots, q$, we get from the definition of $E_{k+p}$ that 
	$$\begin{aligned}
		E_{k+p} =& \max \{f(x^{t}) + \frac{\rho_g + \rho_h-\gamma_t}{4}\| x^{t} - x^{t-1}\|^2: t=\max\{0,k+p-q\},\ldots, k+p\} \\
		\leq& \max \{f(x^{t}) + \frac{\rho_g + \rho_h-\gamma_t}{4}\| x^{t} - x^{t-1}\|^2: t=\max\{0,k-q\},\ldots, k+p\} \\
		=& \max\{E_k,  \max \{f(x^{t}) + \frac{\rho_g + \rho_h-\gamma_t}{4}\| x^{t} - x^{t-1}\|^2: t=k+1,\ldots, k+p\}\}= E_k,
	\end{aligned}$$
where the last equality is due to the assumption that $E_k\geq f(x^{t}) + \frac{\rho_g + \rho_h-\gamma_t}{4}\| x^{t} - x^{t-1}\|^2$ for all $t=k+1,\ldots, k+p$.
\end{proof}

\section{DC formulations for \eqref{eq:aeicp}}\label{sec:DCformulations}
In this section, we will present four equivalent DC formulations for \eqref{eq:aeicp} using a novel DC decomposition technique for polynomials, the difference-of-convex-sums-of-squares (DC-SOS) decomposition, introduced in \cite{niu2018difference}. This entails representing any polynomial as difference of convex sums-of-square polynomials. 

\subsection{First DC formulation}

Consider the nonlinear programming (NLP) formulation of \eqref{eq:aeicp} presented in \cite{Judice09,Niu12} as
\begin{equation}\label{prob:nlp1}
	0 = \min\{f_1(x,y,w,z) := \|y-z x\|^2 + x^{\top}w :  (x,y,w,z)\in \setC_1\},\tag{NLP1}
\end{equation}
where $$\setC_1: = \{(x,y,w,z)\in \R^n_+ \times \R^n_+ \times \R^n_+ \times \R_+ : w = Bx -A y, e^{\top}x = 1, e^{\top}y=z\}$$
and $e$ denotes the vector of ones. The gradient of $f_1$ is computed by 
\begin{equation}
	\label{eq:df1}
	\begin{cases}
		\nabla_x f_1(x,y,w,z) = 2 z(zx-y) + w,\\
		\nabla_y f_1(x,y,w,z) = 2(y-zx),\\
		\nabla_w f_1(x,y,w,z) = x,\\
		\nabla_z f_1(x,y,w,z) = 2x^{\top}(zx-y).
	\end{cases}
\end{equation}
It is known in \cite[Theorem 3.1]{Judice09} that for any global optimal solution $(\bar{x},\bar{y},\bar{w},\bar{z})$ of \eqref{prob:nlp1} with zero optimal value, we have $\bar{y}=\bar{z}\bar{x}$ and $$(\bar{x}, \frac{1}{\bar{z}})\in \AEiCP{A,B}.$$
Unlike the SEiCP, a stationary point of \eqref{prob:nlp1} is not necessarily to be a solution of $\AEiCP{A,B}$. A further discussion in \cite[Theorem 3.2]{Judice09} shows that a stationary point of \eqref{prob:nlp1} is a solution of $\AEiCP{A,B}$ if and only if the Lagrange multipliers associated with the linear equalities $e^{\top}x=1$ and $e^{\top}y = z$ equal $0$.

A DC-SOS decomposition for $\|y-z x\|^2$ is given by $g(x,y,z)-h(x,y,z)$ where
\begin{equation}
	\label{eq:ded_y-zx}
	\begin{cases}
		g(x,y,z) = &\|y\|^2 + \frac{((z+1)^2 + \|y-x\|^2)^2 + ((z-1)^2 + \|y+x\|^2)^2}{16} + \frac{(z^2 	+ \|x\|^2)^2}{2},\\
		h(x,y,z) = &\frac{((z+1)^2 + \|y+x\|^2)^2 + ((z-1)^2 + \|y-x\|^2)^2}{16} + \frac{z^4 + \|x\|^4}{2},
	\end{cases}
\end{equation}
and a DC-SOS decomposition for $x^{\top}w$ reads
$$x^{\top}w = \frac{\|x+w\|^2}{4} - \frac{\|x-w\|^2}{4}.$$
Hence, $f_1$ has a DC-SOS decomposition $G_1-H_1$ where  
\begin{equation}
	\label{eq:G1H1}
	\begin{cases}
		G_1(x,y,w,z) = &\|y\|^2 + \frac{((z+1)^2 + \|y-x\|^2)^2 + ((z-1)^2 + \|y+x\|^2)^2}{16} + \frac{(z^2 	+ \|x\|^2)^2}{2} + \frac{\|x+w\|^2}{4},\\
		H_1(x,y,w,z) = &\frac{((z+1)^2 + \|y+x\|^2)^2 + ((z-1)^2 + \|y-x\|^2)^2}{16} + \frac{z^4 + \|x\|^4}{2} + \frac{\|x-w\|^2}{4}.
	\end{cases}
\end{equation}
Thus, problem \eqref{prob:nlp1} has a DC formulation as
\begin{equation}
	\label{prob:dcp1}
	0 = \min\{G_1(x,y,w,z) - H_1(x,y,w,z) : (x,y,w,z)\in \setC_1\}. \tag{DCP1}
\end{equation}
The gradient of $H_1$ is computed by 
\begin{equation}\label{eq:dH1}
	\begin{cases}
		\nabla_x H_1(x,y,w,z) &=\frac{((z+1)^2 + \|y+x\|^2)(x+y)}{4} + \frac{((z-1)^2 + \|y-x\|^2)(x-y)}{4} + \frac{x-w}{2} + 2\|x\|^2 x,\\
		\nabla_y H_1(x,y,w,z) &=\frac{((z+1)^2 + \|y+x\|^2)(x+y)}{4} - \frac{((z-1)^2 + \|y-x\|^2)(x-y)}{4},\\
		\nabla_w H_1(x,y,w,z) &=\frac{w-x}{2},\\
		\nabla_z H_1(x,y,w,z) &=\frac{((z+1)^2 + \|y+x\|^2)(z+1)}{4} + \frac{((z-1)^2 + \|y-x\|^2)(z-1)}{4} + 2z^3.
	\end{cases}
\end{equation}

\subsection{Second DC formulation} 

Replacing $w$ by $Bx-Ay$ in \eqref{prob:nlp1}, we get another NLP formulation
\begin{equation}\label{prob:nlp2}
	0 = \min\{f_2(x,y,z) := \|y-z x\|^2 + x^{\top}(Bx-Ay) :  (x,y,z)\in \setC_2\},\tag{NLP2}
\end{equation}
where $$\setC_2: = \{(x,y,z)\in \R^n_+ \times \R^n_+ \times \R_+ : Bx -A y\geq 0, e^{\top}x = 1, e^{\top}y = z\}.$$
The gradient of $f_2$ is computed by
\begin{equation}
	\label{eq:df2}
	\begin{cases}
		\nabla_x f_2(x,y,z) = 2 z(zx-y) + (B+B^{\top}) x - Ay,\\
		\nabla_y f_2(x,y,z) = 2(y-zx) - A^{\top}x,\\
		\nabla_z f_2(x,y,z) = 2x^{\top}(zx-y).
	\end{cases}
\end{equation}
In virtue of the equivalence between \eqref{prob:nlp1} and \eqref{prob:nlp2}, it follows immediately from \cite[Theorem 3.1]{Judice09} that for any global optimal solution $(\bar{x},\bar{y},\bar{z})$ of \eqref{prob:nlp2} with zero optimal value, we have $\bar{y}=\bar{z}\bar{x}$ and $$(\bar{x}, \frac{1}{\bar{z}})\in \AEiCP{A,B}.$$

A DC-SOS decomposition for $x^{\top}Ay$ reads
\begin{equation}
	\label{eq:ded_xAy}
	x^{\top}Ay = \frac{\|x + Ay\|^2}{4} - \frac{\|x - Ay\|^2}{4},
\end{equation}
a DC-SOS decomposition for $\|y-zx\|^2$ is given in \eqref{eq:ded_y-zx}, and $x^{\top}B x$ is convex since $B\in \PD$. Hence, $f_2$ has a DC-SOS decomposition $G_2-H_2$ where  
{\footnotesize 
	\begin{equation}
		\label{eq:GHfornlp2}
		\begin{cases}
			G_2(x,y,z) = &\|y\|^2 + \frac{((z+1)^2 + \|y-x\|^2)^2 + ((z-1)^2 + \|y+x\|^2)^2}{16} + \frac{(z^2 	+ \|x\|^2)^2}{2} + x^{\top}Bx + \frac{\|x - Ay\|^2}{4}, \\
			H_2(x,y,z) = &\frac{((z+1)^2 + \|y+x\|^2)^2 + ((z-1)^2 + \|y-x\|^2)^2}{16} + \frac{z^4 + \|x\|^4}{2} + \frac{\|x + Ay\|^2}{4}.
		\end{cases}
\end{equation}}
Then problem \eqref{prob:nlp2} has a DC formulation as
\begin{equation}
	\label{prob:dcp2}
	0 = \min\{G_2(x,y,z) - H_2(x,y,z) : (x,y,z)\in \setC_2\}. \tag{DCP2}
\end{equation}
The gradient of $H_2$ is computed by 
\begin{equation}
	\begin{cases}
		\nabla_x H_2(x,y,z) &=\frac{((z+1)^2 + \|y+x\|^2)(x+y)}{4} + \frac{((z-1)^2 + \|y-x\|^2)(x-y)}{4} + \frac{x+Ay}{2} + 2\|x\|^2 x,\\
		\nabla_y H_2(x,y,z) &=\frac{((z+1)^2 + \|y+x\|^2)(x+y)}{4} - \frac{((z-1)^2 + \|y-x\|^2)(x-y)}{4} + \frac{A^{\top}(Ay+x)}{2},\\
		\nabla_z H_2(x,y,z) &=\frac{((z+1)^2 + \|y+x\|^2)(z+1)}{4} + \frac{((z-1)^2 + \|y-x\|^2)(z-1)}{4} + 2z^3.
	\end{cases}
\end{equation}

\subsection{Third DC formulation} 
An NLP formulation proposed in \cite{Niu12} reads
\begin{equation}\label{prob:nlp3}
	0 = \min\{f_3(x,y,w) := \|y\|^2 + x^{\top}w - \frac{(x^{\top}y)^2}{\|x\|^2} :  (x,y,w)\in \setC_3\},\tag{NLP3}
\end{equation}
where $$\setC_3: = \{(x,y,w)\in \R^n_+ \times \R^n_+ \times \R^n_+ : w = Bx -A y, e^{\top}x = 1\}.$$
The gradient of $f_3$ is computed by
\begin{equation}
	\label{eq:df3}
	\begin{cases}
		\nabla_x f_3(x,y,w) = w - 2\frac{x^{\top}y}{\|x\|^2} y + 2\frac{(x^{\top}y)^2}{\|x\|^4} x,\\
		\nabla_y f_3(x,y,w) = 2y -  2\frac{x^{\top}y}{\|x\|^2} x,\\
		\nabla_w f_3(x,y,w) = x.
	\end{cases}
\end{equation}
\begin{theorem}\label{thm:NLP3}
	Under \Cref{hyp1}, let $(\bar{x},\bar{y},\bar{w})$ be a global optimal solution of \eqref{prob:nlp3} with zero optimal value, then $$(\bar{x},\frac{\|\bar{x}\|}{\|\bar{y}\|})\in \AEiCP{A,B}.$$  
\end{theorem}
\begin{proof}
	Let $(\bar{x},\bar{y},\bar{w})$ be a global optimal solution of \eqref{prob:nlp3} with zero optimal value, we have 
	\begin{equation}
		\label{eq:thm:NLP3_eq1}
		\begin{cases}
			\|\bar{y}\|^2 + \bar{x}^{\top}\bar{w} - \frac{(\bar{x}^{\top}\bar{y})^2}{\|\bar{x}\|^2} = 0,\\
			\bar{w} = B \bar{x} - A\bar{y}\geq 0, e^{\top} \bar{x} = 1, \bar{x}\geq 0, \bar{y}\geq 0.
		\end{cases}
	\end{equation}
	Then 
	$$0 = \|\bar{y}\|^2 + \bar{x}^{\top}\bar{w} - \frac{(\bar{x}^{\top}\bar{y})^2}{\|\bar{x}\|^2} \geq \|\bar{y}\|^2 + \bar{x}^{\top}\bar{w} - \frac{\|\bar{x}\|^2\|\bar{y}\|^2}{\|\bar{x}\|^2} =  \bar{x}^{\top}\bar{w} \geq 0,$$
	where the first inequality is due to Cauchy-Schwartz $|\bar{x}^{\top}\bar{y}|\leq \|\bar{x}\|\|\bar{y}\|$, and the last inequality is due to $\bar{x}\geq 0$ and $\bar{w}\geq 0$. Hence, 
	\begin{equation}
		\label{eq:thm:NLP3_eq2}
		\bar{x}^{\top}\bar{w} = 0 \quad \text{ and } \quad |\bar{x}^{\top}\bar{y}| = \|\bar{x}\|\|\bar{y}\|.
	\end{equation}
	It follows from $|\bar{x}^{\top}\bar{y}| = \|\bar{x}\|\|\bar{y}\|, \bar{x}\geq 0$, $\bar{y}\geq 0$ and $e^{\top}\bar{x}=1$ that there exists a positive scalar $\bar{z}>0$ such that 
	\begin{equation}
		\label{eq:thm:NLP3_eq3}		
		\bar{y}=\bar{z}\bar{x}, \quad \bar{x}\neq 0. 
	\end{equation}
	Therefore,
	\begin{equation}
		\label{eq:thm:NLP3_eq4}
		\bar{w} = B \bar{x} - A\bar{y} = B \bar{x} - \bar{z}A\bar{x}\geq 0.
	\end{equation}
	Combining \eqref{eq:thm:NLP3_eq1}, \eqref{eq:thm:NLP3_eq2} and \eqref{eq:thm:NLP3_eq4}, we proved that  
	$$(\bar{x},\frac{1}{\bar{z}})\in \AEiCP{A,B}.$$  
	$\bar{z}$ is computed by injecting \eqref{eq:thm:NLP3_eq3} to $|\bar{x}^{\top}\bar{y}| = \|\bar{x}\|\|\bar{y}\|$ as
	$\bar{z} = \|\bar{y}\|/\|\bar{x}\|.$
\end{proof}
Consider the objective function of \eqref{prob:nlp3}: $$f_3(x,y,w) = \|y\|^2 + x^{\top}w - \frac{(x^{\top}y)^2}{\|x\|^2}.$$
Let $\varphi(x,y):= (x^{\top}y)^2/\|x\|^2$. Then $\varphi$ is a smooth non-convex function over $\Omega\times \R_+^n$, with its gradient and Hessian computed by 
$$\nabla_x \varphi(x,y) = 2\frac{x^{\top}y}{\|x\|^2} y - 2\frac{(x^{\top}y)^2}{\|x\|^4} x, \quad \nabla_y \varphi(x,y) = 2\frac{x^{\top}y}{\|x\|^2} x,$$
and 
$$\nabla^2 \varphi(x,y) = \begin{bmatrix}
	\nabla_{x,x}^2 \varphi(x,y) & \nabla_{x,y}^2 \varphi(x,y)\\
	\nabla_{y,x}^2 \varphi(x,y) & \nabla_{y,y}^2 \varphi(x,y)
\end{bmatrix},$$
where $$\begin{cases}
	\nabla_{x,x}^2 \varphi (x,y) = \frac{2}{\|x\|^2}(yy^{\top}) + \frac{8(x^{\top}y)^2}{\|x\|^6}(xx^{\top}) - \frac{2(x^{\top}y)^2}{\|x\|^4}I_n - \frac{4x^{\top}y}{\|x\|^4}(xy^{\top}+yx^{\top}),\\
	\nabla_{x,y}^2 \varphi (x,y) = (\nabla_{y,x}^2 \varphi (x,y))^{\top} = \frac{2}{\|x\|^2}(xy^{\top}) + \frac{2x^{\top}y}{\|x\|^2}I_n -\frac{4(x^{\top}y)}{\|x\|^4}(xx^{\top}),\\
	\nabla_{y,y}^2 \varphi(x,y) = \frac{2}{\|x\|^2} (xx^{\top}).
\end{cases}$$
The next proposition gives a DC decomposition for $\varphi(x,y)$.
\begin{proposition}\label{prop:dcdofphi} The function $\varphi(x,y)= (x^{\top}y)^2/\|x\|^2$ has a DC decomposition over $\setC_3$ as 
	$$\varphi(x,y) = \left[\frac{\eta}{2}\|(x,y)\|^2 + \varphi(x,y)\right] - \frac{\eta}{2}\|(x,y)\|^2,$$
	for any 
	\begin{equation}
		\label{eq:estimateeta}
		\eta \geq 3.2 + 20 n M^2,
	\end{equation} 
	where $M$ is computed by solving the linear program
	\begin{equation}
		\label{prob:lpforM}
			M=\max\{ e^{\top}y: (x,y,w)
		\in \setC_3\}.
	\end{equation}
\end{proposition}
\begin{proof}
	We can compute a large enough $\eta$ by estimating an upper bound of the spectral radius of $\nabla^2 \varphi(x,y)$ over $\Omega\times \R_+^n$. Since $(y^{\top}x)$ is the only possible nonzero eigenvalue of the rank-one matrix $xy^{\top}$ with the associated eigenvector $x$, and $\|xy^{\top}\| = \|x\|\|y\|$. Then
	$$\|\nabla_{x,x}^2\varphi(x,y)\| \leq 20\frac{\|y\|^2}{\|x\|^2}, \quad \|\nabla_{x,y}^2\varphi(x,y)\| \leq 8\frac{\|y\|}{\|x\|}, \quad \|\nabla_{y,y}^2\varphi(x,y)\| =2.$$
	Hence $$\begin{aligned}
		\rho(\nabla^2 \varphi(x,y)) &= \|\nabla^2 \varphi(x,y)\|\\
		&\leq \sqrt{\|\nabla_{x,x}^2 \varphi(x,y)\|^2 + 2 \|\nabla_{x,y}^2 \varphi(x,y)\|^2 + \|\nabla_{y,y}^2 \varphi(x,y)\|^2}\\
		& \leq 20\frac{\|y\|^2}{\|x\|^2} + 3.2.
	\end{aligned}$$ 
	The term $\|y\|/\|x\|$ over $\setC_3$ is upper bounded by 
	$$\frac{\|y\|}{\|x\|} \leq \frac{e^{\top}y}{\|x\|} \leq \sqrt{n} e^{\top}y \leq \sqrt{n} \max_{(x,y,z)
		\in \setC_3} e^{\top}y.$$
	Then, taking any 
	$$\eta \geq 3.2 + 20 n \left(\max_{(x,y,z)
		\in \setC_3} e^{\top}y\right)^2,$$
	we get the desired DC decomposition for $\varphi(x,y)$ over $\setC_3$.
\end{proof}
We derive from \Cref{prop:dcdofphi} a DC-SOS decomposition for $f_3$ as $G_3-H_3$: \begin{equation}
	\label{eq:G3H3}
	\begin{cases}
		G_3(x,y,w) = (1+\frac{\eta}{2})\|y\|^2 +  \frac{\eta}{2}\|x\|^2 + \frac{\|x+w\|^2}{4},\\
		H_3(x,y,w) = \frac{\eta}{2}\|(x,y)\|^2 + \frac{(x^{\top}y)^2}{\|x\|^2} + \frac{\|x-w\|^2}{4}.
	\end{cases}
\end{equation}
Then problem \eqref{prob:nlp3} has a DC formulation as
\begin{equation}
	\label{prob:dcp3}
	0 = \min\{G_3(x,y,w) - H_3(x,y,w) : (x,y,w)\in \setC_3\}.\tag{DCP3}
\end{equation}
The gradient of $H_3$ is computed by 
\begin{equation}
	\begin{cases}
		\nabla_x H_3(x,y,w) &=\eta x + 2\frac{x^{\top}y}{\|x\|^2} y - 2\frac{(x^{\top}y)^2}{\|x\|^4} x  + \frac{(x-w)}{2},\\
		\nabla_y H_3(x,y,w) &=\eta y + 2\frac{x^{\top}y}{\|x\|^2} x,\\
		\nabla_w H_3(x,y,w) &=\frac{w-x}{2}.
	\end{cases}
\end{equation}

\section{BDCA and DCA for solving \eqref{eq:aeicp}} \label{sec:BDCAforAEiCP}
The DC formulations proposed in the previous section, namely \eqref{prob:dcp1}, \eqref{prob:dcp2} and \eqref{prob:dcp3}, belong to the convex constrained DC program \eqref{prob:P}. In this section, we will discuss how to apply BDCA and DCA for solving these DC formulations. The standard DCA is regarded as a special case of BDCA without line search. The solution methods for convex subproblems required in these DC algorithms and the line search procedure (exact and inexact) required in BDCA will be discussed. 
\subsection{BDCA and DCA for \eqref{prob:dcp1}} BDCA for \eqref{prob:dcp1} is summarized in \Cref{alg:BDCA_for_DCP1}. Some comments are given below: 
\begin{algorithm}[H]
	\caption{BDCA for \eqref{prob:dcp1}}
	\label{alg:BDCA_for_DCP1}
	\begin{algorithmic}[1]
		\REQUIRE $x^0\in \Omega$, $y^0\in \R_+^n$, $\bar{\alpha}>0$;
		\STATE initialize $X^0 \leftarrow (x^0,y^0, Bx^0 - Ay^0,e^{\top}y^0)$; \label{codeline:bdca_dcp1_X0}
		\FOR{$k=0,1,\ldots$}
		\STATE solve \eqref{prob:CP1} to get a solution $V^k$; \label{codeline:bdca_dcp1_cp1}
		\STATE initialize $X^{k+1}\leftarrow V^k$; \label{codeline:bdca_dcp1_initXk1}
		\STATE set $D^k\leftarrow V^k - X^k$;\label{codeline:bdca_dcp1_Dk}
		\IF{$\setA_1(V^k)\subset \setA_1(X^k)$ and $\langle \nabla f_1(V^k),D^k\rangle < 0 $}\label{codeline:bdca_dcp1_cond}
		\STATE set $X^{k+1}\leftarrow $LineSearch$(V^k,D^k,\bar{\alpha})$; \label{codeline:bdca_dcp1_linesearch}
		\ENDIF\label{codeline:bdca_dcp1_endif}
		\ENDFOR
	\end{algorithmic}
\end{algorithm}

\noindent$\bullet$ DCA for \eqref{prob:dcp1} is just BDCA without the codes between line \ref{codeline:bdca_dcp1_Dk} to \ref{codeline:bdca_dcp1_endif}.

\noindent$\bullet$ The convex subproblem in line \ref{codeline:bdca_dcp1_cp1} is defined by \begin{equation}
	\label{prob:CP1}
	V^k\in \underset{(x,y,w,z)\in \setC_1}{\argmin} \{ G_1(x,y,w,z) - \langle (x,y,w,z), \nabla H_1(x^k,y^k,w^k,z^k) \rangle \},\tag{CP1}
\end{equation}
where $G_1$, $H_1$ and $\nabla H_1$ are given in \eqref{eq:G1H1} and \eqref{eq:dH1} respectively. We can of course introduce a strongly quadratic term $\frac{\rho}{2}\|(x,y,w,z)\|^2$ (for any $\rho>0$) into both $G_1$ and $H_1$ to ensure that the DC components are $\rho$ strongly convex. The problem \eqref{prob:CP1} can be obviously solved via many first- and second-order nonlinear optimization approaches such as gradient-type methods, Newton-type methods, interior-point methods. Some NLP solvers are available such as IPOPT, KNITRO, FILTERSD, CVX and MATLAB FMINCON, but these solvers may not be really efficient when dealing with some large-scale and ill-conditioned instances. 

To enhance the efficiency of solving the subproblem, especially for better handling large-scale cases, we propose reformulating \eqref{prob:CP1} as a quadratic programming (QP) problem, which can be addressed using more efficient QP or Second-Order Cone Programming (SOCP) solvers such as MOSEK \cite{mosek}, GUROBI \cite{Gurobi} and CPLEX \cite{Cplex}. The QP formulation is presented below.\\
\textbf{QP formulation for \eqref{prob:CP1}:}
By introducing the additional variable 
\begin{equation}
	\label{eq:u}
	u := (u_x, u_z, u_{z+1}, u_{z-1}, u_{y+x}, u_{y-x})\in \R^6
\end{equation}
and $6$ associated quadratic  constraints
\begin{equation}\label{eq:QC1}
	\setQC_1=\begin{cases}
		\|x\|^2\leq u_x,\\
		\|z\|^2\leq u_z,\\
		(z + 1)^2\leq u_{z+1},\\
		(z - 1)^2\leq u_{z-1},\\
		\|y + x\|^2\leq u_{y+x},\\
		\|y - x\|^2\leq u_{y-x},
	\end{cases}
\end{equation}
then \eqref{prob:CP1} is formulated as the quadratic  program:
\begin{equation}\label{prob:QP1}
	\underset{(x,y,w,z,u)\in \setC_1 \cap \setQC_1}{\min} \bar{G}_1(x,y,w,u) - \langle (x,y,w,z), \nabla H_1(x^k,y^k,w^k,z^k) \rangle,
	\tag{QP1}
\end{equation}
where $$\bar{G}_1(x,y,w,u) : = \|y\|^2 + \frac{(u_{z+1}+u_{y-x})^2 + (u_{z-1}+u_{y+x})^2}{16} + \frac{(u_z+u_x)^2}{2} + \frac{\|x+w\|^2}{4}$$
is derived from $G_1$ by replacing some squares by the corresponding terms of $u$ in \eqref{eq:QC1}. For the strongly convex DC formulation with additional term $\frac{\rho}{2}\|(x,y,w,z)\|^2$ in both $G_1$ and $H_1$, its convex subproblem has a similar QP formulation as 
\begin{equation}\label{prob:QP1-b}
	\underset{(x,y,w,z,u)\in \setC_1 \cap \setQC_1}{\min} \bar{G}_1(x,y,w,u) + \frac{\rho}{2}\|(x,y,w,z)\|^2- \langle (x,y,w,z), \xi^k \rangle,
\end{equation}
where $\xi^k = \rho(x^k,y^k,w^k,z^k) + \nabla H_1(x^k,y^k,w^k,z^k)$.
Note that the constraint $\setC_1\cap \setQC_1$ and the function $\bar{G}_1$ in \eqref{prob:QP1} and \eqref{prob:QP1-b} does not depend on the iteration $k$. Hence, we can generate them before starting the first iteration.

The equivalence between \eqref{prob:CP1} and \eqref{prob:QP1} is established below:
\begin{theorem}\label{thm:QP1}
	Let $(\bar{x},\bar{y},\bar{w},\bar{z},\bar{u})$ be an optimal solution of \eqref{prob:QP1}, then $(\bar{x},\bar{y},\bar{w},\bar{z})$ is an optimal solution of \eqref{prob:CP1}. Conversely, let $(\bar{x},\bar{y},\bar{w},\bar{z})$ be an optimal solution of \eqref{prob:CP1}, then  $(\bar{x},\bar{y},\bar{w},\bar{z},\bar{u})$ with $$\bar{u} =(\|\bar{x}\|^2, \|\bar{z}\|^2, (\bar{z}+1)^2,(\bar{z}-1)^2,\|\bar{y}+\bar{x}\|^2,\|\bar{y}-\bar{x}\|^2)$$
	is an optimal solution of \eqref{prob:QP1}.
\end{theorem}
\begin{proof}
	We just need to show that for any optimal solution of \eqref{prob:QP1}, we have equalities for all quadratic  constraints in $\setQC_1$. By contradiction, let $(\bar{x},\bar{y},\bar{w},\bar{z},\bar{u})$ be an optimal solution of \eqref{prob:QP1} such that there exists a constraint in $\setQC_1$ with strictly inequality, e.g., suppose that $\|\bar{x}\|^2< \bar{u}_x$, then we can take $\bar{u}_x = \|\bar{x}\|^2$ and keep the same values for all other variables, which leads to a better feasible solution with a smaller value of $\bar{G}_1$, hence the optimality assumption is violated. 
\end{proof}

\noindent$\bullet$ $\setA_1(X)$ denotes the active set of $\setC_1$ at $X:=(x,y,w,z)$ defined as $$\setA_1(X) = \{i\in \{1,\ldots,3n+1\}: X_i=0\}.$$ Since $\setC_1$ is polyhedral convex, then $\setA_1(V^k)\subset \setA_1(X^k)$ is a necessary and sufficient condition for $D^k$ to be a descent direction at $V^k$.

\noindent$\bullet$ LineSearch$(V^k,D^k,\bar{\alpha})$ in line \ref{codeline:bdca_dcp1_linesearch} can be either exact or inexact. Let us denote $$V^k=(V_x^k,V_y^k,V_w^k,V_z^k)\quad \text{ and }\quad D^k = (D_x^k,D_y^k,D_w^k,D_z^k).$$  \\
\textbf{Exact line search:} We can simplify  
$$f_1(V^k+\alpha D^k) = \frac{a_1}{4}\alpha^4 + \frac{a_2}{3}\alpha^3 + \frac{a_3}{2}\alpha^2 + a_4\alpha + a_5,$$
where 
\begin{equation}
	\label{eq:coefpolydcp1}
	\left\{\begin{aligned}
		a_1 =& 4(D_z^k)^2\|D_x^k\|^2,\\
		a_2 =& -6\langle D_{z}^{k} D_{x}^{k}, D_{y}^k - D_{z}^k V_{x}^k - V_{z}^k D_{x}^k \rangle, \\
		a_3 =& 2\left( \|D_{y}^k - D_{z}^k V_{x}^k - V_{z}^k D_{x}^k\|^2+\langle D_{w}^{k},D_{x}^{k}\rangle -2\langle D_{z}^{k} D_{x}^{k},V_{y}^k - V_{z}^kV_{x}^k\rangle\right), \\
		a_4 =& 2\langle V_{y}^k - V_{z}^kV_{x}^k, D_{y}^k - D_{z}^k V_{x}^k - V_{z}^k D_{x}^k \rangle +\langle V_{w}^{k}, D_{x}^{k}\rangle + \langle V_{x}^{k},D_{w}^{k}\rangle,\\
		a_5 =& \langle V_{x}^k,V_{w}^k\rangle + \|V_{y}^k - V_{z}^kV_{x}^k\|^2.	
	\end{aligned}\right.
\end{equation}
The exact line search (with upper bounded stepsize $\bar{\alpha}$) at $V^k$ along $D^k$ is to solve
\begin{equation}
	\alpha_k = \argmin \{f_1(V^k + \alpha D^k): 0\leq \alpha \leq \bar{\alpha}\}.
\end{equation}
This problem is equivalent to 
\begin{equation}
	\label{prob:computealphak_dcp1}
	\alpha_k = \argmin \{f_1(V^k + \alpha D^k): \alpha\in \{0,\min\{\bar{\alpha}_k,\bar{\alpha}\}\}\cup \mathcal{Z}\},
\end{equation}
where \begin{equation}
	\label{eq:alphabar}
	\bar{\alpha}_k = \min\left\{ -(V^k)_i/(D^k)_i, i\in \mathcal{I}^k \right\} \text{ with } \mathcal{I}^k=\{i\in \{1,\ldots,3n+1\}: (D^k)_i<0\}
\end{equation} 
under the assumption that $\min \emptyset = \infty$ where $(D^k)_i$ is the $i$-th component of the vector $D^k$, and $\mathcal{Z}$ is the set of all real roots of the cubic polynomial
$$q(\alpha)= \frac{\mathrm{d} f_1(V^k + \alpha D^k)}{\mathrm{d} \alpha}  = a_1~\alpha^3 + a_2~\alpha^2 + a_3~\alpha + a_4.$$

Note that $\mathcal{Z}$ has at most $3$ distinct real roots which can be all computed by the renowned Cardano-Tartaglia formula. Hence problem \eqref{prob:computealphak_dcp1} can be explicitly solved by computing all real roots of $q(\alpha)$ and checking the values of $f_1(V^k+\alpha D^k)$ at $5$ different values of $\alpha$ at most (possibly three real roots of $q(\alpha)$, $0$ and $\min\{\bar{\alpha}_k,\bar{\alpha}\}$). \\
\textbf{Inexact line search:} The Armijo-type inexact line search is described below, where the initial stepsize $\alpha$ is suggested to be $\min\{\bar{\alpha}_k,\bar{\alpha}\}$ with $\bar{\alpha}_k$ defined in \eqref{eq:alphabar}.
\begin{proc}{ht!}{Procedure}
	\caption{Armijo Line Search}
	\label{alg:Armijo}
	\begin{algorithmic}[1]
		\REQUIRE $V$, $D$, $\bar{\alpha}$.
		\ENSURE $Z$.
		\STATE \textbf{Initialization:} reduction factor $\beta\in (0,1)$; initial stepsize $\alpha \in (0,\bar{\alpha}]$; parameter $\sigma\in (0,1)$; tolerance for line search $\varepsilon>0$;
		\WHILE{$\alpha > \varepsilon/\|D\|$}		
		\STATE $Z \leftarrow V + \alpha D$;
		\STATE $\Delta \leftarrow f_1(V) - f_1(Z) - \sigma \alpha^2 \|D\|^2$;
		\IF{$\Delta \geq 0$ and $Z\in \setC_1$}
		\RETURN $Z$;
		\ENDIF
		\STATE $\alpha \leftarrow \beta \alpha$;
		\ENDWHILE
		\STATE $Z \leftarrow V$;			
		\RETURN $Z$.
	\end{algorithmic}
\end{proc}

\subsection{BDCA and DCA for \eqref{prob:dcp2}} BDCA for \eqref{prob:dcp2} is very similar to Algorithm \ref{alg:BDCA_for_DCP1}. We outline the differences as follows: 	

\noindent$\bullet$ The initialization in line \ref{codeline:bdca_dcp1_X0} of Algorithm \ref{alg:BDCA_for_DCP1} is changed to
$$X^0 \leftarrow (x^0,y^0,e^{\top}y^0).$$

\noindent$\bullet$ The convex subproblem in line \ref{codeline:bdca_dcp1_cp1} of Algorithm \ref{alg:BDCA_for_DCP1} is changed to \begin{equation}
	\label{prob:CP2}
	V^k\in \underset{(x,y,z)\in \setC_2}{\argmin} \{ G_2(x,y,z) - \langle (x,y,z), \nabla H_2(x^k,y^k,z^k) \rangle \},\tag{CP2}
\end{equation}
which has a similar structure as \eqref{prob:CP1}. So it can also be solved by QP approach. A QP formulation for \eqref{prob:CP2} is given by \begin{equation}\label{prob:QP2}
	\min\{ \bar{G}_2(x,y,u) - \langle (x,y,z), \nabla H_2(x^k,y^k,z^k) \rangle: (x,y,z,u)\in \setC_2 \cap \setQC_1\},
	\tag{QP2}
\end{equation}
where $u$ is defined in \eqref{eq:u}, $\setQC_1$ is given by \eqref{eq:QC1}, and {\footnotesize$$\bar{G}_2(x,y,u) = \|y\|^2 + \frac{(u_{z+1} + u_{y-x})^2 + (u_{z-1} + u_{y+x})^2}{16} + \frac{(u_z 	+ u_x)^2}{2} + x^{\top}Bx + \frac{\|x - Ay\|^2}{4}.$$} A similar QP formulation for the strongly convex DC decomposition by adding  $\frac{\rho}{2}\|(x,y,z)\|^2$ into $G_2$ and $H_2$ can be established accordingly. The equivalence between \eqref{prob:CP2} and \eqref{prob:QP2} is described in Theorem \ref{thm:QP2}, whose proof shares similarities with that of Theorem \ref{thm:QP1} and is therefore omitted.
\begin{theorem}\label{thm:QP2}
	Let $(\bar{x},\bar{y},\bar{z},\bar{u})$ be an optimal solution of \eqref{prob:QP2}, then $(\bar{x},\bar{y},\bar{z})$ is an optimal solution of \eqref{prob:CP2}. Conversely, let $(\bar{x},\bar{y},\bar{z})$ be an optimal solution of \eqref{prob:CP2}, then  $(\bar{x},\bar{y},\bar{z},\bar{u})$ with $$\bar{u} =(\|\bar{x}\|^2, \|\bar{z}\|^2, (\bar{z}+1)^2,(\bar{z}-1)^2,\|\bar{y}+\bar{x}\|^2,\|\bar{y}-\bar{x}\|^2)$$
	is an optimal solution of \eqref{prob:QP2}.
\end{theorem}

\noindent$\bullet$ Two conditions checked in line \ref{codeline:bdca_dcp1_cond} of Algorithm \ref{alg:BDCA_for_DCP1} are changed to $$\setA_2(V^k)\subset \setA_2(X^k)\quad \text{and} \quad \langle \nabla f_2(V^k),D^k \rangle < 0,$$ where the active set $\setA_2(X)$ of $\setC_2$ at $X:=(x,y,z)$ is defined by {\small$$\setA_2(X) := \{i\in \{1,\ldots,2n+1\}: X_i=0\}\cup\{2n+1+i: (Bx-Ay)_i=0, i\in \{1,\ldots,n\}\}.$$}

\noindent$\bullet$ LineSearch$(V^k,D^k,\bar{\alpha})$ will be either exact or inexact. The exact line search is similarly computed as follows: Let $V^k=(V_x^k,V_y^k,V_z^k)$ and $D^k = (D_x^k,D_y^k,D_z^k)$. Then
\begin{equation}
	\alpha_k = \argmin_{\alpha} \{f_2(V^k + \alpha D^k): \alpha\in \{0,\min\{\bar{\alpha},\bar{\alpha}_k\}\}\cup \mathcal{Z}\},
\end{equation}
where \begin{equation}
	\bar{\alpha}_k = \min\left\{ -(V^k)_i/(D^k)_i, \forall i\in \mathcal{I}^k, -(B V_{x}^{k} -A V_{y}^{k})_j/(B D_{x}^{k} -A D_{y}^{k})_j, \forall j\in \mathcal{J}^k\right\}
\end{equation} 
with $\mathcal{I}^k=\{i\in \{1,\ldots,2n+1\}: (D^k)_i<0\}$ and $\mathcal{J}^k = \{j\in \{1,\ldots,n\}: (B D_{x}^{k} -A D_{y}^{k})_j < 0\}$,
where the assumption $\min \emptyset = \infty$ is adopted and $\mathcal{Z}$ is the set of all real roots of the cubic polynomial
$$q(\alpha)= \frac{\mathrm{d} f_2(V^k + \alpha D^k)}{\mathrm{d} \alpha}  = a_1~\alpha^3 + a_2~\alpha^2 + a_3~\alpha + a_4,$$
with coefficients $a_1,a_2,a_3,a_4$ given in \eqref{eq:coefpolydcp1} by changing $D_{w}^k$ (resp. $V_{w}^k$) to $B D_{x}^{k} -A D_{y}^{k}$ (resp. $B V_{x}^{k} -A V_{y}^{k}$). The Armijo-type inexact line search is the same as for \eqref{prob:dcp1} by changing all $f_1$ to $f_2$.

\subsection{BDCA and DCA for \eqref{prob:dcp3}} BDCA for \eqref{prob:dcp3} is also similar to Algorithm \ref{alg:BDCA_for_DCP1}. The differences are summarized below: 	

\noindent$\bullet$ The initialization in line \ref{codeline:bdca_dcp1_X0} of Algorithm \ref{alg:BDCA_for_DCP1} is changed to
$$X^0 \leftarrow (x^0,y^0,Bx^0-Ay^0).$$

\noindent$\bullet$ The convex subproblem in line \ref{codeline:bdca_dcp1_cp1} of \Cref{alg:BDCA_for_DCP1} is changed to \begin{equation}
	\label{prob:CP3}
	V^k\in \underset{(x,y,w)\in \setC_3}{\argmin} \{ G_3(x,y,w) - \langle (x,y,w), \nabla H_3(x^k,y^k,w^k) \rangle \},\tag{CP3}
\end{equation}
which is a convex QP and can be solved by some efficient QP solvers such as MOSEK, GUROBI and CPLEX.

\noindent$\bullet$ The parameter $\eta$ required in $G_3$ and $H_3$ verifies the inequality $\eta \geq 3.2 + 20 n M^2$ where $M$ is computed by solving the linear program
$$M=\max\{ e^{\top}y: (x,y,w)
\in \setC_3\}$$
via a linear programming solver such as MOSEK, GUROBI and CPLEX.

\noindent$\bullet$ Two conditions checked in line \ref{codeline:bdca_dcp1_cond} of \Cref{alg:BDCA_for_DCP1} are changed to $$\setA_3(V^k)\subset \setA_3(X^k)\quad \text{and} \quad \langle \nabla f_3(V^k),D^k \rangle < 0,$$ where the active set $\setA_3(X)$ of $\setC_3$ at $X:=(x,y,w)$ is defined by $$\setA_3(X) := \{i\in \{1,\ldots,3n\}: X_i=0\}.$$

\noindent$\bullet$ LineSearch$(V^k, D^k, \bar{\alpha})$ is suggested to be the Armijo inexact line search by substituting $f_1$ with $f_3$. The exact line search is too complicated and thus not recommended due to the non-polynomial term $(x^{\top}y)^2/\|x\|^2$ in $f_3$ (In fact, performing the exact line search here amounts to finding all real roots of a quintic equation -- a problem without closed-form formula due to the well-known Abel-Ruffini theorem).

\section{ADCA, InDCA and HDCA for solving \eqref{eq:aeicp}}\label{sec:ADCA&InDCAforAEiCP}
Comparing the convex subproblems required in BDCA and DCA \Cref{alg:BDCA} line \ref{codeline:bdca_P_convexsubprob}:
$$z^{k}\in \argmin\{g(x)-\langle x, \nabla h(x^k) \rangle :x\in \setC\}$$
with ADCA \Cref{alg:ADCA} line \ref{codeline:adca_P_convexsubprob}:
$$x^{k+1}\in \argmin\{g(x)-\langle x, \nabla h(v^k) \rangle :x\in \setC\},$$
InDCA \Cref{alg:InDCA} line \ref{codeline:indca_P_convexsubprob}:
$$x^{k+1}\in \argmin\{g(x)-\langle x, \nabla h(x^k) + \gamma (x^k - x^{k-1}) \rangle:x\in \setC\},$$
as well as HDCA-LI \Cref{alg:HDCA-LI} line \ref{codeline:hdca-li_P_convexsubprob}: 
$$z^{k}\in \argmin\{g(x)-\langle x, \nabla h(x^k) + \gamma (x^k - x^{k-1}) \rangle:x\in \setC\},$$
and HDCA-NI \Cref{alg:HDCA-NI} line \ref{codeline:hdca-ni_P_convexsubprob}:
$$x^{k+1}\in \argmin\{g(x)-\langle x, \nabla h(v^k) + \gamma_k (x^k - x^{k-1}) \rangle:x\in \setC\},$$
we can observe that the difference among these subproblems is only related to the coefficient vector of $x$ in the scalar product $\langle x, \cdot\rangle$. Hence, all of these subproblems can be solved using the same method previously discussed in Section \ref{sec:BDCAforAEiCP}. 

Moreover, in each of these DCA-type algorithms (DCA, BDCA, ADCA, InDCA, HDCA-LI and HDCA-NI), solving the convex subproblem is the most computationally demanding step. Consequently, the computational time per iteration of these algorithms should be fairly comparable. Hence, in numerical simulations, we can focus on comparing the quality of solutions obtained by these methods for each DC formulation with a fixed number of iterations.

\section{Numerical Simulations}\label{sec:Simulations}
In this section, we conduct numerical experiments for 7 DCA-type algorithms: the classical DCA, two BDCA variants (BDCAe with exact line search and BDCAa with Armijo inexact line search), ADCA, InDCA, HDCA-LI, and HDCA-NI, to solve three DC formulations \eqref{prob:dcp1}, \eqref{prob:dcp2} and \eqref{prob:dcp3} of \eqref{eq:aeicp}. Our codes (available on GitHub\footnote{\url{https://github.com/niuyishuai/HDCA}}) are implemented in MATLAB 2022b and tested on a laptop equipped with a 64-bit Windows 10, i7-10870H 2.20GHz CPU, and 32 GB of RAM.

We first compare the performance of these DCA-type algorithms. Then, we compare the best-performing DCA-type algorithm with the state-of-the-art optimization solvers $\IPOPT$ v3.12.9 \cite{Ipopt}, $\KNITRO$ v11.1.0 \cite{knitro} and $\FILTERSD$ v1.0 \cite{filtersd} on three NLP formulations \eqref{prob:nlp1}, \eqref{prob:nlp2} and \eqref{prob:nlp3}. 

\paragraph{\textbf{AEiCP datasets:}} Two datasets are considered.
\begin{itemize}[leftmargin=12pt]
	\item In the first dataset, we generate random AEiCP test problems in a similar way to \cite{judice2022solution}. The matrix $A$ is asymmetric and positive definite, generated by
	$$A = T + \mu I_n $$
	where $T$ is randomly generated with elements uniformly distributed in the interval $[-1,1]$ and $\mu > |\min\{0,\lambda_{\min}(T+T^{\top})\}|.$ These matrices $A$ exhibit good conditioning, as their condition numbers are less than $4$. The matrix $B$ is a symmetric, strictly diagonally dominant matrix with elements of the form
	$$\begin{cases}
		B_{i,i} = 10, i=1,\ldots,n\\
		B_{i,j} = -1, i=1,\ldots,n, j=i+1,\ldots, \min\{i+4,n\},\\
		B_{i,j} = -1, i=1,\ldots,n, j=\max\{1,i-4\},\ldots, i-1.
	\end{cases}$$		
	We generate $10$ random problems for each $n \in \{10, 100, 500\}$ and denote \RAND{n} as the set of $10$ problems.
	\item In the second dataset, the matrix $A$ is taken from the \emph{Matrix Market} repository NEP (Non-Hermitian Eigenvalue Problem) collection. We choose 22 asymmetric matrices with orders $n$ ranging from 62 to 968, where $n$ is indicated in the problem name (e.g., $n=968$ for \verb|rdb968|). These matrices originate from various fields of real applications, and most of them are ill-conditioned (see \url{https://math.nist.gov/MatrixMarket} for more information). The matrix $B$ is set as the identity matrix, and we transform $A$ to $\PD$ by adding $\mu B$ such that $A + \mu B \in \PD$ with $\mu = |\min\{0, \lambda_{\min}(A + A^{\top})\}| + 1$.
\end{itemize}

\paragraph{\textbf{Experimental setup:}} The setups for DCA-type algorithms and the compared optimization solvers are summarized below:
\begin{itemize}[leftmargin=12pt]
	\item Initialization:
	\begin{itemize}[leftmargin=12pt]
		\item For DCA-type algorithms, we set $x^0 = \verb|rand|(n,1)$ and normalize $x^0$ by $x^0 = x^0/\verb|sum|(x^0)$ to obtain a vector on the simplex $\Omega$. Then we set $y^0 = \verb|rand|(n,1)$, $w^0 = Bx^0 - Ay^0$ and $z^0 = \verb|sum|(y^0)$. Note that all compared DCA-type algorithms use the same initial point for fairness.
		\item The optimization solvers IPOPT, KNITRO, FILTERSD employ the same initial point as in DCA-type algorithms.
	\end{itemize}
	\item Other settings: 
	\begin{itemize}[leftmargin=12pt]
		\item We employ MOSEK to solve the convex subproblems and the linear problem \eqref{prob:lpforM} for computing $M$. MOSEK's termination is achieved by setting the tolerances \texttt{MSK\_DPAR\_INTPNT\_QO\_TOL\_REL\_GAP}, \texttt{MSK\_DPAR\_INTPNT\_QO\_TOL\_PFEAS}, and \texttt{MSK\_DPAR\_INTPNT\_QO\_TOL\_DFEAS} to $10^{-8}$.
		\item The parameter $\eta=3.2 + 20 n M^2$ is used in \eqref{prob:dcp3}. 
		\item The parameter $q=10$ is used for both ADCA and HDCA-NI.
		\item An additional strongly convex regularization, $\frac{\rho}{2}\|\cdot\|^2$, is introduced to each DC formulation. We set $\rho=0.1$ by default, ensuring that $\rho\leq \min\{\rho_g,\rho_h\}$.
		\item The parameter $\bar{\beta}=0.99$ for HDCA-NI, and thus $\delta = (1-\bar{\beta}^2)\rho/2 \approx 9.95\times 10^{-4}$.
		\item The stepsize for the line search is upper bounded by $\bar{\alpha} = 10$.
		\item The stepsize for the inertial force is $\gamma = \rho$ for InDCA, $\gamma = 2\rho/(1+(1+\bar{\alpha})^2)$ for HDCA-LI and $\gamma_{k} = (2\rho(1-\beta_k^2) - 4 \delta)/(3-\beta_k^2)$ for HDCA-NI. 
	\end{itemize}
\end{itemize}

\subsection{Numerical results of DCA-type algorithms}\label{subsec:resultsfordcas}
\paragraph{\textbf{Tests on the \RAND{n} dataset:}}
We terminate all DCA-type algorithms (DCA, BDCAe, BDCAa, InDCA, ADCA, HDCA-LI and HDCA-NI) for a fixed number of iterations $\MaxIT$ (say $200$), and evaluate each DC formulation and DCA-type algorithm by comparing the trend of the average objective value for $10$ test problems in each dataset \texttt{RAND(n)} where $n =10, 100, 500$. The numerical results are shown in \Cref{fig:numericalresultsDCP1_rand2,fig:numericalresultsDCP2_rand2,fig:numericalresultsDCP3_rand2}. The left column depicts the trend of the average objective value versus the number of iterations, while the right column presents the average CPU time (in seconds) for each DCA-type algorithm. Note that the line search in HDCA-LI for \eqref{prob:dcp1} and \eqref{prob:dcp2} is exact, whereas the line search in HDCA-LI for \eqref{prob:dcp3} is inexact. Next, we summarize some observations as follows:

\begin{itemize}[leftmargin=12pt]
	\item For the model \eqref{prob:dcp1}, we observe from \Cref{fig:numericalresultsDCP1_rand2} that:
	\begin{figure}[ht!]
		\subfigure[\texttt{RAND(10)}]{ 
			\begin{minipage}{.50\textwidth}
				\centering
				\includegraphics[width=\linewidth]{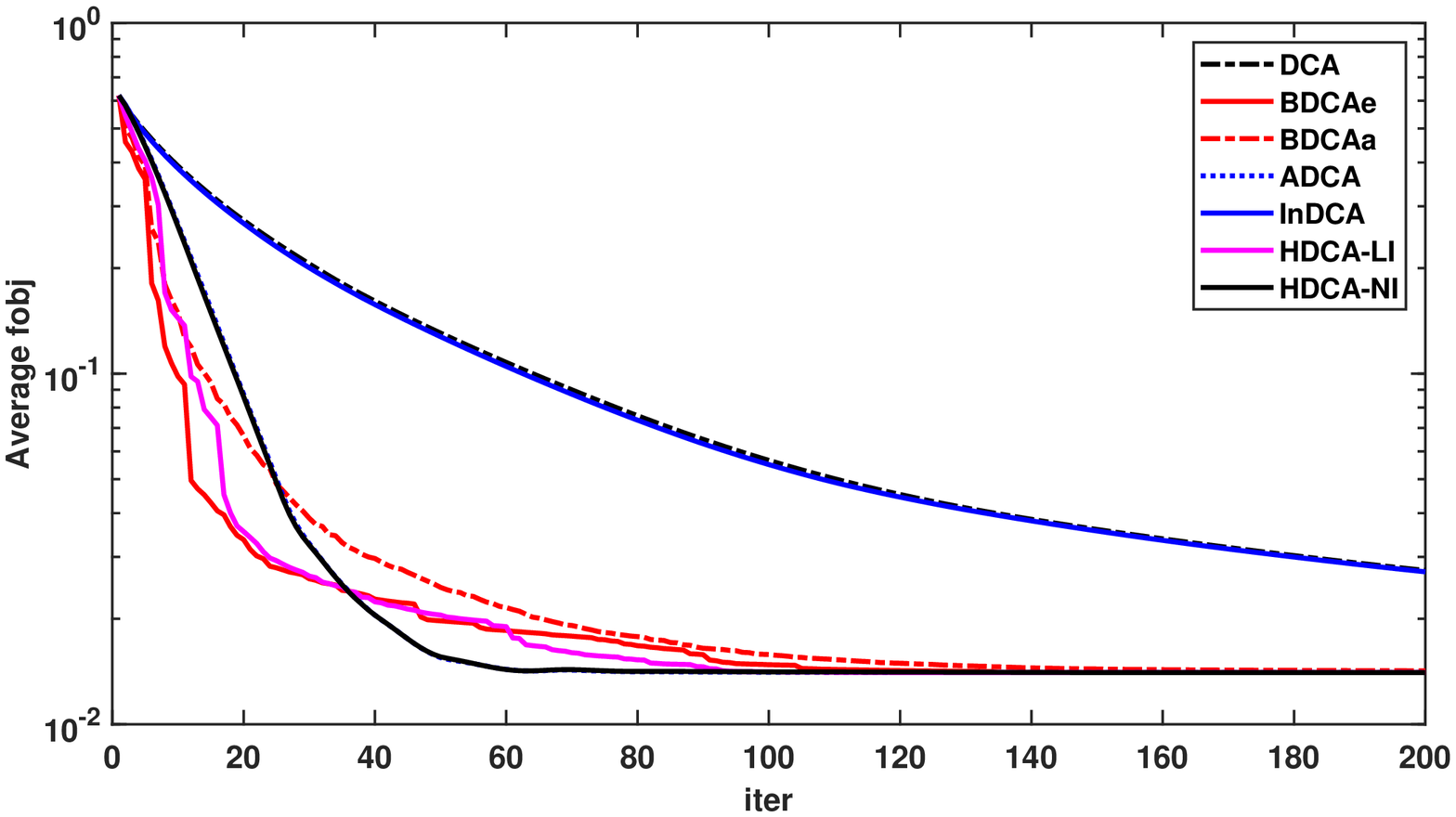}
			\end{minipage}\hfill
			\begin{minipage}{.30\textwidth}
				\centering
				\includegraphics[width=\linewidth]{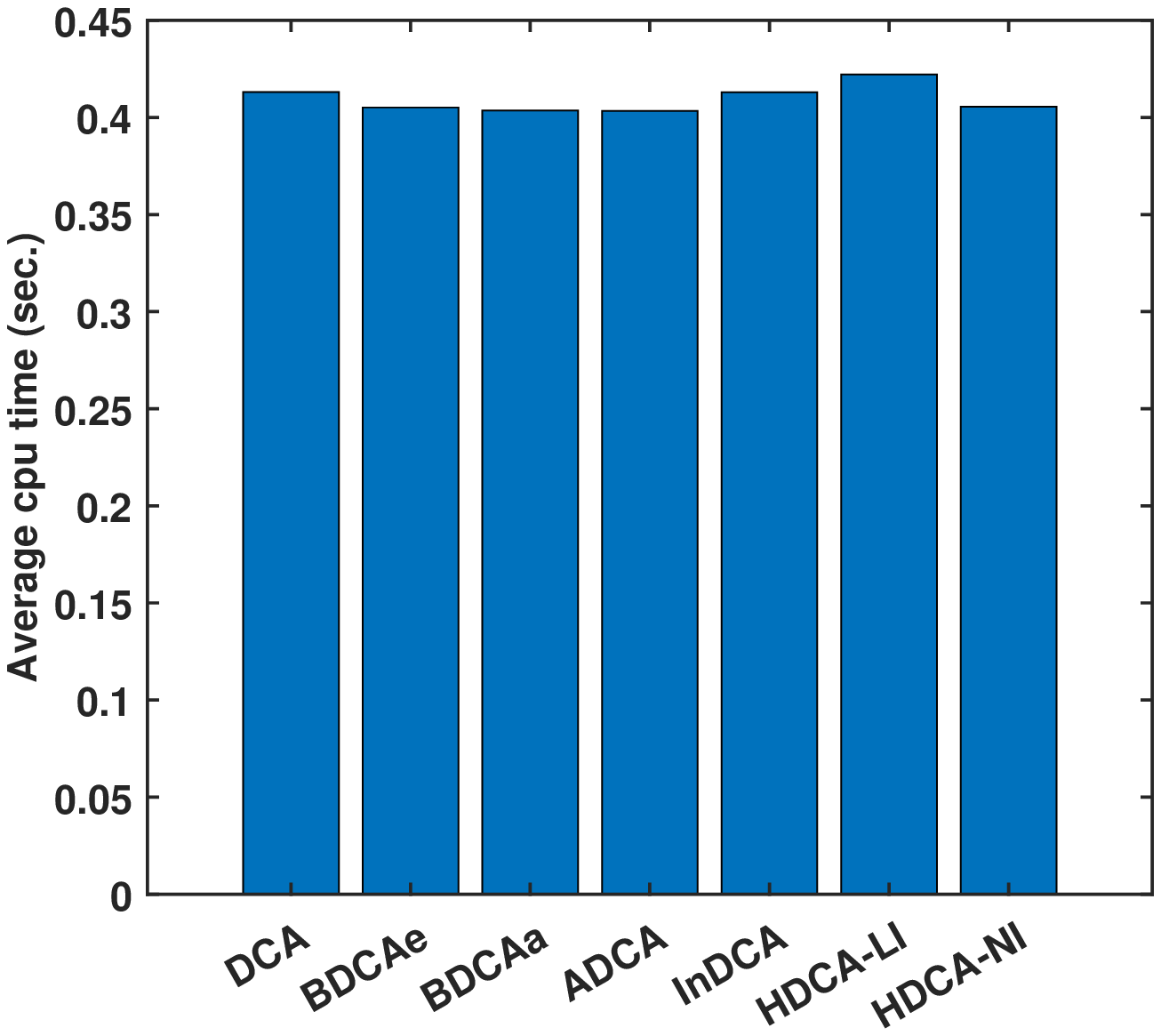}
			\end{minipage}
		}
		
		\subfigure[\texttt{RAND(100)}]{ 
			\begin{minipage}{.50\textwidth}
				\centering
				\includegraphics[width=\linewidth]{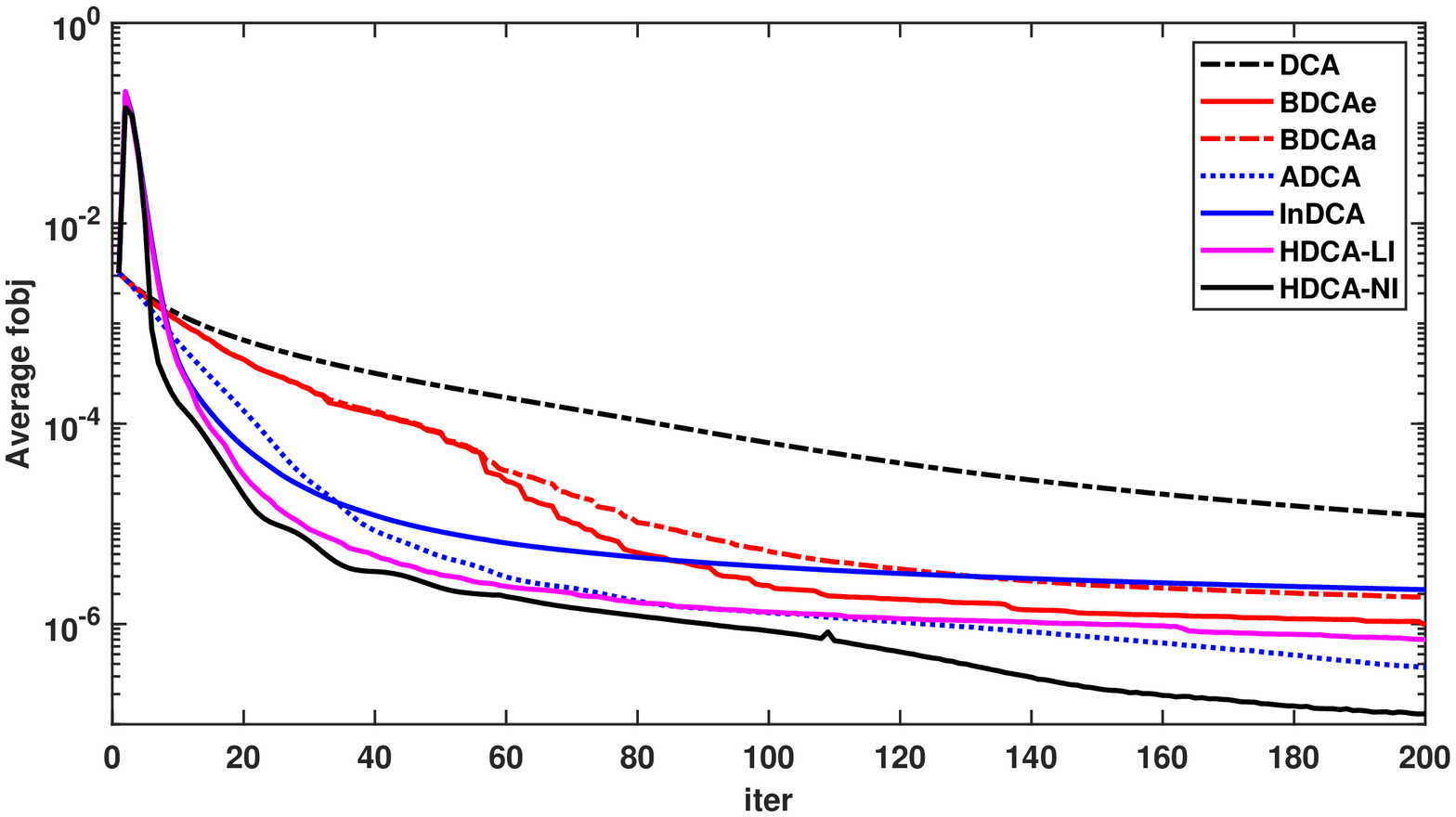}
			\end{minipage}\hfill
			\begin{minipage}{.30\textwidth}
				\centering
				\includegraphics[width=\linewidth]{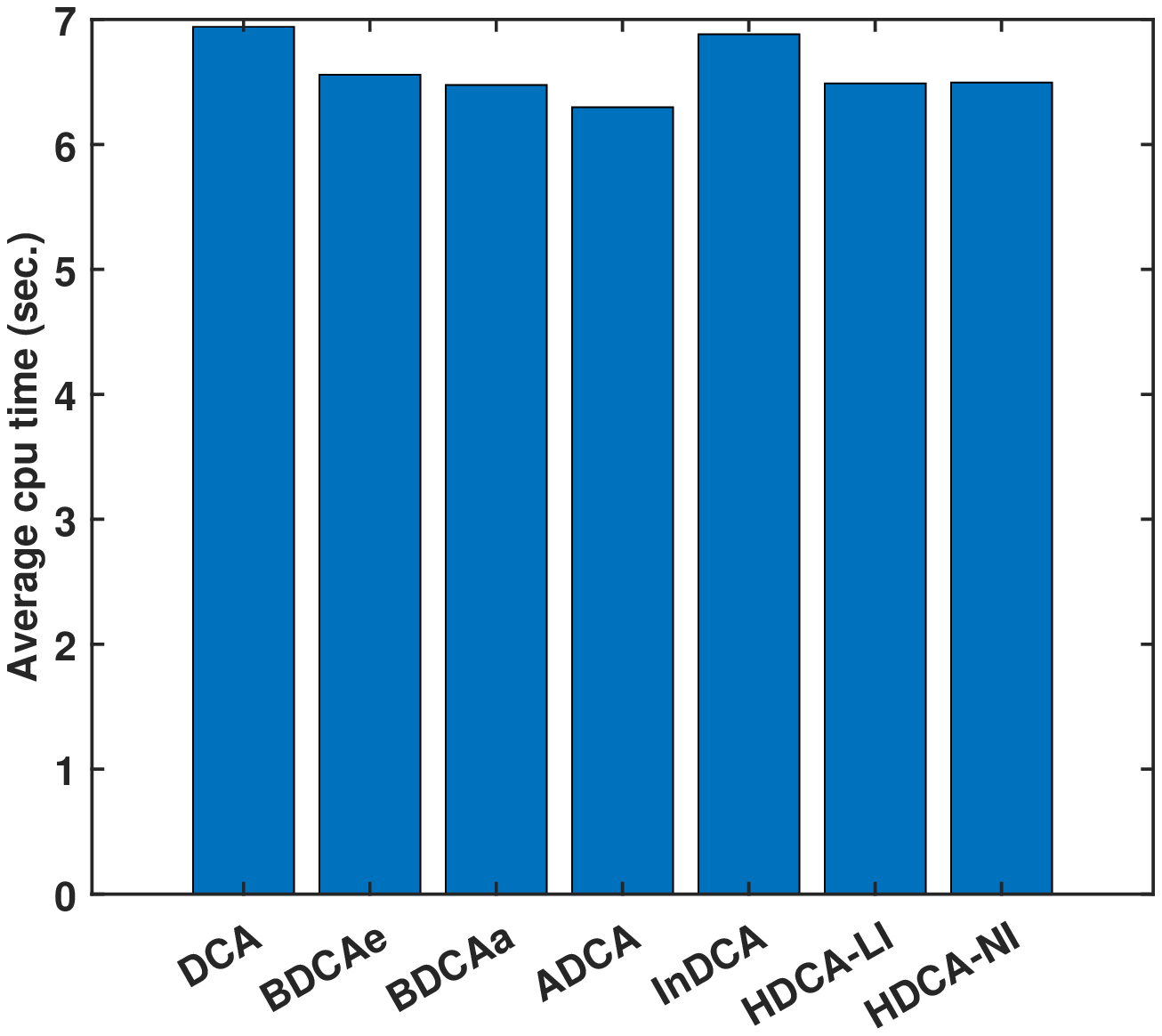}
			\end{minipage}
		}
		
		\subfigure[\texttt{RAND(500)}]{ 
			\begin{minipage}{.50\textwidth}
				\centering
				\includegraphics[width=\linewidth]{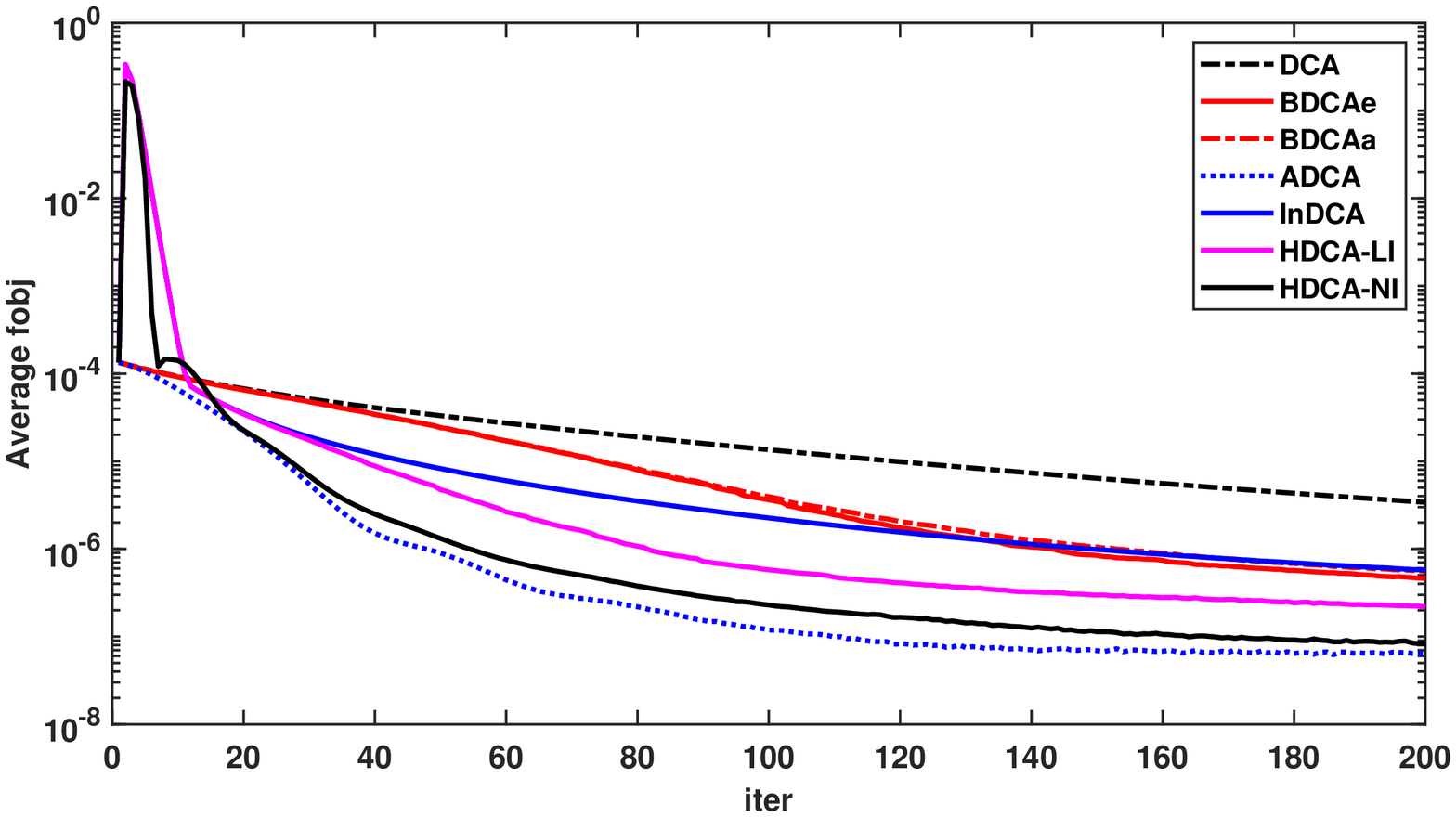}
			\end{minipage}\hfill
			\begin{minipage}{.30\textwidth}
				\centering
				\includegraphics[width=\linewidth]{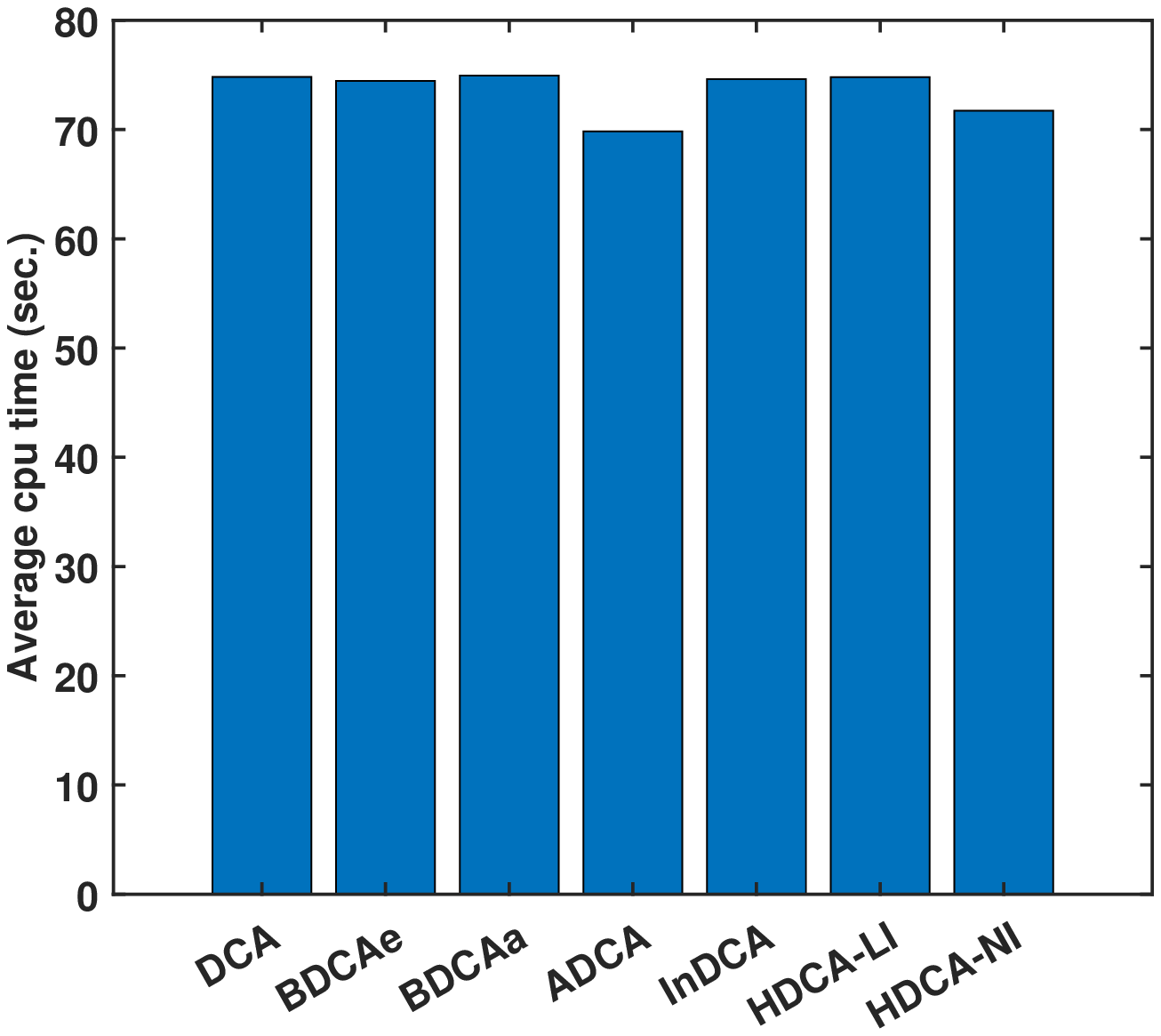}
			\end{minipage}
		}
		
		\caption{Numerical results of DCA, BDCAe, BDCAa, ADCA, InDCA, HDCA-LI and HDCA-NI for solving \eqref{prob:dcp1} on the test datasets \texttt{RAND(n)} with $n\in\{10,100,500\}$.}
		\label{fig:numericalresultsDCP1_rand2}
	\end{figure}
	\begin{itemize}[leftmargin=12pt]
		\item The accelerated variants of DCA (HDCA-NI, HDCA-LI, BDCAe, BDCAa, ADCA, InDCA) consistently outperform the classical DCA.
		\item The hybrid method HDCA-NI yields the best numerical result in terms of the average objective value for the majority of tested cases, with the exception on the dataset \RAND{500}, where ADCA emerges as the best performer.
		\item HDCA-LI secures the second-best average objective value for the dataset \RAND{10}, while ADCA holds the position for the dataset \RAND{100}.
		\item For accelerated DCA without hybridization (i.e., BDCAe, BDCAa, InDCA and ADCA), it appears that ADCA outperforms BDCAe, which in turn outperforms BDCAa, while InDCA ranks as the second-worst algorithm.
		\item As anticipated, the average CPU time per iteration remains nearly identical for all tested DCA-type algorithms.
	\end{itemize}
	\begin{figure}[ht!]
		\subfigure[\texttt{RAND(10)}]{ 
			\begin{minipage}{.50\textwidth}
				\centering
				\includegraphics[width=\linewidth]{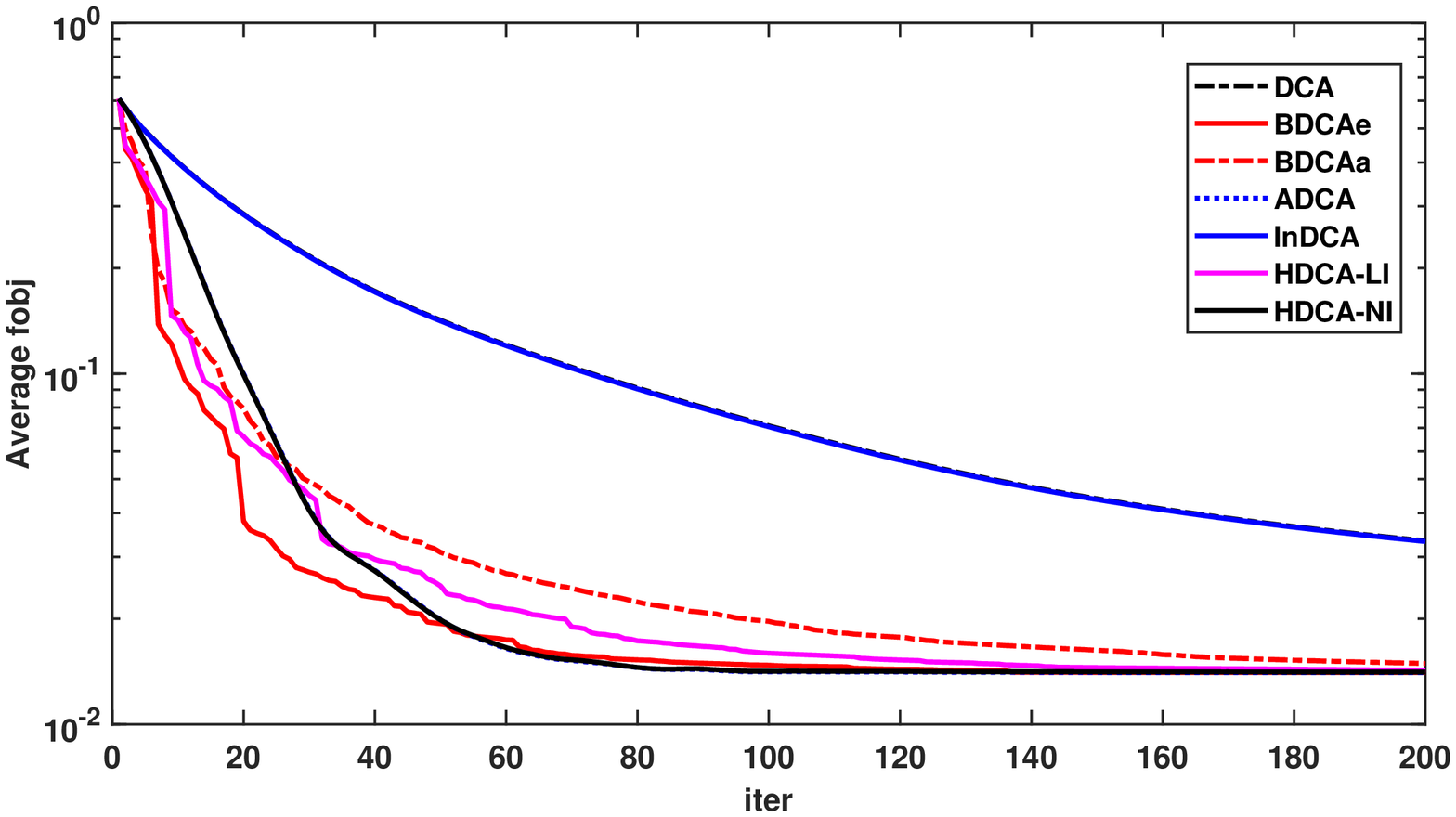}
			\end{minipage}\hfill
			\begin{minipage}{.30\textwidth}
				\centering
				\includegraphics[width=\linewidth]{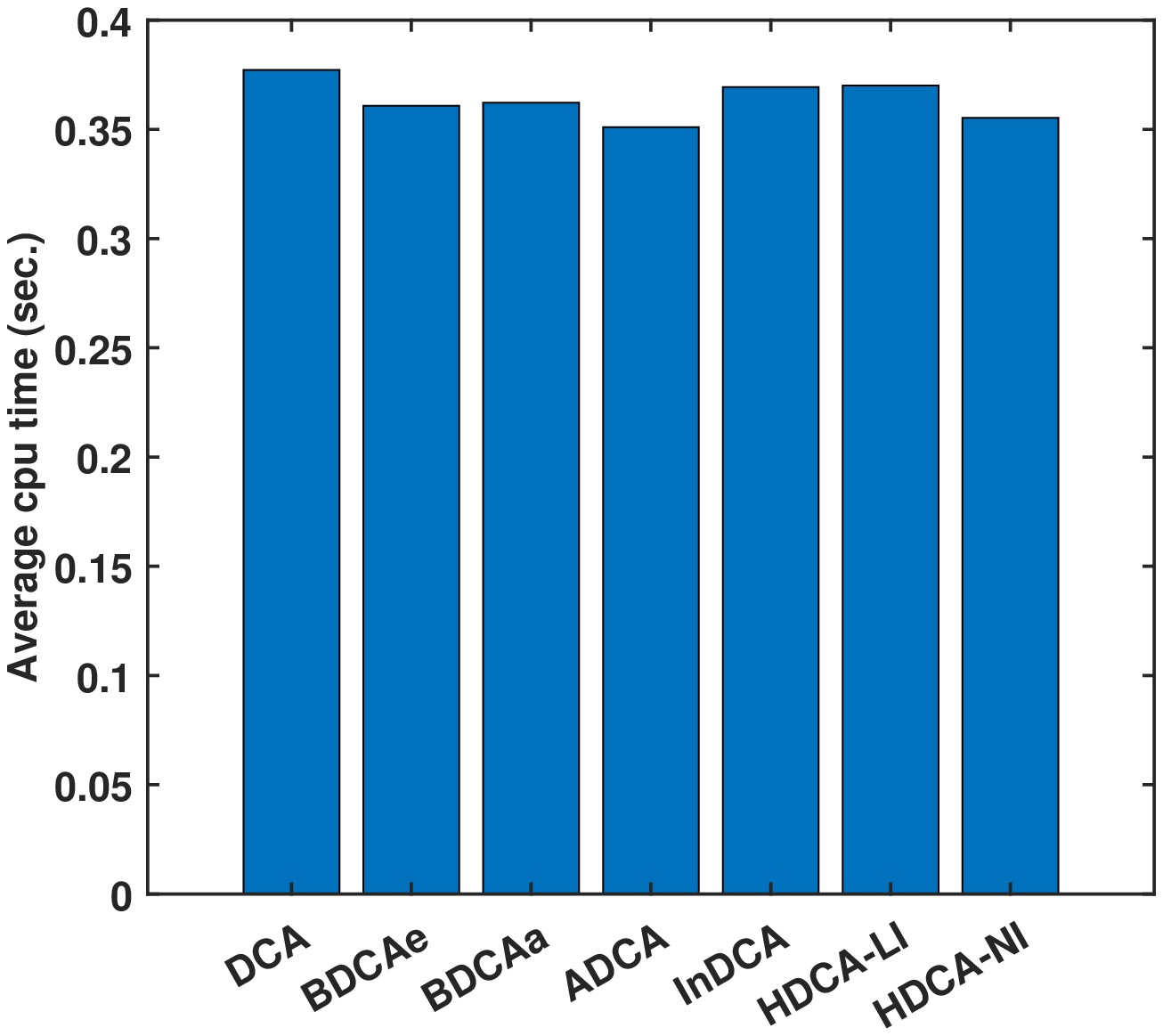}
			\end{minipage}
		}
		
		\subfigure[\texttt{RAND(100)}]{ 
			\begin{minipage}{.50\textwidth}
				\centering
				\includegraphics[width=\linewidth]{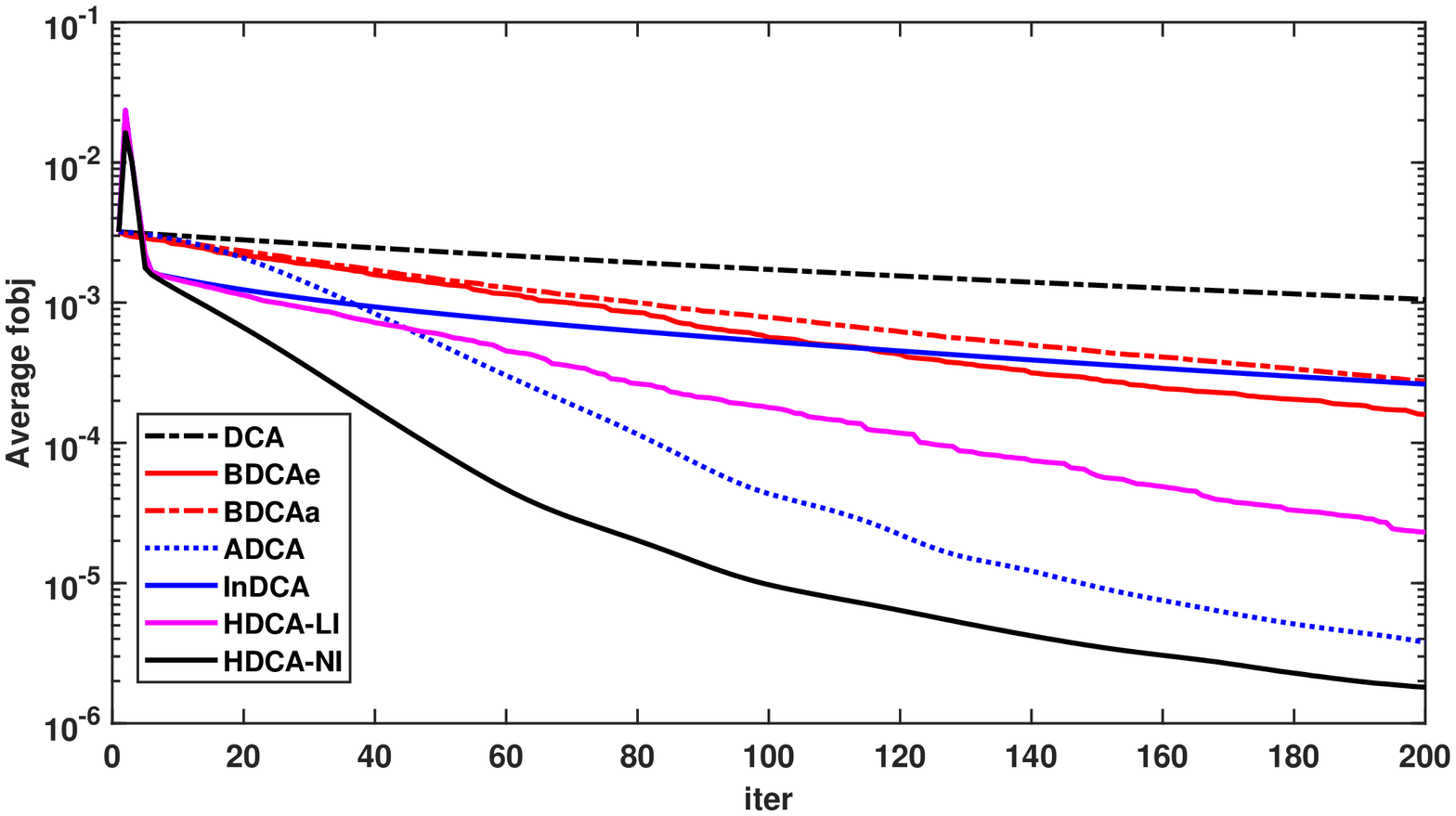}
			\end{minipage}\hfill
			\begin{minipage}{.30\textwidth}
				\centering
				\includegraphics[width=\linewidth]{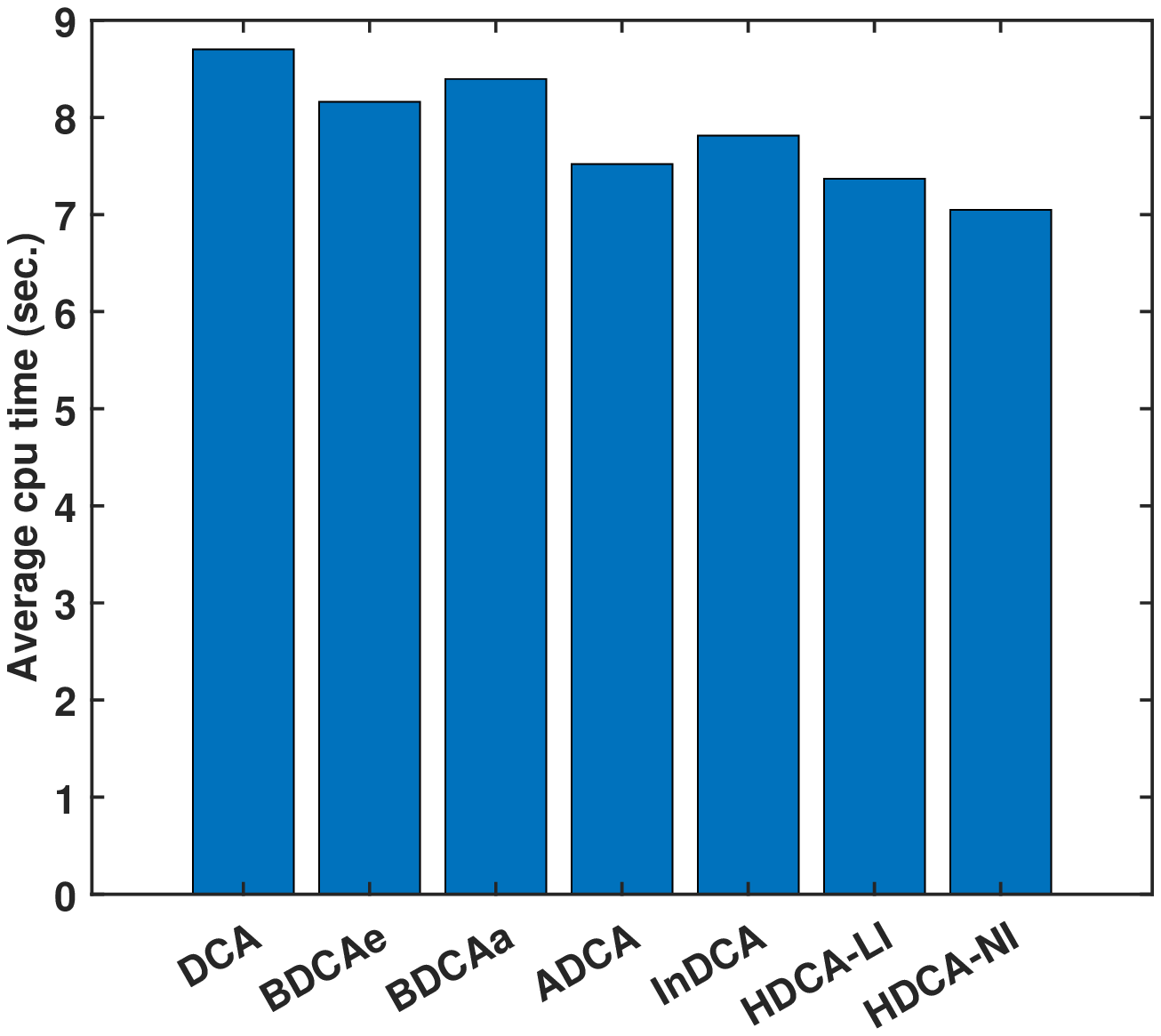}
			\end{minipage}
		}
		
		\subfigure[\texttt{RAND(500)}]{ 
			\begin{minipage}{.50\textwidth}
				\centering
				\includegraphics[width=\linewidth]{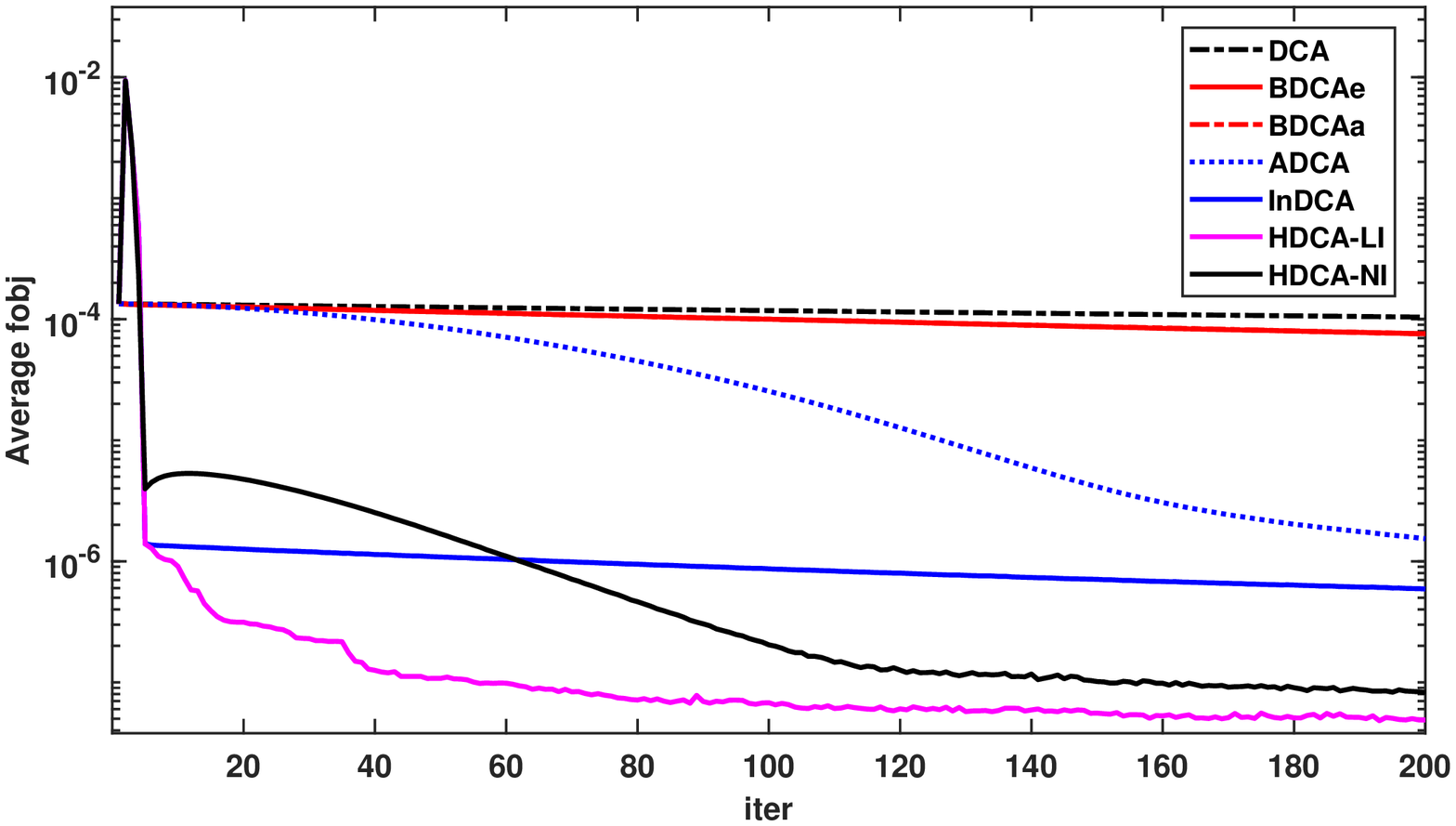}
			\end{minipage}\hfill
			\begin{minipage}{.30\textwidth}
				\centering
				\includegraphics[width=\linewidth]{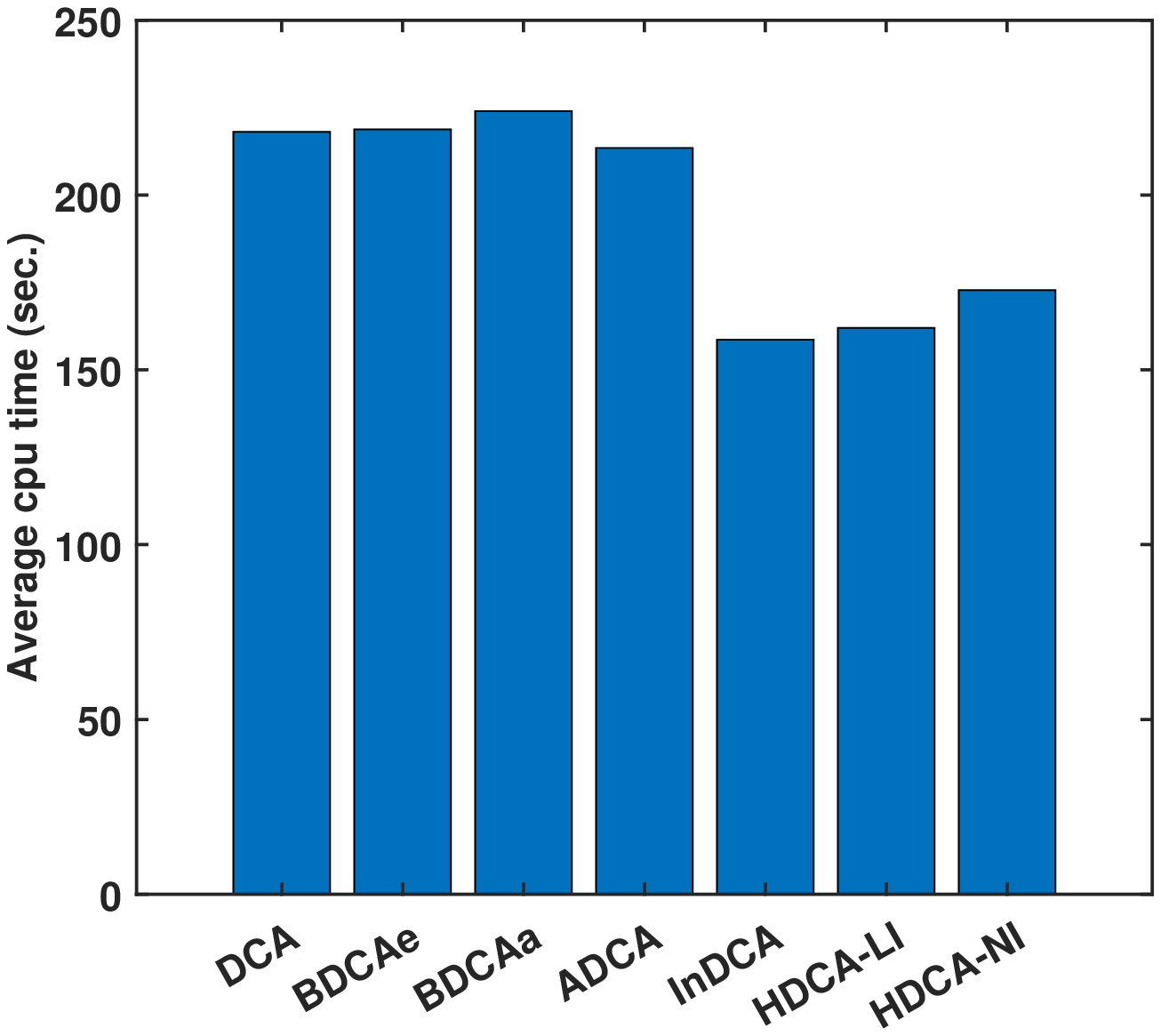}
			\end{minipage}
		}
		
		\caption{Numerical results of DCA, BDCAe, BDCAa, ADCA, InDCA, HDCA-LI and HDCA-NI for solving \eqref{prob:dcp2} on the datasets \texttt{RAND(n)} with $n\in\{10,100,500\}$.}
		\label{fig:numericalresultsDCP2_rand2}
	\end{figure}
	\item For \eqref{prob:dcp2}, we once again observe in \Cref{fig:numericalresultsDCP2_rand2} that all accelerated variants of DCA outshine the classical DCA in terms of the average objective value. Among the top performers, HDCA-NI, HDCA-LI and ADCA consistently stand out, then followed by BDCAe and BDCAa. The performance of InDCA, however, varies significantly, as illustrated by the contrasting performance of InDCA on \RAND{10} and \RAND{500}, where it performs remarkably on \RAND{500} but merely matches the classical DCA on \RAND{10}. 
	\item \begin{figure}[ht!]
		\subfigure[\texttt{RAND(10)}]{ 
			\begin{minipage}{.50\textwidth}
				\centering
				\includegraphics[width=\linewidth]{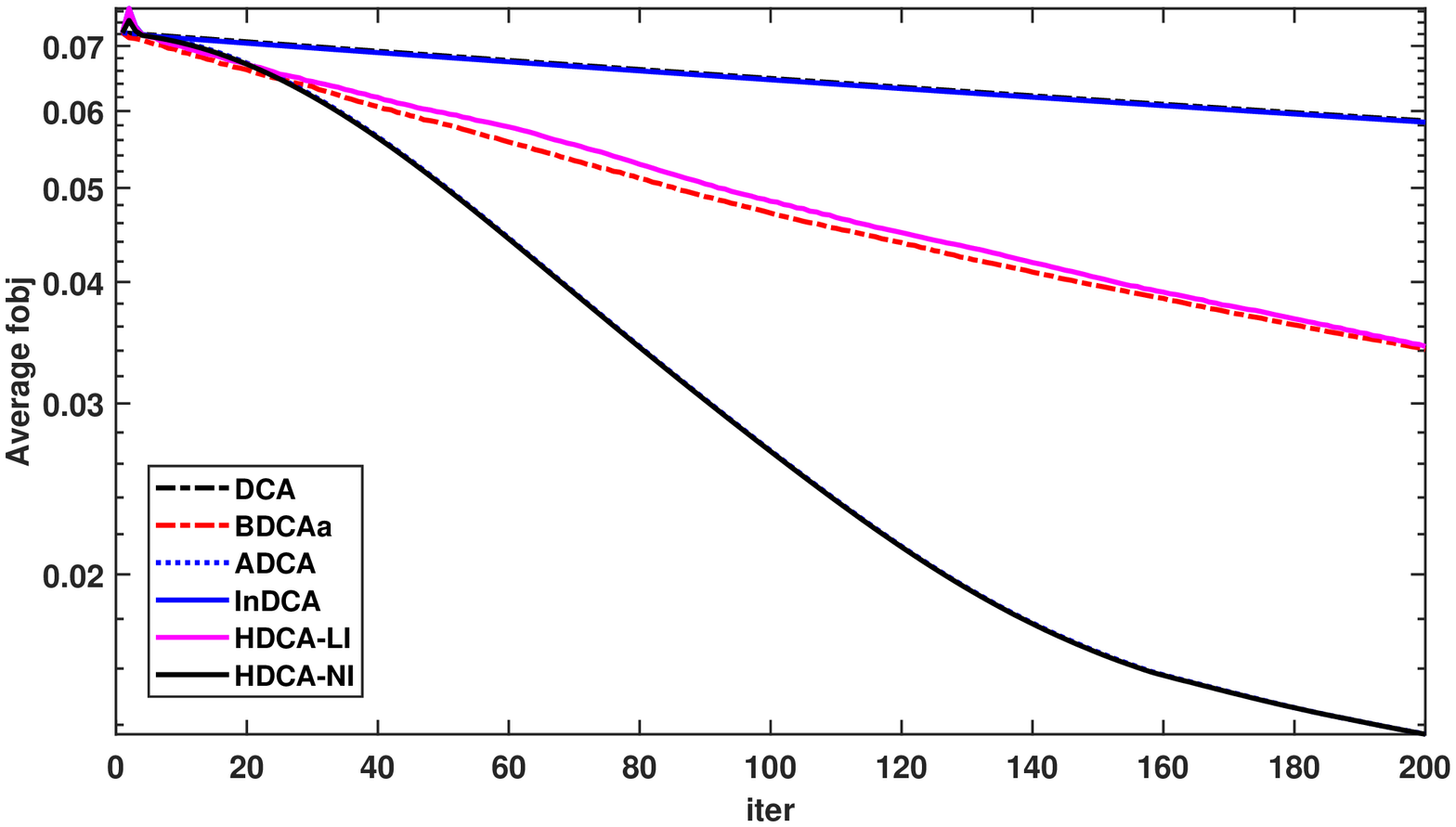}
			\end{minipage}\hfill
			\begin{minipage}{.30\textwidth}
				\centering
				\includegraphics[width=\linewidth]{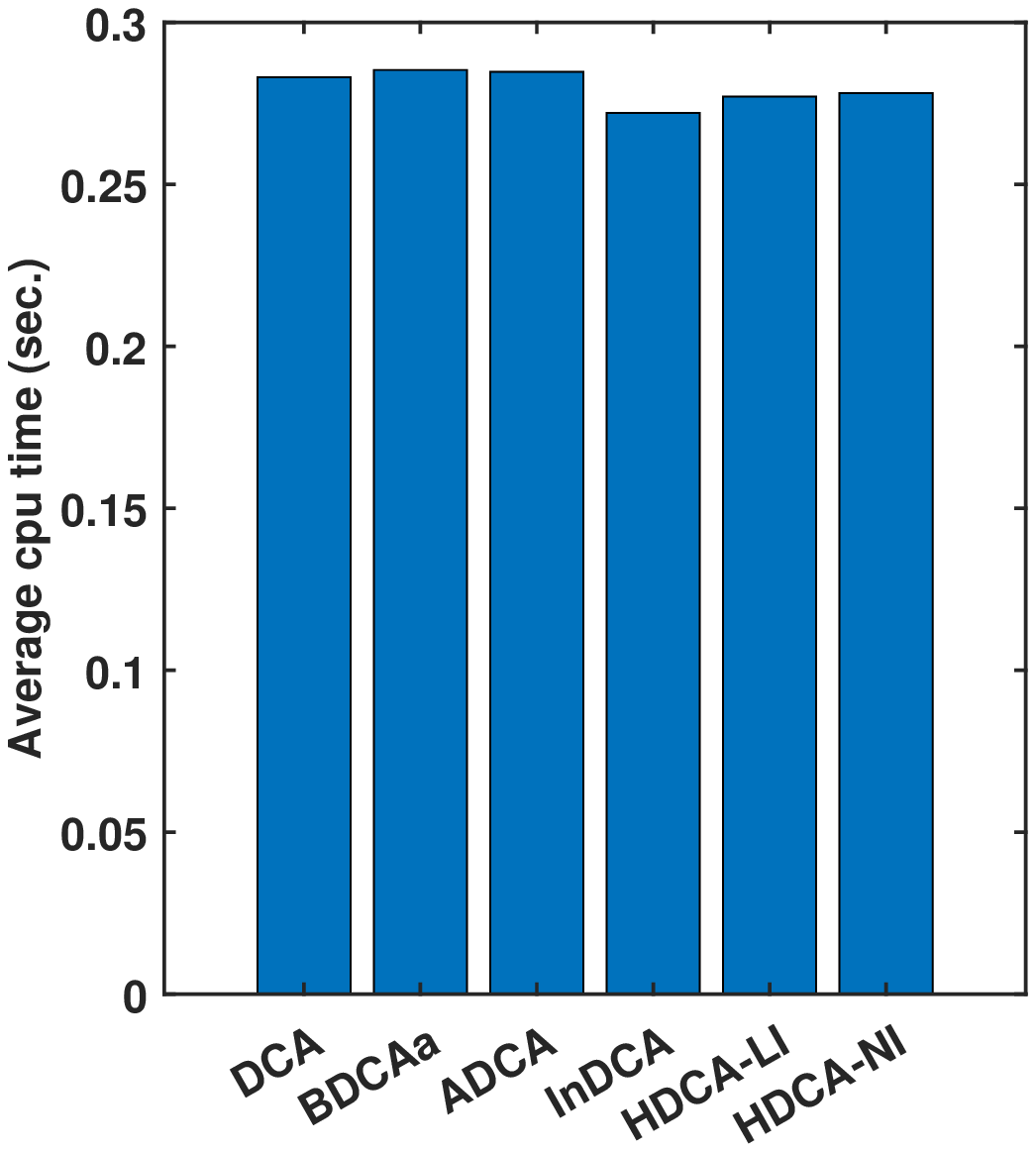}
			\end{minipage}
		}
		
		\subfigure[\texttt{RAND(100)}]{ 
			\begin{minipage}{.50\textwidth}
				\centering
				\includegraphics[width=\linewidth]{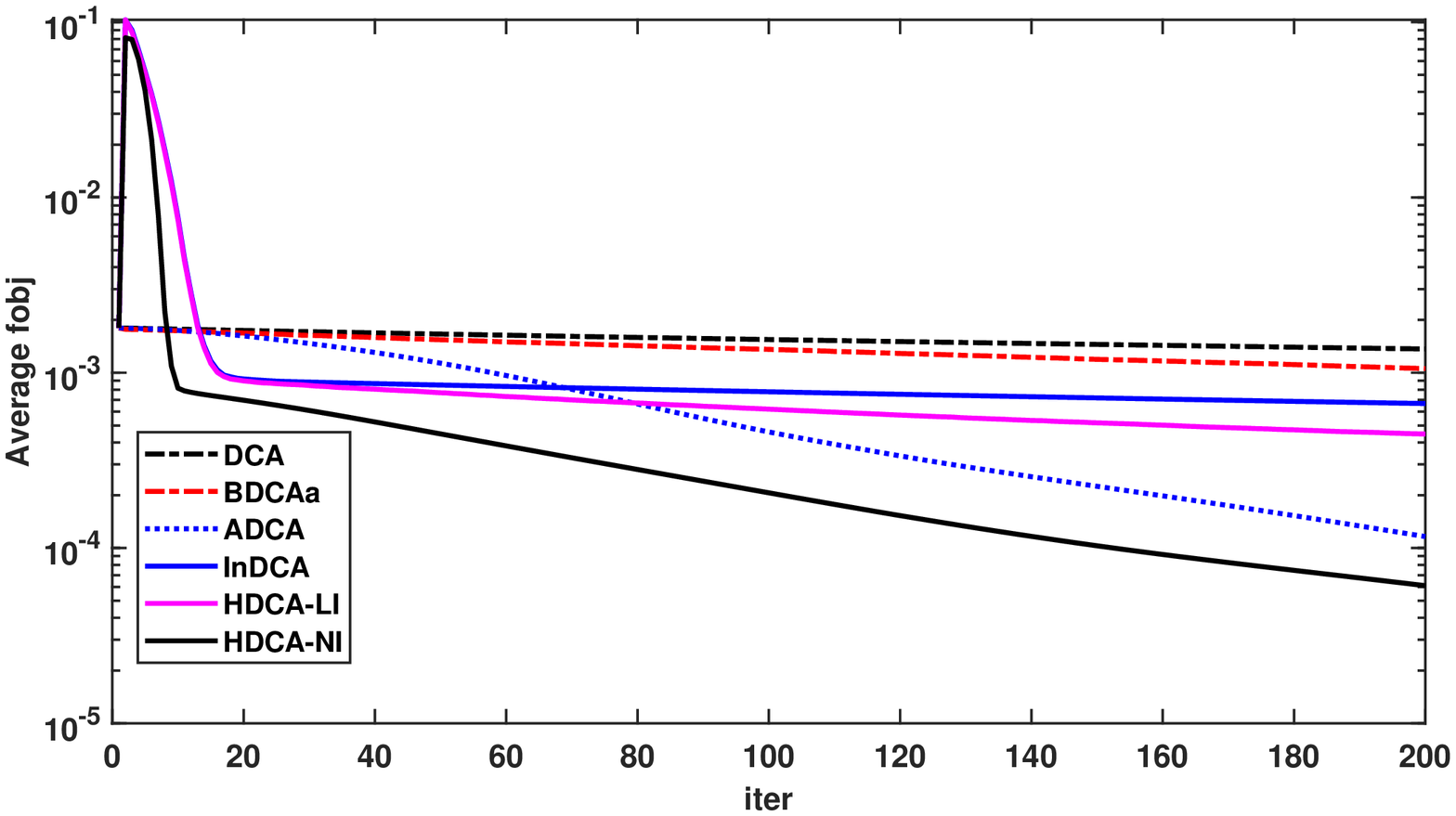}
			\end{minipage}\hfill
			\begin{minipage}{.30\textwidth}
				\centering
				\includegraphics[width=\linewidth]{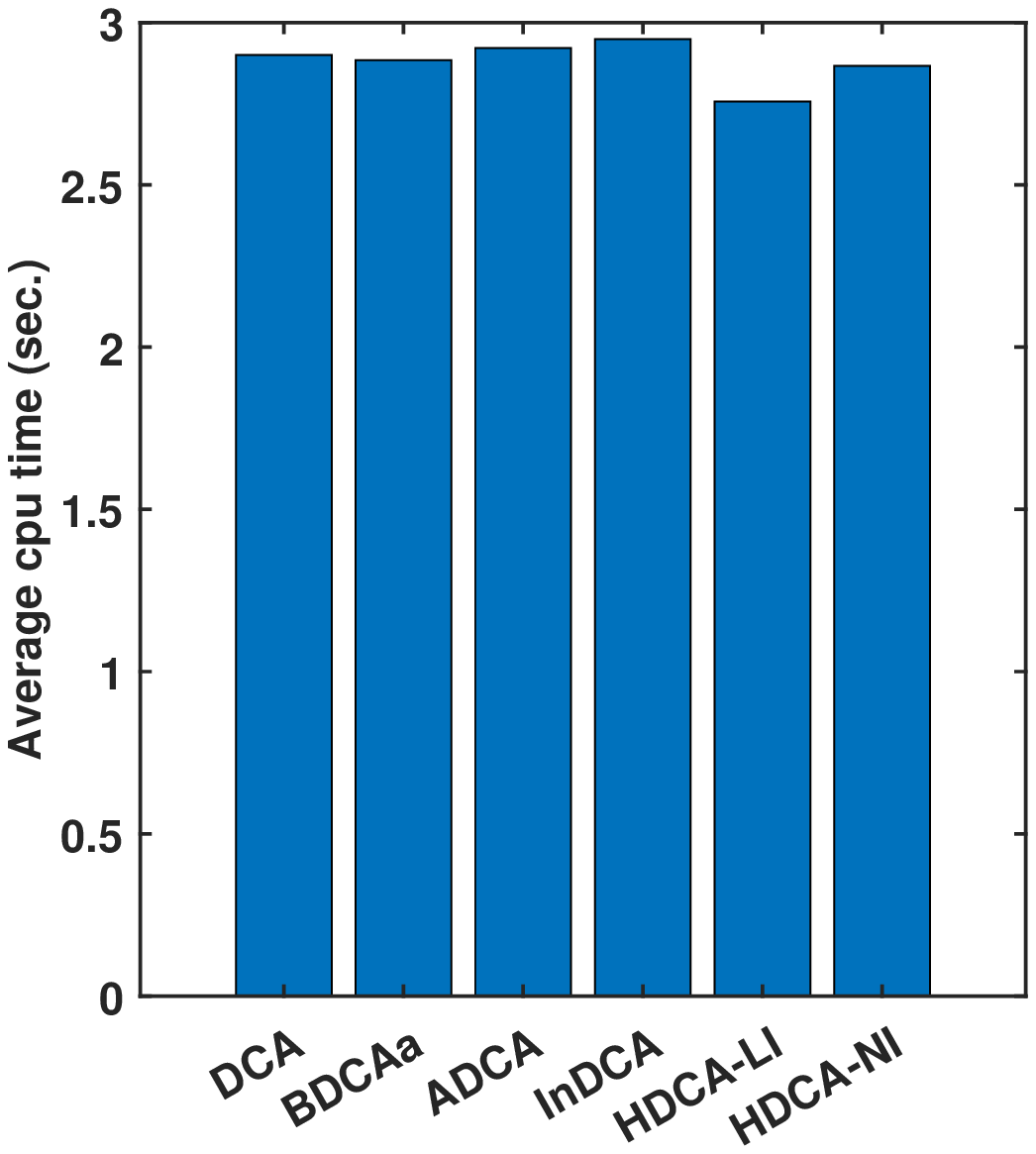}
			\end{minipage}
		}
		
		\subfigure[\texttt{RAND(500)}]{ 
			\begin{minipage}{.50\textwidth}
				\centering
				\includegraphics[width=\linewidth]{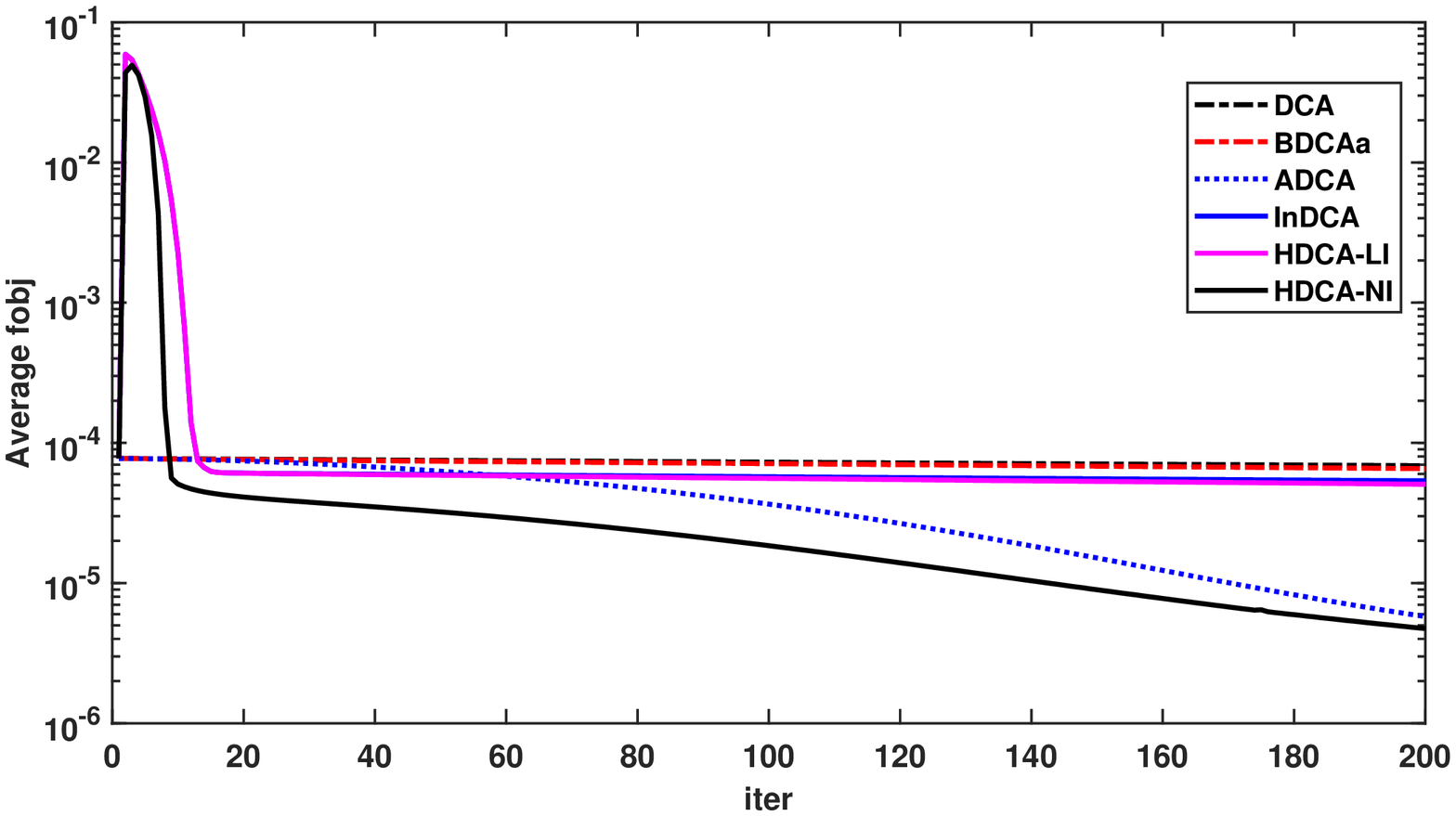}
			\end{minipage}\hfill
			\begin{minipage}{.30\textwidth}
				\centering
				\includegraphics[width=\linewidth]{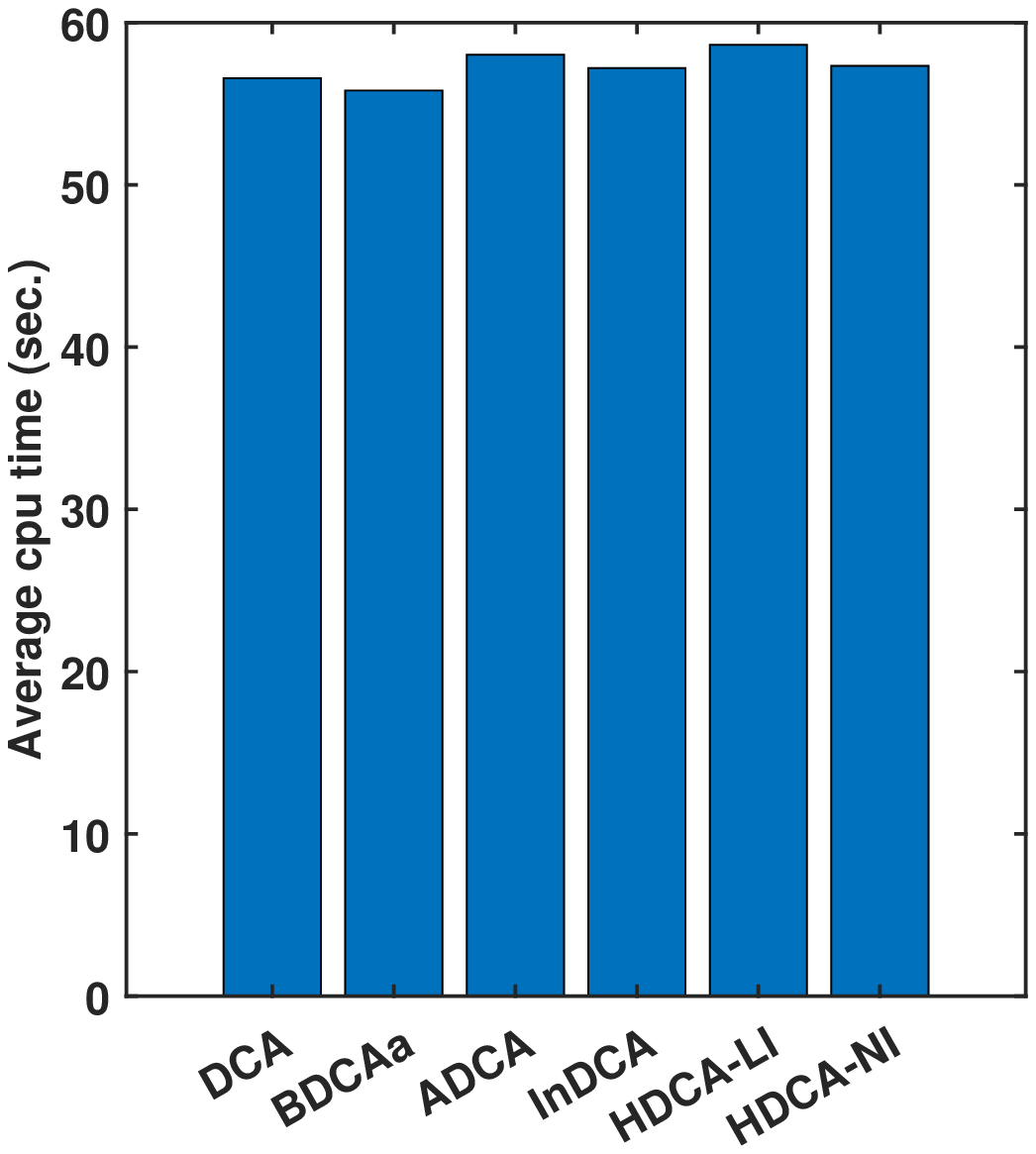}
			\end{minipage}
		}
		
		\caption{Numerical results of DCA, BDCAa, ADCA, InDCA, HDCA-LI and HDCA-NI for solving \eqref{prob:dcp3} on the datasets \texttt{RAND(n)} with $n\in\{10,100,500\}$.}
		\label{fig:numericalresultsDCP3_rand2}
	\end{figure} For \eqref{prob:dcp3}, we observe a similar result in \Cref{fig:numericalresultsDCP3_rand2}, with all accelerated DCA variants outperforming the classical DCA. The best result is consistently provided by HDCA-NI, followed by ADCA, HDCA-LI, BDCAe, and InDCA. It is worth noting that the choice of parameter $\bar{\alpha}=10$ is a conservative setting; it seems that increasing this value often leads to better numerical results for BDCAe and BDCAa. \Cref{fig:alphabarvsfobj} illustrates an example of the impact of $\bar{\alpha}$ on the computed objective value for BDCAa when solving \eqref{prob:dcp3} within 200 iterations on the dataset \RAND{100}. We observe that the optimal result is achieved with $\bar{\alpha}=200$. Beyond this value, the performance of BDCAa plateaus, as $\bar{\alpha}_k < 200$ for all $k \in \N$. Notably, this observation is consistent across all DC formulations.
	\begin{figure}[tbhp]
		\centering
		\includegraphics[width=0.6\textwidth]{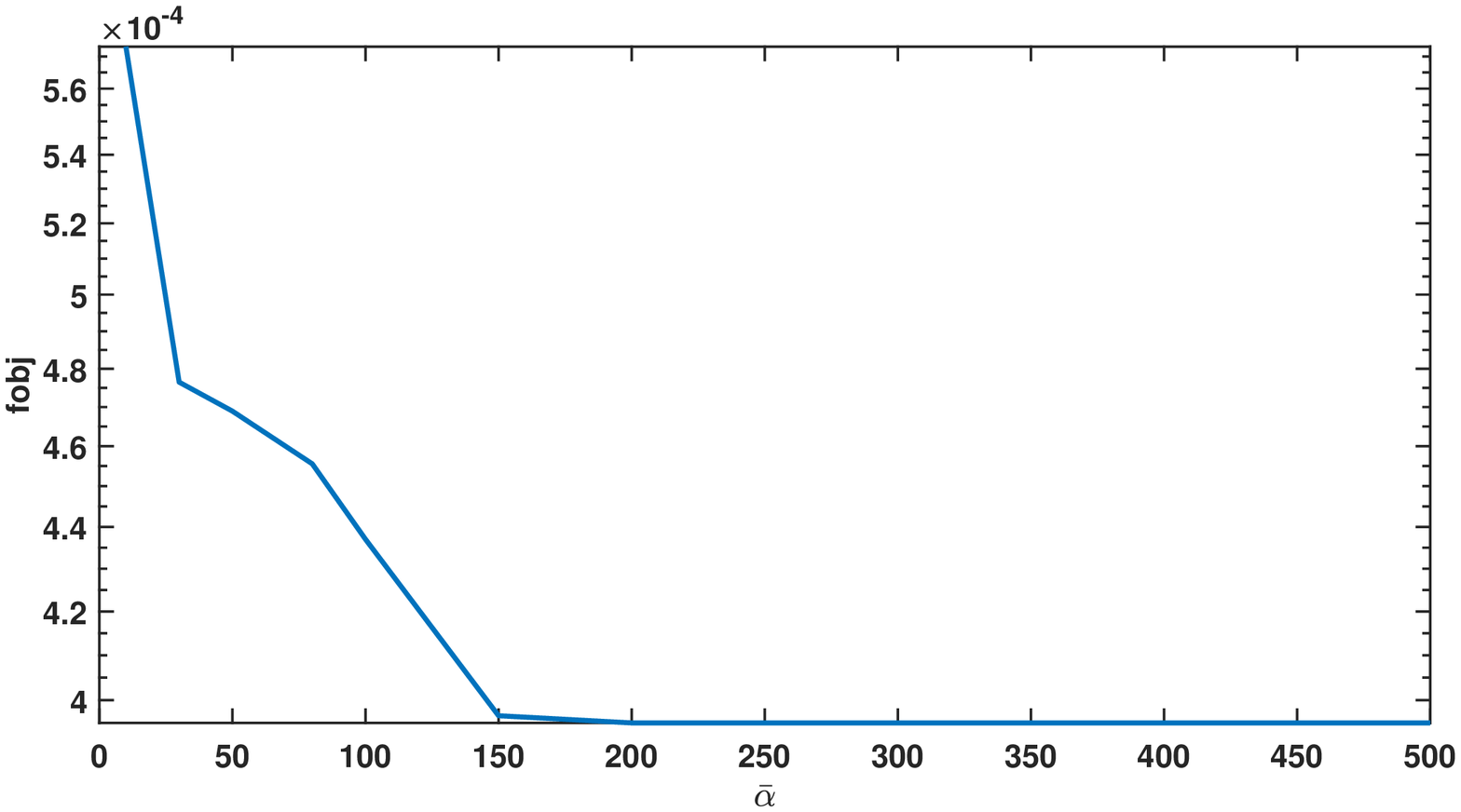}
		\caption{Example to illustrate the influence of $\bar{\alpha}$ on the objective value for BDCAa.}
		\label{fig:alphabarvsfobj}
	\end{figure}
\end{itemize}

\paragraph{\textbf{Tests on the NEP dataset:}}
The numerical results of the DCA-type algorithms on the NEP dataset are summarized in \Cref{tab:performcomp_allDCAs_dcp1,tab:performcomp_allDCAs_dcp2,tab:performcomp_allDCAs_dcp3} for the three DC formulations. In these tables, we adopt the following notations: 
\begin{itemize}[leftmargin=12pt]
	\item cond($A$) - the condition number of the matrix $A$;
	\item $f$ - the objective value; The smaller the value of $f$, the higher the quality of the computed solution;
	\item $c$ - the feasibility measure of the computed solution defined by 
	$$c := -\log(\|[x]_-\| + \|[ w]_-\| + |w^{\top}x|),$$
	where $[x]_-=\min\{x,0\}$ and $w = \frac{x^{\top}Ax}{x^{\top}Bx} Bx - Ax$. The larger the value of $c$, the higher the quality of the computed solution;
	\item avg - average results regarding to $f$ and $c$. 
\end{itemize}
\begin{remark}
	The reasoning behind reporting both $f$ and $c$ is that a small value in $f$ may not necessarily imply a large value in $c$, especially if either the matrix $A$ or $B$ is ill-conditioned. For example, in \eqref{prob:dcp1}, let's assume that we obtain an approximate solution $x + \Delta x$ for the exact solution $x$. When the value of $f$ is small, then the error $\Delta x$ is also small. We get from the relation  
	$$w = \frac{x^{\top}Ax}{x^{\top}Bx} Bx - Ax$$ for the exact solution $(x,w)$ that, if $A$ is ill-conditioned and $B$ is well-conditioned, then $B(x+\Delta x)\approx B x$ and
	$$\frac{(x+\Delta x)^{\top}A(x+\Delta x)}{(x+\Delta x)^{\top}B(x+\Delta x)} B(x+\Delta x) - A(x+\Delta x)\approx w + \underbrace{\frac{2x^{\top}A\Delta x + \Delta x^{\top}A\Delta x}{x^{\top}B x} Bx - A \Delta x}_{\approx \Delta w}.$$
	Consequently, when dealing with an ill-conditioned $A$, we will get a large $A\Delta x$, indicating that $\Delta w$ could be substantial as well. This, in turn, could result in a considerable $\|[w + \Delta w]_{-}\|$ and ultimately lead to a small $c$. 
\end{remark}

\begin{table}[tbhp]
	\caption{Numerical results of DCA, BDCAe, BDCAa, ADCA, InDCA, HDCA-LI and HDCA-NI for solving \eqref{prob:dcp1} on the NEP dataset with \texttt{MaxIT}=200.}
	\label{tab:performcomp_allDCAs_dcp1}
	\centering
	\resizebox{\columnwidth}{!}{
		\begin{tabular}{l|c|cc|cc|cc|cc|cc|cc|cc} 
			\toprule
			\multirow{2}{*}{Prob} & \multirow{2}{*}{cond($A$)} & \multicolumn{2}{c|}{DCA} & \multicolumn{2}{c|}{BDCAe} & \multicolumn{2}{c|}{BDCAa} & \multicolumn{2}{c|}{ADCA} & \multicolumn{2}{c|}{InDCA} & \multicolumn{2}{c|}{HDCA-LI} & \multicolumn{2}{c}{HDCA-NI}\\
			& & $f$ & $c$ & $f$ & $c$ & $f$ & $c$ & $f$ & $c$ & $f$ & $c$ & $f$ & $c$ & $f$ & $c$ \\
			\midrule
			bfw398a & 7.58e+03 & 1.30e-04 & 1.32 & 7.86e-05 & 1.40 & 1.13e-04 & 1.34 & 1.42e-06 & 2.35 & 1.88e-06 & 1.93 & 8.87e-07 & 2.20 & 9.42e-07 & 2.21\\
			bfw62a & 1.48e+03 & 4.30e-06 & 1.55 & 1.57e-06 & 1.80 & 3.26e-06 & 1.61 & 9.88e-07 & 2.20 & 8.85e-07 & 1.94 & 6.42e-07 & 2.10 & 4.06e-07 & 2.23\\
			bfw782a & 4.62e+03 & 1.98e-04 & 0.94 & 1.27e-04 & 1.04 & 1.83e-04 & 0.96 & 3.45e-06 & 2.00 & 2.77e-06 & 1.77 & 1.30e-06 & 1.89 & 3.59e-07 & 2.04\\
			bwm200 & 2.93e+03 & 2.21e-08 & -2.47 & 2.09e-08 & -2.47 & 3.08e-08 & -2.47 & 3.07e-08 & -2.47 & 2.35e-08 & -2.47 & 3.03e-08 & -2.47 & 2.01e-08 & -2.36\\
			dwa512 & 3.72e+04 & 4.65e-04 & 0.57 & 2.12e-05 & 1.02 & 3.79e-04 & 0.60 & 1.94e-06 & 1.44 & 2.43e-06 & 1.42 & 1.59e-06 & 1.49 & 8.21e-08 & 2.44\\
			dwb512 & 4.50e+00 & 3.56e-04 & 1.65 & 2.20e-05 & 2.15 & 3.53e-04 & 1.65 & 7.63e-05 & 1.91 & 1.64e-05 & 2.35 & 2.34e-06 & 2.67 & 1.17e-07 & 3.01\\
			lop163 & 3.42e+07 & 3.06e-04 & 0.94 & 4.59e-06 & 1.92 & 2.78e-04 & 0.96 & 2.42e-06 & 2.04 & 3.47e-05 & 1.58 & 1.62e-05 & 1.68 & 2.98e-07 & 2.32\\
			mhd416a & 2.41e+25 & 1.80e-09 & 0.18 & 4.55e-09 & 0.18 & 3.62e-09 & 0.18 & 4.56e-09 & 0.18 & 8.75e-09 & 0.24 & 3.34e-10 & 0.27 & 1.67e-08 & 1.62\\
			mhd416b & 5.05e+09 & 9.55e-08 & 3.57 & 9.04e-08 & 3.57 & 9.72e-08 & 3.57 & 1.20e-07 & 3.57 & 5.67e-05 & 2.57 & 2.71e-08 & 5.04 & 3.29e-09 & 5.05\\
			odep400a & 8.31e+05 & 3.50e-04 & 0.10 & 2.53e-04 & 0.15 & 3.35e-04 & 0.10 & 1.17e-05 & 0.74 & 9.26e-05 & 0.41 & 4.41e-05 & 0.54 & 2.98e-08 & 1.97\\
			olm100 & 2.78e+04 & 2.00e-08 & -1.20 & 1.65e-08 & -1.20 & 2.82e-08 & -1.20 & 2.38e-08 & -1.20 & 4.31e-08 & -1.49 & 4.42e-08 & -1.49 & 3.41e-08 & -1.55\\
			olm500 & 7.65e+05 & 1.85e-08 & -2.80 & 3.14e-07 & -2.80 & 2.02e-08 & -2.80 & 1.41e-08 & -2.81 & 2.61e-08 & -2.86 & 4.38e-06 & -2.91 & 1.12e-08 & -2.81\\
			rbs480a & 1.35e+05 & 2.74e-08 & -1.17 & 2.83e-08 & -1.17 & 2.84e-08 & -1.17 & 3.10e-08 & -1.17 & 6.81e-08 & -2.09 & 5.38e-08 & -2.09 & 1.48e-08 & -1.85\\
			rbs480b & 1.63e+05 & 4.01e-08 & -1.57 & 3.96e-08 & -1.57 & 3.99e-08 & -1.57 & 3.81e-08 & -1.57 & -4.37e-08 & -2.08 & 3.81e-08 & -2.08 & 9.70e-08 & -1.98\\
			rdb200 & 8.32e+02 & 1.24e-05 & -0.78 & 1.23e-05 & -0.78 & 1.24e-05 & -0.78 & 1.22e-05 & -0.78 & 3.48e-07 & -0.11 & 2.46e-07 & -0.11 & 2.46e-07 & -0.15\\
			rdb450 & 1.64e+03 & 2.49e-06 & -1.01 & 2.47e-06 & -1.01 & 2.52e-06 & -1.01 & 2.49e-06 & -1.01 & 8.67e-08 & -0.49 & 8.44e-08 & -0.49 & 1.21e-07 & -0.55\\
			rdb968 & 2.91e+01 & 9.80e-06 & -0.41 & 8.82e-06 & -0.41 & 9.48e-06 & -0.41 & 7.74e-06 & -0.40 & 3.13e-06 & -0.23 & 2.68e-06 & -0.21 & 3.41e-07 & 0.05\\
			rw136 & 1.49e+05 & 2.96e-05 & 1.48 & 2.95e-05 & 1.48 & 2.95e-05 & 1.48 & 2.77e-05 & 1.48 & 2.26e-05 & 1.68 & 5.52e-06 & 1.89 & 1.29e-06 & 2.01\\
			rw496 & 1.14e+10 & 6.59e-06 & 1.78 & 6.62e-06 & 1.78 & 7.31e-06 & 1.78 & 6.59e-06 & 1.78 & 1.99e-05 & 1.86 & 3.42e-06 & 2.13 & 1.67e-08 & 2.85\\
			tols340 & 2.35e+05 & 2.35e-08 & 4.44 & 1.08e-07 & 0.50 & 3.22e-08 & 1.92 & 5.78e-08 & -3.26 & 1.02e-08 & 0.97 & 1.73e-08 & 0.81 & 1.57e-08 & 2.58\\
			tols90 & 2.49e+04 & 1.87e-08 & -2.05 & 1.83e-08 & -2.05 & 1.98e-08 & -2.05 & 2.71e-08 & -2.05 & 3.00e-08 & -2.29 & 3.14e-08 & -2.28 & 2.25e-08 & -2.31\\
			tub100 & 2.36e+04 & 1.54e-08 & -2.67 & 2.17e-08 & -2.67 & 1.24e-08 & -2.67 & 1.41e-08 & -2.67 & 2.51e-08 & -2.67 & 9.46e-09 & -2.67 & 1.90e-08 & -2.56\\
			\midrule
			\multicolumn{2}{c|}{avg} & 8.49e-05 & 0.11 & 2.58e-05 & 0.04 & 7.76e-05 & 0.00 & 7.06e-06 & 0.01 & 1.16e-05 & 0.09 & 3.80e-06 & 0.27 & 2.05e-07 & 0.65 \\
			\bottomrule
	\end{tabular}}	
\end{table}

\begin{table}[tbhp]
	\caption{Numerical results of DCA, BDCAe, BDCAa, ADCA, InDCA, HDCA-LI and HDCA-NI for solving \eqref{prob:dcp2} on the NEP dataset with \texttt{MaxIT}=200.}
	\label{tab:performcomp_allDCAs_dcp2}
	\centering
	\resizebox{\columnwidth}{!}{
		\begin{tabular}{l|c|cc|cc|cc|cc|cc|cc|cc} 
			\toprule
			\multirow{2}{*}{Prob} & \multirow{2}{*}{cond($A$)} & \multicolumn{2}{c|}{DCA} & \multicolumn{2}{c|}{BDCAe} & \multicolumn{2}{c|}{BDCAa} & \multicolumn{2}{c|}{ADCA} & \multicolumn{2}{c|}{InDCA} & \multicolumn{2}{c|}{HDCA-LI} & \multicolumn{2}{c}{HDCA-NI}\\
			& & $f$ & $c$ & $f$ & $c$ & $f$ & $c$ & $f$ & $c$ & $f$ & $c$ & $f$ & $c$ & $f$ & $c$ \\
			\midrule
			bfw398a & 7.58e+03 & 1.73e-04 & 1.28 & 1.32e-04 & 1.32 & 1.52e-04 & 1.30 & 7.13e-06 & 1.96 & 1.53e-06 & 1.99 & 1.00e-06 & 2.12 & 4.40e-07 & 2.20\\
			bfw62a & 1.48e+03 & 1.59e-06 & 1.82 & 2.53e-07 & 2.31 & 9.28e-07 & 1.95 & 6.10e-07 & 2.11 & 3.69e-07 & 2.24 & 2.32e-07 & 2.49 & 1.49e-07 & 2.63\\
			bfw782a & 4.62e+03 & 2.33e-04 & 0.92 & 1.84e-04 & 0.98 & 2.19e-04 & 0.93 & 1.92e-05 & 1.50 & 2.96e-06 & 1.83 & 1.36e-06 & 1.96 & 5.71e-07 & 2.11\\
			bwm200 & 2.93e+03 & 3.44e-06 & -2.46 & 3.96e-06 & -2.45 & 3.47e-06 & -2.46 & 4.09e-06 & -2.46 & 3.34e-06 & -2.46 & 4.34e-06 & -2.46 & 2.87e-06 & -2.45\\
			dwa512 & 3.72e+04 & 5.87e-04 & 0.55 & 2.37e-05 & 1.02 & 3.74e-04 & 0.60 & 3.01e-06 & 1.35 & 5.54e-05 & 0.98 & 2.08e-05 & 1.08 & 5.46e-08 & 2.14\\
			dwb512 & 4.50e+00 & 3.60e-04 & 1.65 & 2.25e-05 & 2.13 & 3.49e-04 & 1.65 & 9.88e-05 & 1.86 & 7.15e-05 & 1.99 & 2.49e-05 & 2.16 & 7.21e-06 & 2.23\\
			lop163 & 3.42e+07 & 3.24e-04 & 0.93 & 1.68e-05 & 1.54 & 2.56e-04 & 0.98 & 3.75e-06 & 1.91 & 7.21e-05 & 1.35 & 2.79e-05 & 1.51 & 1.66e-06 & 2.09\\
			mhd416a & 2.41e+25 & 1.01e-06 & 1.00 & 7.37e-07 & 1.07 & 1.02e-06 & 1.00 & 1.10e-06 & 0.98 & 1.03e-06 & 0.99 & 1.04e-06 & 0.99 & 3.39e-06 & 0.74\\
			mhd416b & 5.05e+09 & 1.24e-07 & 3.57 & 1.35e-07 & 3.57 & 1.29e-07 & 3.57 & 1.64e-07 & 3.57 & 2.37e-04 & 1.80 & 1.33e-07 & 4.32 & 6.28e-09 & 5.42\\
			odep400a & 8.31e+05 & 3.68e-04 & 0.09 & 4.16e-05 & 0.47 & 3.25e-04 & 0.11 & 2.23e-05 & 0.57 & 1.76e-04 & 0.28 & 4.63e-05 & 0.50 & 1.64e-06 & 1.12\\
			olm100 & 2.78e+04 & 9.22e-07 & -1.18 & 1.95e-06 & -1.17 & 9.27e-07 & -1.18 & 1.73e-06 & -1.17 & 2.89e-06 & -1.17 & 2.89e-06 & -1.18 & 3.19e-06 & -1.17\\
			olm500 & 7.65e+05 & 6.61e-06 & -2.85 & 6.61e-06 & -2.85 & 6.61e-06 & -2.85 & 6.87e-06 & -2.85 & 6.45e-06 & -2.85 & 6.45e-06 & -2.85 & 1.02e-05 & -2.89\\
			rbs480a & 1.35e+05 & -6.50e-07 & -1.17 & -7.77e-07 & -1.16 & -1.14e-06 & -1.16 & -9.94e-07 & -1.16 & -7.38e-08 & -1.18 & -7.38e-08 & -1.18 & -5.44e-07 & -1.20\\
			rbs480b & 1.63e+05 & -1.69e-06 & -1.56 & -3.85e-07 & -1.56 & -6.90e-07 & -1.55 & -8.60e-07 & -1.55 & -5.78e-07 & -1.50 & -5.78e-07 & -1.50 & -9.76e-07 & -1.50\\
			rdb200 & 8.32e+02 & 1.25e-05 & -0.78 & 1.23e-05 & -0.78 & 1.25e-05 & -0.78 & 1.25e-05 & -0.78 & 1.21e-05 & -0.78 & 1.21e-05 & -0.78 & 1.21e-05 & -0.77\\
			rdb450 & 1.64e+03 & 2.60e-06 & -1.01 & 2.59e-06 & -1.01 & 2.61e-06 & -1.01 & 2.60e-06 & -1.01 & 2.91e-06 & -1.01 & 2.75e-06 & -1.00 & 3.09e-06 & -1.00\\
			rdb968 & 2.91e+01 & 9.01e-06 & -0.41 & 8.61e-06 & -0.41 & 8.85e-06 & -0.41 & 8.40e-06 & -0.40 & 6.30e-06 & -0.25 & 4.56e-06 & -0.25 & 4.92e-06 & -0.29\\
			rw136 & 1.49e+05 & 2.97e-05 & 1.48 & 2.95e-05 & 1.48 & 2.95e-05 & 1.48 & 2.91e-05 & 1.48 & 5.22e-05 & 1.56 & 2.12e-05 & 1.62 & 2.49e-06 & 1.91\\
			rw496 & 1.14e+10 & 6.75e-06 & 1.78 & 6.68e-06 & 1.78 & 6.67e-06 & 1.78 & 6.69e-06 & 1.78 & 4.18e-05 & 1.78 & 1.77e-05 & 1.85 & 7.17e-08 & 2.77\\
			tols340 & 2.35e+05 & 5.10e-05 & -0.57 & -5.80e-06 & -0.61 & -8.79e-07 & -0.85 & 1.30e-04 & -1.10 & 5.83e-05 & -0.67 & -7.83e-06 & -0.63 & 8.72e-05 & -1.20\\
			tols90 & 2.49e+04 & 2.07e-05 & -2.28 & 2.34e-05 & -2.28 & 1.08e-05 & -2.27 & 2.46e-05 & -2.24 & 1.88e-05 & -2.29 & 2.16e-05 & -2.25 & 2.32e-05 & -2.18\\
			tub100 & 2.36e+04 & 5.75e-06 & -2.66 & 5.73e-06 & -2.66 & 5.75e-06 & -2.66 & 5.75e-06 & -2.66 & 5.73e-06 & -2.66 & 5.73e-06 & -2.66 & 5.75e-06 & -2.66\\
			\midrule
			\multicolumn{2}{c|}{avg} & 9.98e-05 & -0.08 & 2.34e-05 & 0.03 & 8.01e-05 & -0.08 & 1.76e-05 & 0.08 & 3.76e-05 & -0.00 & 9.75e-06 & 0.18 & 7.67e-06 & 0.37 \\
			\bottomrule
	\end{tabular}}	
\end{table}

\begin{table}[tbhp]
	\caption{Numerical results of DCA, BDCAa, ADCA, InDCA, HDCA-LI and HDCA-NI for solving \eqref{prob:dcp3} on the NEP dataset with \texttt{MaxIT}=200.}
	\label{tab:performcomp_allDCAs_dcp3}
	\centering
	\resizebox{\columnwidth}{!}{
		\begin{tabular}{l|c|cc|cc|cc|cc|cc|cc} 
			\toprule
			\multirow{2}{*}{Prob} & \multirow{2}{*}{cond($A$)} & \multicolumn{2}{c|}{DCA} &  \multicolumn{2}{c|}{BDCAa} & \multicolumn{2}{c|}{ADCA} & \multicolumn{2}{c|}{InDCA} & \multicolumn{2}{c|}{HDCA-LI} & \multicolumn{2}{c}{HDCA-NI}\\
			& & $f$ & $c$ & $f$ & $c$ & $f$ & $c$ & $f$ & $c$ & $f$ & $c$ & $f$ & $c$ \\
			\midrule
			bfw398a & 7.58e+03 & 1.18e-02 & 0.85 & 9.35e-03 & 0.85 & 2.03e-03 & 0.86 & 1.18e-02 & 0.85 & 9.24e-03 & 0.85 & 2.02e-03 & 0.86\\
			bfw62a & 1.48e+03 & 5.22e-03 & 0.78 & 1.54e-03 & 0.79 & 7.53e-04 & 0.88 & 4.98e-03 & 0.78 & 1.54e-03 & 0.79 & 7.65e-04 & 0.87\\
			bfw782a & 4.62e+03 & 1.36e-03 & 0.88 & 1.21e-03 & 0.88 & 6.29e-04 & 0.88 & 1.34e-03 & 0.88 & 1.24e-03 & 0.88 & 6.22e-04 & 0.88\\
			bwm200 & 2.93e+03 & 1.39e-08 & -2.47 & 1.38e-08 & -2.47 & 1.69e-08 & -2.47 & 1.40e-08 & -2.47 & 1.38e-08 & -2.47 & 1.25e-08 & -2.31\\
			dwa512 & 3.72e+04 & 2.16e-05 & 0.97 & 2.15e-05 & 0.97 & 2.15e-05 & 0.97 & 2.13e-05 & 0.98 & 2.15e-05 & 0.97 & 2.13e-05 & 0.98\\
			dwb512 & 4.50e+00 & 3.74e-05 & 1.99 & 3.73e-05 & 1.99 & 3.75e-05 & 1.99 & 3.74e-05 & 1.99 & 3.78e-05 & 1.99 & 3.77e-05 & 1.99\\
			lop163 & 3.42e+07 & 1.75e-04 & 0.94 & 1.65e-04 & 0.95 & 8.51e-05 & 1.03 & 1.72e-04 & 0.95 & 1.61e-04 & 0.96 & 8.27e-05 & 1.04\\
			mhd416a & 2.41e+25 & 4.30e-09 & 0.70 & 3.58e-09 & 0.70 & 1.36e-08 & 0.70 & 2.41e-08 & 0.23 & 3.18e-09 & 0.70 & 2.70e-11 & -0.09\\
			mhd416b & 5.05e+09 & 9.95e-06 & 2.54 & 1.01e-05 & 2.54 & 1.01e-05 & 2.54 & 1.07e-05 & 2.54 & 1.04e-05 & 2.54 & 1.09e-05 & 2.54\\
			odep400a & 8.31e+05 & 9.86e-05 & 0.30 & 9.84e-05 & 0.30 & 5.28e-05 & 0.44 & 9.81e-05 & 0.30 & 9.81e-05 & 0.30 & 4.23e-05 & 0.48\\
			olm100 & 2.78e+04 & 4.58e-09 & -1.20 & 1.34e-08 & -1.20 & 4.43e-09 & -1.20 & 3.28e-09 & -1.44 & 1.51e-08 & -1.16 & 3.38e-09 & -1.51\\
			olm500 & 7.65e+05 & 1.44e-08 & -2.81 & 1.82e-08 & -2.81 & 1.28e-08 & -2.81 & 9.81e-09 & -2.82 & 2.04e-08 & -2.66 & 1.52e-08 & -2.78\\
			rbs480a & 1.35e+05 & 1.70e-08 & -1.95 & 2.31e-08 & -1.95 & 1.46e-08 & -1.95 & 2.04e-08 & -2.09 & 2.42e-08 & -1.89 & 1.43e-08 & -2.00\\
			rbs480b & 1.63e+05 & 3.55e-08 & -2.01 & 2.05e-08 & -2.01 & 2.41e-08 & -2.01 & 3.11e-08 & -2.08 & 1.99e-08 & -1.95 & 1.82e-08 & -2.01\\
			rdb200 & 8.32e+02 & 4.97e-07 & -0.49 & 5.00e-07 & -0.49 & 5.12e-07 & -0.49 & 7.61e-07 & -0.33 & 5.77e-07 & -0.49 & 4.79e-07 & -0.18\\
			rdb450 & 1.64e+03 & 2.11e-07 & -0.86 & 2.16e-07 & -0.86 & 2.22e-07 & -0.86 & 1.95e-07 & -0.74 & 2.95e-07 & -0.85 & 9.34e-08 & -0.59\\
			rdb968 & 2.91e+01 & 5.61e-06 & -0.47 & 5.62e-06 & -0.47 & 4.80e-06 & -0.44 & 5.37e-06 & -0.46 & 5.61e-06 & -0.47 & 3.48e-06 & -0.37\\
			rw136 & 1.49e+05 & 2.79e-04 & 0.84 & 2.73e-04 & 0.84 & 2.03e-04 & 0.90 & 2.72e-04 & 0.84 & 2.67e-04 & 0.84 & 1.97e-04 & 0.91\\
			rw496 & 1.14e+10 & 1.85e-04 & 0.89 & 1.85e-04 & 0.89 & 1.67e-04 & 0.91 & 1.83e-04 & 0.89 & 1.81e-04 & 0.90 & 1.65e-04 & 0.92\\
			tols340 & 2.35e+05 & 4.18e-08 & 4.79 & 4.81e-08 & 4.76 & 5.42e-09 & 5.26 & 2.00e-09 & -2.90 & 1.46e-08 & -0.48 & 2.13e-09 & -2.85\\
			tols90 & 2.49e+04 & 1.15e-08 & -2.05 & 9.65e-09 & -2.05 & 2.25e-08 & -2.05 & 5.14e-09 & -2.27 & 4.66e-09 & -2.05 & 2.06e-08 & -2.31\\
			tub100 & 2.36e+04 & 5.71e-09 & -2.67 & 4.84e-09 & -2.67 & 5.57e-09 & -2.67 & 9.87e-09 & -2.67 & 5.16e-09 & -2.67 & 3.54e-09 & -2.50\\
			\midrule
			\multicolumn{2}{c|}{avg} & 8.73e-04 & -0.02 & 5.86e-04 & -0.02 & 1.81e-04 & 0.02 & 8.59e-04 & -0.41 & 5.82e-04 & -0.25 & 1.81e-04 & -0.36\\
			\bottomrule
	\end{tabular}}	
\end{table}

\begin{table}[tbhp]
	\caption{Numerical results of InDCA, HDCA-LI and HDCA-NI using the conservative inertial strategy for solving \eqref{prob:dcp3} on the NEP dataset with \texttt{MaxIT}=200.}
	\label{tab:performcomp_allDCAs_dcp3_conserv}
	\centering
	\resizebox{0.6\columnwidth}{!}{
		\begin{tabular}{l|c|cc|cc|cc} 
			\toprule
			\multirow{2}{*}{Prob} & \multirow{2}{*}{cond($A$)} & \multicolumn{2}{c|}{InDCA} & \multicolumn{2}{c|}{HDCA-LI} & \multicolumn{2}{c}{HDCA-NI}\\
			& & $f$ & $c$ & $f$ & $c$ & $f$ & $c$  \\
			\midrule
			bfw398a & 7.58e+03 & 1.18e-02 & 0.85 & 9.29e-03 & 0.85 & 2.03e-03 & 0.86\\
			bfw62a & 1.48e+03 & 5.21e-03 & 0.78 & 1.72e-03 & 0.79 & 7.54e-04 & 0.88\\
			bfw782a & 4.62e+03 & 1.36e-03 & 0.87 & 1.20e-03 & 0.88 & 6.27e-04 & 0.88\\
			bwm200 & 2.93e+03 & 1.37e-08 & -2.47 & 1.45e-08 & -2.47 & 3.25e-08 & -2.47\\
			dwa512 & 3.72e+04 & 2.16e-05 & 0.97 & 2.15e-05 & 0.97 & 2.15e-05 & 0.97\\
			dwb512 & 4.50e+00 & 3.74e-05 & 1.99 & 3.76e-05 & 1.99 & 3.76e-05 & 1.99\\
			lop163 & 3.42e+07 & 1.75e-04 & 0.94 & 1.63e-04 & 0.95 & 8.48e-05 & 1.03\\
			mhd416a & 2.41e+25 & 5.16e-09 & 0.70 & 5.49e-09 & 0.70 & 1.41e-08 & 0.70\\
			mhd416b & 5.05e+09 & 1.06e-05 & 2.54 & 9.73e-06 & 2.54 & 1.04e-05 & 2.54\\
			odep400a & 8.31e+05 & 9.87e-05 & 0.30 & 9.67e-05 & 0.30 & 5.29e-05 & 0.44\\
			olm100 & 2.78e+04 & 1.74e-08 & -1.20 & 1.75e-08 & -1.20 & 1.31e-08 & -1.20\\
			olm500 & 7.65e+05 & 2.24e-08 & -2.81 & 1.78e-08 & -2.81 & 4.75e-07 & -2.81\\
			rbs480a & 1.35e+05 & 2.78e-08 & -1.95 & 1.48e-08 & -1.95 & 1.85e-08 & -1.95\\
			rbs480b & 1.63e+05 & 2.38e-08 & -2.01 & 1.92e-08 & -2.01 & 3.15e-08 & -2.01\\
			rdb200 & 8.32e+02 & 5.04e-07 & -0.49 & 5.00e-07 & -0.49 & 5.13e-07 & -0.49\\
			rdb450 & 1.64e+03 & 2.15e-07 & -0.86 & 2.13e-07 & -0.86 & 2.17e-07 & -0.86\\
			rdb968 & 2.91e+01 & 5.61e-06 & -0.47 & 5.61e-06 & -0.47 & 4.81e-06 & -0.43\\
			rw136 & 1.49e+05 & 2.79e-04 & 0.84 & 2.72e-04 & 0.84 & 2.02e-04 & 0.90\\
			rw496 & 1.14e+10 & 1.85e-04 & 0.89 & 1.85e-04 & 0.89 & 1.69e-04 & 0.91\\
			tols340 & 2.35e+05 & 1.50e-07 & 4.57 & 8.07e-09 & 5.06 & 5.21e-08 & 4.76\\
			tols90 & 2.49e+04 & 1.76e-08 & -2.05 & 7.84e-09 & -2.05 & 2.14e-08 & -2.05\\
			tub100 & 2.36e+04 & 5.37e-09 & -2.67 & 5.73e-09 & -2.67 & 5.49e-09 & -2.67\\
			\midrule
			\multicolumn{2}{c|}{avg} & 8.73e-04 & -0.03 & 5.91e-04 & -0.01 & 1.82e-04 & -0.00\\
			\bottomrule
	\end{tabular}}	
\end{table}

As seen in \Cref{tab:performcomp_allDCAs_dcp1} for \eqref{prob:dcp1}, the best average value in $f$ is achieved by HDCA-NI, followed by HDCA-LI, InDCA, BDCAe, BDCAa, and DCA. The best average value in $c$ is also obtained by HDCA-NI, with the subsequent ranking being HDCA-LI, DCA, InDCA, BDCAe, ADCA, and BDCAa. Similar results are observed in \Cref{tab:performcomp_allDCAs_dcp2} for \eqref{prob:dcp2}.  It is worth noting that the inertial-based methods (InDCA, HDCA-LI, and HDCA-NI) perform poorly for \eqref{prob:dcp3} when $A$ is ill-conditioned, as seen in \Cref{tab:performcomp_allDCAs_dcp3}. For example, in the case \texttt{tols340} where cond($A$)=2.35e+05, the value of $c$ is $-2.85$ for HDCA-NI, $-2.90$ for InDCA, $-0.48$ for HDCA-BI, $4.79$ for DCA, $4.76$ for BDCAa, and $5.26$ for ADCA. To some extent, this can be alleviated by using a `conservative' inertial strategy, which introduces the inertial force $\gamma (x^k-x^{k-1})$ when $k\geq 2$ instead of $k\geq 0$. Then, we will get improvement in $c$ as $4.57$ for InDCA, $5.00$ for HDCA-LI, and $4.76$ for HDCA-NI. The numerical results of InDCA, HDCA-LI, and HDCA-NI using the conservative inertial strategy are summarized in \Cref{tab:performcomp_allDCAs_dcp3_conserv}, and we can observe that the average value of $c$ is greatly improved at the cost of a slight increase in the average value of $f$.

We conclude that the best DCA-type algorithm is always given by HDCA-NI, hence in the next subsection, we will compare it with other optimization solvers. 

\subsection{Numerical results of other optimization solvers}\label{subsec:resultsforothers}
In this section, we present the numerical results of the optimization solvers IPOPT, KNITRO, and FILTERSD in comparison with HDCA-NI. Our primary focus is to evaluate their performance based on the objective value $f$, the feasibility measure $c$, and the CPU time (in seconds). To determine the CPU time for HDCA-NI, we employ the following stopping criteria:
$$|f^{k+1} - f^{k}|\leq (1+|f^{k+1}|)\varepsilon,$$
where the tolerance $\varepsilon$ is set to $10^{-8}$, and $f^k$ represents the objective value at the $k$-th iteration. The other solvers use their default termination settings. The numerical results for solving \eqref{prob:nlp1},\eqref{prob:nlp2}, and \eqref{prob:nlp3} are summarized in \Cref{tab:perform_others_nlp1,tab:perform_others_nlp2,tab:perform_others_nlp3}, respectively. Furthermore, instead of providing the results for all $10$ instances within each dataset \RAND{n}, we only present their averages.

\begin{itemize}[leftmargin=12pt]
	\item For the results of \eqref{prob:nlp1} in \Cref{tab:perform_others_nlp1}, we observe that IPOPT performs best in terms of average values for $f$, $c$, and CPU time. HDCA-NI is the second-best method, followed by KNITRO and FILTERSD. It is important to note that, in contrast to IPOPT, FILTERSD and HDCA-NI, the solver KNITRO demonstrates considerable instability and frequently encounters ``\texttt{LCP solver problem}" issues, in ill-conditioned NEP instances. This results in significantly different outcomes in each run. Therefore, we consider the best result for FILTERSD in terms of $f$ across three runs. Nevertheless, FILTERSD still exhibits poor performance in terms of average $f$, $c$, and CPU time.
	\item For the results of \eqref{prob:nlp2} in \Cref{tab:perform_others_nlp2} and \eqref{prob:nlp3} in \Cref{tab:perform_others_nlp3}, IPOPT remains the best solver in terms of average $f$, $c$, and CPU time. HDCA-NI often yields good average values for $f$, but with poor average values for $c$. FILTERSD consistently provides the worst results for $f$.
\end{itemize}

\begin{table}[tbhp]
	\caption{Numerical results of IPOPT, KNITRO, FILTERSD and HDCA-NI  for solving \eqref{prob:nlp1} on \RAND{n} and NEP datasets.}
	\label{tab:perform_others_nlp1}
	\centering
	\resizebox{0.9\columnwidth}{!}{
		\begin{tabular}{l|ccc|ccc|ccc|ccc} \toprule
			\multirow{2}{*}{Dataset} & \multicolumn{3}{c|}{HDCA-NI} & \multicolumn{3}{c|}{IPOPT} & \multicolumn{3}{c|}{KNITRO} & \multicolumn{3}{c}{FILTERSD} \\
			& $f$ & $c$ & CPU & $f$ & $c$ & CPU & $f$ & $c$ & CPU & $f$ & $c$ & CPU \\
			\midrule
			RAND(10) & 1.40e-02 & 1.88 & 0.451 & 4.71e-03 & 6.19 & 0.022 & 9.14e-03 & 4.17 & 0.008 & 3.79e-03 & 3.17 & 0.006\\
			RAND(100) & 1.41e-03 & 1.03 & 2.793 & 4.71e-04 & 3.45 & 0.816 & 9.15e-04 & 2.32 & 1.646 & 2.81e+01 & 0.15 & 0.670\\
			RAND(500) & 1.41e-04 & 0.55 & 28.157 & 4.79e-05 & 1.50 & 5.854 & 6.61e-04 & 0.37 & 17.710 & 5.69e+03 & -0.39 & 32.241\\
			\midrule
			bfw398a & 1.39e-06 & 2.14 & 8.173 & 1.68e-05 & 2.88 & 0.428 & 1.82e-03 & 1.39 & 0.856 & 6.08e+01 & 1.27 & 1.839\\
			bfw62a & 6.42e-07 & 2.03 & 0.533 & 3.82e-10 & 4.54 & 0.315 & 8.43e-07 & 2.52 & 0.602 & 1.88e-07 & 2.09 & 0.659\\
			bfw782a & 1.00e-06 & 1.86 & 11.712 & 3.36e-05 & 2.65 & 0.492 & 7.20e-04 & 1.16 & 6.106 & 1.25e+02 & 1.38 & 6.422\\
			bwm200 & 1.07e-08 & -2.45 & 0.408 & 3.98e-08 & -0.79 & 0.029 & 5.50e-06 & -2.01 & 0.855 & 3.13e+01 & -1.23 & 0.091\\
			dwa512 & 8.78e-07 & 1.74 & 4.129 & 1.69e-08 & 2.73 & 0.840 & 1.00e-04 & 1.93 & 2.908 & 8.21e-04 & 0.50 & 68.637\\
			dwb512 & 2.50e-07 & 2.91 & 4.415 & 3.83e-08 & 3.60 & 0.647 & 1.85e-05 & 2.18 & 9.078 & 3.75e-04 & 1.65 & 4.769\\
			lop163 & 6.59e-07 & 2.15 & 6.083 & 6.42e-13 & 4.99 & 0.790 & 8.13e-07 & 2.35 & 0.256 & 4.03e-04 & 0.89 & 2.824\\
			mhd416a & -1.29e-09 & 1.67 & 0.256 & 1.10e-07 & -1.11 & 0.092 & 1.75e-04 & -0.46 & 6.693 & 6.57e+01 & -0.37 & 0.787\\
			mhd416b & 2.06e-07 & 3.12 & 3.537 & 4.52e-07 & 3.42 & 0.681 & 1.51e-03 & 1.89 & 0.936 & 6.85e+01 & 1.79 & 0.522\\
			odep400a & 2.91e-07 & 1.54 & 4.077 & 3.71e-08 & 3.14 & 0.290 & 8.27e-05 & 1.78 & 1.087 & 4.16e-04 & 0.07 & 25.959\\
			olm100 & 2.82e-09 & -1.55 & 0.115 & 3.48e-09 & -0.88 & 0.022 & 1.63e-05 & -0.88 & 0.111 & 1.24e-09 & -1.41 & 0.053\\
			olm500 & -5.14e-08 & -2.91 & 3.567 & 1.62e-06 & -1.73 & 0.081 & 1.10e-03 & -2.11 & 14.042 & 1.03e-10 & -3.68 & 18.787\\
			rbs480a & 2.52e-08 & -1.92 & 2.548 & 7.03e-08 & -0.73 & 0.455 & 1.17e-03 & -0.72 & 6.043 & 8.18e+01 & -0.85 & 6.370\\
			rbs480b & 2.40e-08 & -2.03 & 1.933 & 8.32e-08 & -1.17 & 0.719 & 1.18e-03 & -1.19 & 6.771 & 4.79e+04 & -2.01 & 25.796\\
			rdb200 & 2.81e-07 & -0.15 & 0.502 & 1.57e-08 & 0.30 & 0.125 & 1.61e-05 & -0.75 & 0.413 & 5.78e+03 & -0.79 & 0.187\\
			rdb450 & 1.33e-07 & -0.53 & 0.681 & 1.54e-08 & 0.50 & 0.065 & 5.00e-04 & 0.38 & 1.977 & 1.84e+04 & -0.99 & 1.954\\
			rdb968 & 4.46e-06 & -0.28 & 2.027 & 8.28e-08 & 1.58 & 0.324 & 3.35e-04 & 1.19 & 18.358 & 1.43e+02 & 0.70 & 11.839\\
			rw136 & 1.43e-06 & 2.00 & 4.554 & 1.54e-13 & 5.32 & 0.865 & 1.78e-07 & 2.57 & 4.571 & 2.97e-05 & 1.48 & 1.797\\
			rw496 & 3.06e-07 & 2.27 & 5.993 & 1.44e-08 & 4.85 & 0.360 & 3.00e-05 & 2.22 & 3.469 & 6.60e-06 & 1.77 & 4.325\\
			tols340 & 2.87e-06 & -0.03 & 1.918 & 6.85e-06 & -2.29 & 0.064 & 1.39e-03 & -2.80 & 4.342 & 2.09e-01 & -2.09 & 1.808\\
			tols90 & 3.51e-09 & -2.32 & 0.086 & 1.68e-07 & -2.04 & 0.036 & 4.36e-05 & -2.19 & 0.257 & 1.25e-15 & -0.34 & 0.161\\
			tub100 & 4.56e-09 & -2.64 & 0.267 & 2.85e-08 & -1.75 & 0.018 & 4.36e-07 & -1.73 & 0.149 & 4.02e-09 & -2.62 & 0.057\\
			\midrule
			avg & 6.24e-06 & 0.29 & 3.827 & 4.32e-06 & 1.18 & 0.544 & 4.35e-04 & 0.28 & 4.304 & 3.13e+03 & -0.13 & 8.715\\
			\bottomrule
	\end{tabular}}	
\end{table}

\begin{table}[tbhp]
	\caption{Numerical results of IPOPT, KNITRO, FILTERSD and HDCA-NI for solving \eqref{prob:nlp2} on \RAND{n} and NEP datasets.}
	\label{tab:perform_others_nlp2}
	\centering
	\resizebox{0.9\columnwidth}{!}{
		\begin{tabular}{l|ccc|ccc|ccc|ccc} \toprule
			\multirow{2}{*}{Dataset} & \multicolumn{3}{c|}{HDCA-NI} & \multicolumn{3}{c|}{IPOPT} & \multicolumn{3}{c|}{KNITRO} & \multicolumn{3}{c}{FILTERSD} \\
			& $f$ & $c$ & CPU & $f$ & $c$ & CPU & $f$ & $c$ & CPU & $f$ & $c$ & CPU \\
			\midrule
			RAND(10) & 1.41e-02 & 2.00 & 0.594 & 3.27e-03 & 6.13 & 0.025 & 5.66e-03 & 5.13 & 0.009 & 6.95e-03 & 2.98 & 0.023\\
			RAND(100) & 1.41e-03 & 0.99 & 5.435 & 3.27e-04 & 3.19 & 0.726 & 5.66e-04 & 2.36 & 0.813 & 1.06e+01 & 0.13 & 0.780\\
			RAND(500) & 1.45e-04 & 0.24 & 30.844 & 3.36e-05 & 1.48 & 8.736 & 1.44e-03 & 0.48 & 3.801 & 1.35e+02 & -0.21 & 6.294\\
			\midrule
			bfw398a & 9.86e-07 & 2.08 & 4.710 & 1.88e-08 & 3.63 & 2.405 & 8.69e-06 & 1.92 & 1.078 & 5.97e+01 & 1.27 & 0.850\\
			bfw62a & 9.22e-06 & 1.64 & 0.113 & 1.24e-08 & 3.52 & 0.230 & 1.88e-06 & 2.18 & 0.157 & 1.30e-06 & 1.95 & 0.764\\
			bfw782a & 1.53e-06 & 1.88 & 8.884 & 2.24e-05 & 2.73 & 1.318 & 5.85e-04 & 1.35 & 1.281 & 1.24e+02 & 1.38 & 6.211\\
			bwm200 & 6.33e-08 & -2.47 & 0.105 & 3.98e-08 & -0.78 & 0.033 & 1.19e-05 & -0.66 & 0.063 & -1.23e+03 & -1.23 & 0.082\\
			dwa512 & 3.72e-07 & 1.79 & 2.908 & 1.70e-08 & 2.72 & 0.736 & 7.13e-04 & 1.98 & 0.411 & 8.21e-04 & 0.50 & 3.533\\
			dwb512 & 1.70e-05 & 2.08 & 1.275 & 1.22e-07 & 3.45 & 1.051 & 1.60e-05 & 2.17 & 1.924 & 3.12e-08 & 3.59 & 33.208\\
			lop163 & 3.48e-06 & 1.90 & 3.343 & -7.70e-09 & 4.40 & 0.314 & 2.07e-03 & 1.93 & 0.030 & 4.03e-04 & 0.89 & 0.105\\
			mhd416a & 1.52e-08 & 2.99 & 0.142 & 1.10e-07 & -0.07 & 0.091 & 3.23e-04 & -0.55 & 0.202 & -1.20e+03 & -0.37 & 0.594\\
			mhd416b & 3.19e-07 & 2.63 & 2.249 & 9.96e-08 & 3.54 & 0.698 & 1.01e-05 & 2.35 & 1.132 & 6.79e+01 & 1.79 & 0.296\\
			odep400a & 1.54e-06 & 1.13 & 5.524 & 5.02e-09 & 3.30 & 0.605 & 8.33e-04 & 1.58 & 0.207 & 5.41e-08 & 1.81 & 23.659\\
			olm100 & 8.65e-07 & -1.17 & 0.206 & 2.48e-09 & -1.03 & 0.026 & 3.63e-07 & -0.82 & 0.017 & 1.98e-07 & -2.24 & 0.040\\
			olm500 & 8.16e-06 & -2.88 & 1.119 & 1.62e-06 & -1.38 & 0.077 & 5.00e-04 & -1.60 & 0.425 & -4.19e-12 & -2.92 & 7.570\\
			rbs480a & -8.66e-07 & -1.20 & 17.931 & 7.02e-08 & -0.61 & 0.313 & 7.25e-04 & -0.73 & 0.311 & -4.12e+02 & -0.85 & 3.955\\
			rbs480b & -4.15e-07 & -1.50 & 8.610 & 4.66e-06 & -1.18 & 0.845 & 1.29e-05 & -1.15 & 0.478 & -3.42e+02 & -1.15 & 8.460\\
			rdb200 & 1.23e-05 & -0.78 & 1.111 & 6.56e-08 & 0.35 & 0.085 & 9.03e-07 & 0.31 & 0.050 & 7.25e-06 & -0.67 & 0.323\\
			rdb450 & 3.01e-06 & -1.01 & 0.443 & 1.53e-08 & 0.49 & 0.078 & 9.54e-05 & 0.05 & 0.380 & 4.78e+04 & -1.23 & 2.466\\
			rdb968 & 5.36e-06 & -0.29 & 24.240 & 8.22e-08 & 1.88 & 0.435 & 2.96e-04 & 1.46 & 1.914 & 1.20e+02 & 0.70 & 14.062\\
			rw136 & 8.37e-06 & 1.66 & 2.067 & -8.18e-09 & 4.40 & 0.487 & 1.26e-07 & 4.43 & 1.682 & 5.27e-07 & 2.18 & 1.883\\
			rw496 & 3.57e-07 & 2.44 & 9.305 & 5.18e-06 & 3.68 & 0.322 & 2.80e-06 & 2.38 & 0.984 & 8.31e-08 & 2.56 & 47.600\\
			tols340 & 1.15e-04 & -2.31 & 17.904 & 6.86e-06 & -2.85 & 0.049 & 8.59e-04 & -2.55 & 0.217 & 2.20e-11 & 8.92 & 3.339\\
			tols90 & 5.37e-06 & -2.15 & 28.596 & 1.68e-07 & -2.09 & 0.025 & 2.87e-05 & -2.16 & 0.015 & 4.58e-09 & -2.03 & 0.047\\
			tub100 & 5.70e-06 & -2.66 & 0.590 & 2.78e-08 & -1.50 & 0.022 & 3.23e-06 & -1.07 & 0.015 & 1.11e-09 & -2.31 & 0.045\\
			\midrule
			avg & 1.37e-05 & 0.16 & 6.889 & 3.01e-06 & 1.12 & 0.759 & 3.41e-04 & 0.53 & 0.671 & 1.80e+03 & 0.49 & 6.615\\
			\bottomrule
	\end{tabular}}	
\end{table}

\begin{table}[tbhp]
	\caption{Numerical results of IPOPT, KNITRO, FILTERSD and HDCA-NI for solving \eqref{prob:nlp3} on \RAND{n} and NEP datasets.}
	\label{tab:perform_others_nlp3}
	\centering
	\resizebox{0.9\columnwidth}{!}{
		\begin{tabular}{l|ccc|ccc|ccc|ccc} \toprule
			\multirow{2}{*}{Dataset} & \multicolumn{3}{c|}{HDCA-NI} & \multicolumn{3}{c|}{IPOPT} & \multicolumn{3}{c|}{KNITRO} & \multicolumn{3}{c}{FILTERSD} \\
			& $f$ & $c$ & CPU & $f$ & $c$ & CPU & $f$ & $c$ & CPU & $f$ & $c$ & CPU \\
			\midrule
			RAND(10) & 6.86e-03 & 1.43 & 0.753 & 4.57e-03 & 6.23 & 0.018 & 4.72e-03 & 5.18 & 0.015 & 6.39e-03 & 2.41 & 0.008\\
			RAND(100) & 6.99e-04 & 0.62 & 3.267 & 4.57e-04 & 3.64 & 0.774 & 4.73e-04 & 3.11 & 1.606 & 3.02e+00 & 0.56 & 1.289\\
			RAND(500) & 7.37e-05 & 0.13 & 47.000 & 4.59e-05 & 1.60 & 6.789 & 4.73e-05 & 1.49 & 91.480 & 4.07e+01 & -0.63 & 36.676\\
			\midrule
			bfw398a & 4.66e-04 & 0.87 & 9.086 & 4.64e-07 & 3.30 & 1.639 & 2.67e-07 & 3.17 & 37.139 & 5.43e+01 & 1.27 & 0.756\\
			bfw62a & 4.27e-05 & 1.51 & 3.626 & 1.92e-11 & 4.75 & 0.352 & 1.97e-06 & 2.24 & 2.598 & 4.02e-13 & 5.09 & 0.014\\
			bfw782a & 3.99e-04 & 0.89 & 13.881 & 4.05e-06 & 3.04 & 0.684 & 5.09e-08 & 3.54 & 120.165 & 1.11e+02 & 1.38 & 10.282\\
			bwm200 & 4.78e-09 & -2.41 & 0.074 & 3.95e-08 & -0.26 & 0.028 & 1.78e-09 & -2.09 & 0.505 & 3.68e-10 & -1.94 & 0.352\\
			dwa512 & 2.12e-05 & 0.98 & 0.580 & 1.02e-08 & 2.76 & 0.442 & 1.02e-07 & 2.43 & 370.415 & 5.61e-07 & 1.63 & 53.139\\
			dwb512 & 3.74e-05 & 1.99 & 0.352 & 8.23e-07 & 3.23 & 0.307 & 1.84e-07 & 3.36 & 51.577 & 1.65e-05 & 2.06 & 57.646\\
			lop163 & 5.15e-05 & 1.11 & 2.869 & 2.50e-12 & 4.74 & 0.611 & 1.20e-08 & 2.91 & 0.960 & 1.13e-05 & 1.40 & 2.412\\
			mhd416a & 5.23e-11 & -0.01 & 0.149 & 1.06e-07 & -0.54 & 0.076 & 9.23e-09 & -0.17 & 3.203 & 5.76e+01 & -0.37 & 1.237\\
			mhd416b & 1.01e-05 & 2.54 & 0.347 & 3.14e-08 & 3.87 & 0.284 & 6.63e-07 & 3.23 & 669.339 & 3.08e-06 & 2.75 & 2.933\\
			odep400a & 4.77e-05 & 0.46 & 2.541 & 3.80e-08 & 3.13 & 0.230 & 8.18e-08 & 1.90 & 17.341 & 7.84e-07 & 1.19 & 27.961\\
			olm100 & 4.05e-08 & -1.51 & 0.044 & 2.04e-09 & -0.62 & 0.020 & 8.10e-08 & -0.92 & 0.117 & 1.17e-10 & -1.04 & 0.045\\
			olm500 & 3.21e-08 & -2.80 & 0.321 & 1.62e-06 & -1.91 & 0.091 & 1.71e-05 & -1.78 & 12.095 & 3.12e+03 & -2.96 & 3.099\\
			rbs480a & 9.96e-09 & -2.04 & 0.972 & 6.70e-08 & -0.79 & 0.415 & 1.51e-05 & -0.68 & 6.660 & 7.39e+01 & -0.85 & 2.929\\
			rbs480b & 8.31e-09 & -2.03 & 1.041 & 7.88e-08 & -1.15 & 0.509 & 1.26e-05 & -1.18 & 7.049 & 9.45e+01 & -1.15 & 10.893\\
			rdb200 & 5.54e-07 & -0.20 & 0.176 & 1.11e-08 & 0.32 & 0.152 & 8.26e-08 & 0.17 & 0.275 & 2.94e+01 & 0.16 & 0.168\\
			rdb450 & 1.02e-07 & -0.62 & 0.288 & 1.27e-04 & 0.51 & 0.364 & 5.96e-06 & 0.40 & 2.293 & 6.48e+01 & 0.17 & 1.605\\
			rdb968 & 3.49e-06 & -0.37 & 1.896 & 6.32e-08 & 1.91 & 0.431 & 4.14e-06 & 1.08 & 17.603 & 1.28e+02 & 0.70 & 17.414\\
			rw136 & 2.47e-05 & 1.30 & 4.793 & 1.58e-14 & 5.99 & 0.645 & 4.88e-08 & 2.60 & 4.538 & 1.31e-06 & 1.92 & 1.696\\
			rw496 & 1.83e-04 & 0.89 & 1.213 & 7.88e-09 & 4.59 & 0.533 & 1.06e-08 & 3.03 & 12.119 & 1.56e-10 & 3.75 & 45.088\\
			tols340 & 8.53e-10 & -2.91 & 0.148 & 6.85e-06 & -3.02 & 0.054 & 2.35e-04 & -2.79 & 3.484 & 3.50e+04 & -0.30 & 0.995\\
			tols90 & 6.20e-09 & -2.32 & 0.090 & 1.68e-07 & -1.90 & 0.026 & 1.52e-04 & -2.11 & 0.294 & 2.13e-14 & -0.82 & 0.113\\
			tub100 & 1.62e-10 & -2.57 & 0.061 & 2.76e-08 & -1.21 & 0.018 & 2.44e-09 & -2.47 & 0.081 & 9.34e-10 & -2.39 & 0.046\\
			\midrule
			avg & 5.44e-05 & -0.28 & 3.662 & 7.49e-06 & 1.29 & 0.588 & 1.97e-05 & 0.69 & 57.253 & 1.55e+03 & 0.44 & 11.100\\
			\bottomrule
	\end{tabular}}	
\end{table}

In conclusion, the best optimization solver for all NLP formulations of \eqref{eq:aeicp} among all compared methods is IPOPT. HDCA-NI is competitive with IPOPT and outperforms KNITRO and FILTERSD in terms of the objective value.

\section{Conclusions}
\label{sec:conclusions}

In this paper, we established three DC programming formulations of \eqref{eq:aeicp} based on the DC-SOS decomposition, which are numerically solved via several accelerated DC programming approaches, including BDCA with exact and inexact line search, ADCA, InDCA and two proposed novel hybrid accelerated DCA (HDCA-LI and HDCA-NI). Numerical simulations were carried out to compare the performance of the proposed DCA-type algorithms against the cutting-edge optimization solvers (IPOPT, KNITRO and FILTERSD). Numerical results indicated noteworthy performance, even on large-scale and ill-conditioned datasets.

Future work will prioritize multiple avenues for enhancement. We aim to refine these DCA-type algorithms to improve the quality of computed results in terms of $c$. A particular focus will be on improving the initialization process, possibly through insightful heuristics, to increase efficiency, especially for large-scale and ill-conditioned instances. Moreover, designing a more efficient algorithm to address the QP formulations of convex subproblems required in all DCA-type methods, instead of relying on a general QP solver will also be a key aspect of our upcoming research efforts.


\bibliographystyle{amsplain}
\bibliography{references}

\end{document}